\documentclass[10pt,oneside]{book}

\usepackage[final]{cl_mystyle}

\docstyle

\newcommand{\nin}{\notin}
\newcommand{\To}{\rightarrow}
\newcommand{\gdw}{\leftrightarrow}
\newcommand{\mult}{\times}

\newcommand{\A}{\frak{A}}
\newcommand{\B}{\frak{B}}
\newcommand{\C}{\frak{C}}
\newcommand{\D}{\frak{D}}
\newcommand{\T}{\frak{T}}
\newcommand{\F}{\frak{F}}
\newcommand{\Q}{\frak{Q}}
\newcommand{\G}{\frak{G}}
\newcommand{\Lfr}{\frak{L}}
\newcommand{\AAA}{{\cal A}}
\newcommand{\FF}{{\cal F}}
\newcommand{\CC}{{\cal C}}
\newcommand{\RR}{{\cal R}}

\newcommand{\VV}{{\cal V}}
\newcommand{\MM}{{\cal M}}
\newcommand{\NN}{{\cal N}}
\newcommand{\SSS}{{\cal S}}

\begin{document}

    \beforepreface

    \thispagestyle{empty}

\null\vskip0.5in

\begin{center}
    \hyphenpenalty=10000\Large\uppercase\expandafter{\textbf{Rosenberg's characterization\\ of maximal clones}}
\end{center}

\vfill

\begin{center}
    \large\rm By\\
    Michael Pinsker\\
    \vskip0.3cm
    \small{marula@gmx.at}
\end{center}

\vfill

\begin{center}
    \footnotesize
    DIPLOMA THESIS \\
    AT THE\\
    \uppercase\expandafter{Vienna university of technology} \\
    \uppercase\expandafter{May 2002}
\end{center}

\vskip0.75in

\newpage

    \prefacesection{Abstract}

We will give a proof of I. G. Rosenberg's characterization of
maximal clones, first published in \cite{Ros70}. The theorem lists
six types of relations on a finite set such that a clone over this
set is maximal if and only if it contains just the functions
preserving one of the relations of the list. In Universal Algebra,
 this translates immediately into a characterization of the
finite preprimal algebras: A finite algebra is preprimal if and
only if its term operations are exactly the functions preserving a
relation of one of the six types listed in the theorem. The
difficult part of the proof is to show that all maximal clones or
preprimal algebras respectively are of that form. This follows
from, and, as we will also demonstrate, is indeed equivalent to, a
characterization of primal algebras: We will show that the primal
algebras are exactly those whose term operations do not preserve
any of the relations on the list.

    \prefacesection{Preface}

A \emph{clone} (\emph{cl}osed \emph{o}peration \emph{ne}twork)
$\CC$ over a set $A$ is a set of operations on this set which
contains the projections and which is closed under compositions.
The set of all clones over $A$ forms a lattice $Clone(A)$ with
respect to inclusion, and a clone is called \emph{maximal} if and
only if it is a dual atom in $Clone(A)$.

It is a fact that if $A$ is finite, then every clone is contained
in a maximal clone and the maximal clones are finite in number. In
his work \cite{Ros70} I. G. Rosenberg gave a characterization of
the maximal clones over a a finite base set in terms of relations:
The theorem lists six types of relations on $A$ such that a clone
is maximal if and only if it is just the set of functions
\emph{preserving} one of the relations of the list.

However, the original proof of this deep theorem is quite
technical and hard to follow. It is the aim of the present work to
provide a shorter and somewhat more understandable proof.

Our proof is based on the one by R. W. Quackenbush in
\cite{Qua71}, who showed the more difficult implication of the
theorem, namely that every maximal clone is of the form described
before. It draws heavily on results of R. W. Quackenbush
\cite{Qua80} on algebras with minimal spectrum, of H. P. Gumm
\cite{Gum79} on algebras in permutable varieties, and of A. Foster
and A. Pixley \cite{FP64} on primality. Also a part of the
original proof of I. G. Rosenberg has been included. We would like
to add that there exists another new proof of the difficult
implication of the theorem by V. A. Buevich in \cite{Bue96}.

This thesis has been divided into three chapters. In the first
chapter, we introduce the theorem and explain the connection
between maximal clones and preprimal algebras. Chapter 2 contains
the proof of half of the equivalence: Every maximal clone is a set
of functions preserving one of the relations listed in the
theorem. Chapter 3 is devoted to the proof of the converse
statement that all relations of the list yield a maximal clone.

All global conventions regarding notation will be made in the
first chapter together with the basic definitions, and additional
conventions will be introduced in Notations \ref{not-ch2},
\ref{not-qua} and \ref{not-ch3}. We tried to keep this work
self-contained, the reader is assumed to be familiar only with the
rudiments of Universal Algebra, lattice theory, and some basic
facts about groups and fields; information on clones can be found
in \cite{Sze86}.

I would like to thank M. Goldstern for his support and many
helpful suggestions.

    \tableofcontents

    \afterpreface

    \chapter{Rosenberg's preprimal algebra characterization}
We will state the characterization of the maximal clones and
provide the reader with the necessary definitions. Moreover, the
connection between another possible viewpoint of the theorem,
namely the characterization of finite preprimal algebras, and the
theorem itself as a statement about clones will be explained.

\begin{defn}
    Let $A$ be a set and denote by $\FF_{n}$ the set of all $n-$ary functions
    on $A$. Then $\FF=\bigcup_{n=0}^{\infty}\FF_{n}$ is the set of
    all functions on $A$ of arbitrary arity. A \textit{clone} is a
    subset of $\FF$ which is closed under compositions and which
    contains all projections. The set of all clones on $A$ form a
    lattice $Clone(A)$ with respect to inclusion. A clone is called
    \textit{maximal} iff it is maximal in $Clone(A)\setminus\{\FF\}$.
\end{defn}

In order to bring these definitions into the context of Universal
Algebra, one can think of a clone $\CC$ on $A$ as the set of term
operations of the algebra $\A=(A,\CC)$. Conversely, given an
algebra $\A=(A,F)$, the term operations $\T(F)$ form a clone over
$A$. This interpretation of clones makes sense, for it provides
the possibility of making use of the existing apparatus of
Universal Algebra, e.g. congruence relations. It is for this
reason that we will talk about algebras rather than about clones
for the biggest part of our proof.

\begin{defn}
    An algebra $\A$ is \textit{primal} iff every function on $A$ is a
    term operation of $\A$; $\A$ is \textit{preprimal} iff it is not
    primal but for any function $f$ not a term operation of $\A$,
    $(A,F \cup \{f\})$ is primal.
\end{defn}
By the previous discussion, maximal clones correspond to preprimal
algebras and vice-versa. Let $\RR_{n}$ be the set of all $n-$ary
relations on $A$; then $\RR=\bigcup_{n=1}^{\infty}\RR_{n}$ is the
set of all relations on $A$ of arbitrary finite arity. We define
for an arbitrary set $R\subseteq\RR$ of relations on $A$ the set
of polymorphisms $Pol(R)$, that is, if we write $R_k$ for the
$k$-ary relations in $R$ and $a_1,...,a_n$ for the coordinates of
an $n$-tuple $a$,
\begin{eqnarray*}
\begin{split}
    Pol(R) = \bigcup_{n=0}^\infty\{&f \in \FF_n : \forall k\geq 0 \,\forall \rho \in
    R_k \,\forall r_1,...,r_n \in \rho
    \\& ((f(r_{11},...,r_{n1}),...,f(r_{1k},...,r_{nk})) \in
    \rho)\}.
\end{split}
\end{eqnarray*}
With this definition, Rosenberg's theorem states that a clone over
a finite set $A$ is maximal iff it is of the form $Pol(\{\rho
\})$, where $\rho$ is a relation in one of six classes to be
specified later. To formulate Rosenberg's theorem in detail, we
need a couple of definitions.

For a function $f$ on $A$ define the \textit{graph} of $f$ to be
the set $\{ (a,f(a)) : a \in A\}$. Sometimes we will talk about a
function and mean the graph of the function as a subset of $A^2$;
confusion is unlikely since things should be clear from
context.\newline A permutation $\pi$ is $prime$ iff all cycles of
$\pi$ have the same prime length.\newline We call a subset $\rho
\subseteq A^{4}$ \textit{affine} iff there is a binary operation
$+$ on $A$ such that $(A,+)$ is an abelian group and $(a,b,c,d)\in
\rho \leftrightarrow a+ b = c + d$ holds. An affine $\rho$ is
\textit{prime} iff $(A,+)$ is an abelian $p$-group for some prime
$p$, that is, all elements of the group have the same prime order
$p$. \newline For $h\geq1$ a subset $\rho\subseteq A^{h}$ is
\textit{totally symmetric} iff for all permutations $\pi$ of
$\{1,...,h\}$ and all tuples $(a_1,...,a_h)\in A^h$,
$(a_1,...,a_h)\in\rho$ iff $(a_{\pi(1)},...,a_{\pi(n)})\in\rho$.
Define $\iota_h^A \subseteq A^h$ by
\[
\iota_h^A=\{\,(a_1,...,a_h)\, |\, \exists i \, \exists j
\,(\,i\neq j \,\wedge \,a_i=a_j\,)\,\}.
\]
Then $\rho$ is called \textit{totally reflexive} iff $\iota_h^A
\subseteq \rho$. Note that for $h=2$, totally reflexive means
reflexive and totally symmetric means symmetric. If $\rho$ is
totally reflexive and totally symmetric we define the
\textit{center} of $\rho$ to be the set
\[
C(\rho)=\{a \in A | \forall \,a_2,...,a_h \in A \,
(a,a_2,...,a_h)\in\rho\}.\] We say that $\rho\subseteq A^h$ is
\emph{central} iff it is totally reflexive, totally symmetric and
has a nonvoid center which is a proper subset of $A$. Note that
$h\leq |A|$ as otherwise we would have $\rho\supseteq
\iota_h^A=A^h$ and the center of $\rho$ would be trivial.
\newline For an arbitrary set $S$ and $1\leq r\leq \lambda$, denote the
$r-$th projection from $S^\lambda$ onto $S$ by $\pi_r^\lambda$.
Now let $h=\{0,1,...,h-1\}$ and define $\omega_\lambda$ to be the
$h$-ary relation on $h^\lambda$ satisfying $(a_1,...,a_h) \in
\omega_\lambda$ iff for all $1\leq r\leq \lambda$,
$(\pi_r^\lambda(a_1),...,\pi_r^\lambda(a_h))\in \iota_h^h$. For
$3\leq h\leq |A|$, we call a $h$-ary relation $\rho$ on $A$
$h$\emph{-regularly generated} iff there exists a $\lambda \geq 1$
and a surjection $\varphi: A \rightarrow h^\lambda$ such that
$\rho=\varphi^{-1}(\omega_\lambda)$. Note that for any relation,
$h$-regularly generated implies totally reflexive and totally
symmetric.

Now here comes the theorem.
\newpage
\begin{thm}[I. G. Rosenberg \cite{Ros70}]
    Let $1<|A|<\aleph_0$. A clone $\CC$ on $A$ is maximal if and only if it
    is of the form $Pol(\rho)$, where $\rho$ is an
    $h$-ary relation belonging to one of the following classes:
    \begin{enumerate}
        \item{The set of all partial orders with least and greatest element}
        \item{The set of all prime permutations}
        \item{The set of all non-trivial equivalence relations}
        \item{The set of all prime-affine relations}
        \item{The set of all central relations}
        \item{The set of all h-regularly generated relations}
    \end{enumerate}
\end{thm}

We will refer to the six classes as Rosenberg's list ($RBL$) from
now on. Then in the terminology of algebras, the theorem sounds
like this.

\begin{cor}
    A finite non-trivial algebra $\A$ is preprimal iff there
    exists a relation $\rho$ in $RBL$ such that $\T(\A)=Pol(\rho)$.
\end{cor}

\begin{rem}
    As with Rosenberg's theorem the maximal clones over a set $A$ with finite cardinality $\kappa$ are known, one can
    calculate their number $\eta_{\kappa}$. That number grows fast with the size $\kappa$. Here are values for a couple of cardinalities $\kappa$:
    \begin{center}
        \begin{tabular}{c|c|c|c|c|c|c}
           $\kappa$&2&3&4&5&6&7\\
          \hline
          $\eta_\kappa$&$5$&$18$&82&643&15182&7848984
        \end{tabular}
    \end{center}
\end{rem}

\begin{rem}
    The clone lattice $Clone(A)$ is countable only for $|A|=2$.
    For $|A|\geq 3$ we have already $|Clone(A)|=2^{\aleph_0}$.
\end{rem}

\chapter{Primal algebra characterization}

We will prove the more difficult part of the equivalence by
proving the following theorem.
\begin{thm}\label{Ros1}
    If a finite non-trivial algebra $\A$ has no subalgebra of a finite
    power of $\A$ belonging to $RBL$, then $\A$ is primal.
\end{thm}
The required implication in Rosenberg's theorem follows indeed.
\begin{cor}
    If a finite non-trivial algebra $\A$ is preprimal then the set of
    term operations of $\A$ is of the form $Pol(\rho)$, where $\rho$
    is a relation in $RBL$.
\end{cor}
\begin{proof}
    Since $\A$ is not primal, by the last theorem there exists a
    subalgebra of a finite power of $\A$ with universe $\rho$ in
    $RBL$; hence, the term operations satisfy $\T(\A)\subseteq
    Pol(\rho)$. But as $Pol(\rho)$ is closed under composition and
    projections and as $\A$ is preprimal, $\T(\A)= Pol(\rho)$.
\end{proof}
The corollary is in fact equivalent to the theorem.
\begin{thm}
    If all finite non-trivial preprimal algebras $\A$ satisfy $\T(\A)= Pol(\rho)$, where $\rho$
    is a relation in $RBL$, then every finite non-trivial algebra which preserves no relation belonging to $RBL$
    is primal.
\end{thm}
\begin{proof}
    Let $\A$ be a finite non-trivial algebra preserving no
    relation belonging to $RBL$. Then the clone $\T (\A)$ is
    contained in no clone of the form $Pol(\rho)$, $\rho\in RBL$.
    But since all maximal clones are of that form and since
    $Clone(A)$ is dually atomic (see \cite{Sze86}), this means that $\T (\A)$ must be
    the greatest element in that lattice and thus the clone of
    all functions on $A$. Hence, $\A$ is primal.
\end{proof}

To prove Theorem \ref{Ros1}, we will first show that the
hypotheses imply that all subalgebras of finite powers of $\A$
have cardinality a power of the cardinality of $\A$, which is a
result by R. Quackenbush in \cite{Qua71}. R. W. Quackenbush also
essentially showed in \cite{Qua80} that then the algebra generates
a congruence permutable variety; we will follow his proof in the
beginning, but then use a slightly different approach to prove
this, combining works of D. Clark and P. Krauss in \cite{CK76} and
of I. Chajda and G. Eigenthaler in \cite{CE01}. Following H. P.
Gumm in \cite{Gum79} and then H. Werner in \cite{Wer74} we will
conclude that all powers of $\A$ can only have factor congruences,
which trivially implies that the equational class generated by
$\A$ is congruence distributive. A criterion for primality due to
A. Foster and A. Pixley \cite{FP64} will finally conclude the
proof. Here is a summary of which implications we will prove; it
might be helpful to look at it from time to time. The notions
which occur in those implications will be defined in the
respective sections.
\begin{itemize}
    \item{If $\A$ is a finite non-trivial algebra having no subalgebra of a
            power of $\A$ belonging to $RBL$, then $\A$ has almost minimal spectrum (Theorem \ref{q-main-3.4}).}
    \item{If $\A$ is a finite non-trivial algebra with almost minimal spectrum, then the
            variety generated by $\A$ is congruence
            permutable (Theorem \ref{q-qua-1.1}).}
    \item{If $\A$ is a finite simple algebra in a permutable variety, then $\A$
            is either prime affine or its powers have only (trivial) factor congruences (Theorem \ref{cl-3}).}
    \item{If $\A$ is a finite simple non-trivial algebra with no proper subalgebras
            and no non-trivial automorphisms, and if $\A$ generates a permutable and distributive
            variety, then $\A$ is primal (Theorem \ref{q-fp3.1}).}
\end{itemize}

\begin{nota}\label{not-ch2}
    Until the end of the chapter, as we will be proving Theorem \ref{Ros1}, we will denote the algebra
    satisfying the hypotheses of the theorem by $\A=(A,F)$. We
    will use the symbol $F$ also for the corresponding operations
    on powers of $\A$. The congruence lattice of $\A$ will play an
    important role and we will write $Con(\A)$ for it. By $0\in
    Con(\A)$ we mean the diagonal $\{(a,a)|\,a\in A\}$ and
    by $1\in Con(\A)$ the trivial congruence $A^2$.
\end{nota}

    \section{$\A$ has almost minimal spectrum}

\begin{defn}
    The \textit{spectrum} $Spec(\VV)$ of a variety $\VV$ is the set of
    all cardinalities of finite members of $\VV$. For a finite algebra
    $\A$ we define $Spec(\A)=Spec(\VV(\A))$, where $\VV(\A)$ denotes the
    variety determined by $\A$. $\A$ is said to have
    $\textit{minimal spectrum}$ iff $Spec(\A)=\{\,|\A|^n\, |\, n \geq
    0 \,\}$.
\end{defn}

The original goal of the author was to prove in this section that
our algebra $\A$ has minimal spectrum. This would have made it
easy to find a title for this section. However, it did not work
out and we will obtain that result later. The following definition
will help us out for the moment.

\begin{defn}
    We say that a finite algebra $\A$ has \emph{almost minimal
    spectrum} iff all subalgebras of finite powers of $\A$ have
    cardinality a power of the cardinality of $\A$.
\end{defn}
\begin{rem}
    Recall that every algebra in $\VV(\A)$ is a homomorphic image
    of a subalgebra of a power of $\A$. The notion of almost
    minimal spectrum is thus weaker than the one of minimal
    spectrum.
\end{rem}
This section is devoted to the proof of the following theorem
which is due to R. W. Quackenbush \cite{Qua71}.
\begin{thm}\label{q-main-3.4}
    Let $\A$ be a finite non-trivial algebra having no subalgebra of a
    power of $\A$ belonging to $RBL$. Then $\A$ has almost minimal
    spectrum.
\end{thm}

The proof will be by contradiction: Suppose $\A$ does not have
almost minimal spectrum;  then there is an $m$ and a subalgebra
$\B$ of $\A^m$ with $|B|$ not a power of $\kappa=|A|$. Choose $m$
minimal in the sense that for all $n<m$ every subalgebra of $\A^n$
has cardinality a power of $\kappa$. As $\A$ has no proper
subalgebras (a proper subalgebra would be a unary central
relation), $\B$ must even be a subdirect product (that is, the
projection of $\B$ on any coordinate is onto); thus clearly,
$m>1$.
\begin{nota}\label{not-qua}
    For the rest of this section (that is, until Theorem
    \ref{q-main-3.4} has been proven), we will extend Notation \ref{not-ch2} and
    use the following conventions: $\A$ will be assumed to satisfy
    all hypotheses of Theorem \ref{q-main-3.4}. The letter $\kappa$ will be reserved for
    the cardinality of $\A$. For the universe of $\A$ we write
    $A=\{\alpha_1,...,\alpha_\kappa\}$. $\B$ and $m$ as just
    defined will not change their meaning.
\end{nota}
\noindent For $E=\{i_1,...,i_j\}\subseteq \{ 1,...,m\}$ where
$i_1<...<i_j$ define the projection
$$
    \pi_E:\quad
    \begin{matrix}
        A^m && \rightarrow && A^j \\
        (a_1,...,a_m) && \mapsto && (a_{i_1},...,a_{i_j}).
    \end{matrix}
$$
Define further for $1\leq i\leq m$ the projection
$\rho_i=\pi_{E(i)}$ where $E(i)=\{1,...,i-1,i+1,...,m\}$.

\begin{lem}\label{q-main-4.1}
    $|\rho_i(B)|=\kappa^{m-1}$ for $1\leq i\leq m$.
\end{lem}

\begin{proof}
    We will proof by induction on $n$ that any
    projection $\pi$ mapping from $B$ to $A^n$, where $n< m$, is onto.
    Since $\A$ has no proper subalgebras, all the projections
    $\pi_i^m$ from $B$ to $A$ are onto and hence our assertion is
    true for $n=1$. Now assume that for all $i<n<m$, if $|E|=i$ then
    $|\pi_E(B)|=|A|^i$ and let $E$ be an $n$-element subset of
    $\{1,...,m\}$. As $\pi_E(B)$ is a subalgebra of $A^n$ and
    $n<m$, there exists an $i$ such that $|\pi_E(B)|=|A|^i$.
    Trivially, $|\pi_E(B)|\leq|A|^n$ and by induction hypothesis,
    $|\pi_E(B)|\geq|A|^{n-1}$ so that either $|\pi_E(B)|=|A|^{n-1}$ or
    $|\pi_E(B)|=|A|^{n}$. Consider in the first case
    $E'=E\setminus \{i\}$ for an arbitrary $i\in E$. By induction
    hypothesis we know that also $|\pi_{E'}(B)|=|A|^{n-1}$. From
    this follows that for $b, b' \in B$, if $b_j=b_j'$ for $j\in
    E'$ then $b_i=b_i'$. Hence for any $b,b'\in B$,
    $\rho_i(b)=\rho_i(b')$ implies that $b=b'$ and $\rho_i$ is
    one-one. But this means that $\rho_i$ embeds $\B$ as a
    subalgebra of $\A^{m-1}$, contradicting our assumption that
    $m$ is minimal with respect to having a subalgebra of
    cardinality not a power of $A$. Thus $|\pi_E(B)|$ must be equal to
    $|A|^n$ and the induction is complete.
\end{proof}

\begin{cor}\label{q-main-4.2}
    $|A|^{m-1}<|B|<|A|^m$.
\end{cor}

Let $P$ be a partition $\{1,...,m\}$, and let $\sim_P$ be the
equivalence relation induced by $P$. Define a subset $B_P$ of $B$
by $B_P=\{b\in B \,|\, i\sim_P j \rightarrow b_i=b_j\}$. Then
clearly, $B_P$ is a subuniverse of $B$. To denote a partition, we
will only list its non-trivial classes; $(i,j)$ denotes the
partition with only one non-trivial class, $\{i,j\}$.

\begin{lem}\label{q-main-4.3}
    Let $m\geq 4$. If for some $1\leq i,j\leq m,\,i\neq j$ we have
    $|B_{(i,j)}|=|A|^{m-1}$, then the same holds for all $1\leq i,j\leq
    m,\,i\neq j$.
\end{lem}

\begin{proof}
    Let $|B_{(i,j)}|=|A|^{m-1}$; it suffices to show that for
    $k\neq i,j$ we have $|B_{(i,k)}|=|A|^{m-1}$. Our assumption
    $|B_{(i,j)}|=|A|^{m-1}$ obviously implies $|B_{(i,j,k)}|=|A|^{m-2}$. Since
    $B_{(i,k)}$ can be embedded into $\A^{m-1}$ by leaving away
    the $k$-th coordinate, $|B_{(i,k)}|$ must be a power of $|A|$.
    Trivially, $|B_{(i,k)}|\leq |A|^{m-1}$ and since
    $|B_{(i,j,k)}|\leq |B_{(i,k)}|$ we have $|B_{(i,k)}|\geq
    |A|^{m-2}$. Suppose now that $|B_{(i,k)}|=
    |A|^{m-2}$; then $|B_{(i,j,k)}|=|B_{(i,k)}|$ and so $b_i=b_k$
    implies $b_i=b_j$ for all $b\in B$. But this contradicts that by the proof of Lemma
    \ref{q-main-4.1}, $|\pi_{\{i,j,k\}}(B)|=|A|^3$. Therefore,
    $|B_{(i,k)}|\neq |A|^{m-2}$ and so $|B_{(i,k)}|= |A|^{m-1}$.
\end{proof}

Define a subset $B'$ of $A^m$ by
$$B'=\{(a_2,a_3,...,a_m,a_m')\,|\,\exists a_1 \in A
((a_1,...,a_m)\in B \wedge (a_1,...,a_{m-1},a_m')\in B)\}.$$ Then
$B'$ is a subuniverse of $\A^m$ and the following holds:

\begin{lem}\label{q-main-4.4}
    Let $m\geq 4$. If $|B_{(2,3)}|=|A|^{m-2}$, then
    \begin{itemize}
    \item{$|A|^{m-1}<|B'|<|A|^m$}
    \item{$|B_{(1,2)}'|=|A|^{m-2}$}
    \item{$|B_{(m-1,m)}'|=|A|^{m-1}$}
    \end{itemize}
\end{lem}

\begin{proof}
    First note that $|B_{(m-1,m)}'|=|\rho_1(B)|=|A|^{m-1}$, the
    latter equality provided by Lemma \ref{q-main-4.1}.
    Furthermore, $|B'|\geq |B_{(m-1,m)}'|=|A|^{m-1}$. Since we
    know that $|\rho_m(B)|=|A|^{m-1}$ but $|B|>|A|^{m-1}$, there
    exist $a_1,...,a_{m-1},a_{m},a_m' \in A$ such that $a_m\neq
    a_m'$ and both $(a_1,...,a_{m-1},a_m)\in B$ and $(a_1,...,a_{m-1},a_m')\in
    B$. Hence, $(a_2,...,a_{m-1},a_m,a_m')\in B'$ and so
    $|B'|>|A|^{m-1}$. Given $a_1,...,a_{m-1}\in A$, $|\rho_m(B)|=|A|^{m-1}$
    implies there exists an $a_m\in A$ such that $(a_1,...,a_m)\in
    B$. Since we assume $|B_{(2,3)}|=|A|^{m-2}$, from $a_2=a_3$ it
    follows that such an $a_m$ is unique so that if we choose any
    $a_m' \neq a_m$, $(a_3,a_3,a_4,...,a_{m-1},a_m,a_m')\notin B'$.
    Thus, $|B'|<|A|^m$, and $|B_{(1,2)}'|=|B_{(2,3)}|=|A|^{m-2}$.
\end{proof}

\begin{lem}\label{q-main-4.5}
    Let $m\geq 4$. Then for $1\leq i < j \leq m$,
    $|B_{(i,j)}|=|A|^{m-1}$.
\end{lem}

\begin{proof}
    If $|B_{(2,3)}|=|A|^{m-1}$ then the assertion follows from Lemma \ref{q-main-4.3}.
    If not, then $|B_{(2,3)}|=|A|^{m-2}$, and so by the last
    lemma $|B'_{(m-1,m)}|=|A|^{m-1}$. Thus, $B'$ satisfies
    the hypotheses on $B$ in Lemmas \ref{q-main-4.1} and
    \ref{q-main-4.3} and application of Lemma
    \ref{q-main-4.3} yields $|B'_{(1,2)}|=|A|^{m-1}$
    contradicting $|B'_{(1,2)}|=|A|^{m-2}$ which we established in
    the previous lemma. Hence, $|B_{(2,3)}|=|A|^{m-2}$ is
    impossible and the lemma follows.
\end{proof}

We can summarize what we have established so far:

\begin{thm}\label{q-main-4.6}
    Let $\B$ be a subalgebra of $\A^m$ with $|B|$ not a power of $\kappa$. If $m\geq 4$, then
    $B$ is totally reflexive.
\end{thm}

Denote by $\B^*=(B^*,F)$ the subalgebra of $\A^m$ generated by
$\iota_m^A$. By considering this algebra we will show that for
$m\geq 4$ we can assume without loss of generality that $\B$ is
totally reflexive and totally symmetric:

\begin{thm}\label{q-main-4.7}
    Let $m \geq 4$. Then $B^*$ is totally reflexive and totally
    symmetric and $|A|^{m-1}<|B^*|<|A|^m$.
\end{thm}

\begin{proof}
    $\iota_m^A$ is both totally reflexive and totally symmetric
    and it is easy to see that $B^*$ inherits those properties. By
    Theorem \ref{q-main-4.6} we have $\iota_m^A \subseteq B$ and so $B^*
    \subseteq B$; hence, $|B^*|\leq |B| < |A|^m$. Moreover,
    $|B_{(1,2)}|=|A|^{m-1}$ by Lemma \ref{q-main-4.5}, and
    $B_{(1,2)}$ is obviously a proper subset of $\iota_m^A$. Thus,
    $|A|^{m-1}<|\iota_m^A|<|B^*|$ and the theorem follows.
\end{proof}

We have shown that in the case $m \geq 4$, we can assume $B$ to be
totally reflexive and totally symmetric by replacing $B$ with
$B^*$ if necessary. Our next step will be to prove the totally
reflexive and totally symmetric possibility absurd; as a result,
$m\geq 4$ cannot occur.

\subsubsection{The totally reflexive and totally symmetric case}

First note that in this case $m \leq \kappa$ since otherwise every
element of $A^m$ would have two equal components and so the total
reflexivity of $B$ would imply $B=A^m$. Now choose $h\leq \kappa$
to be maximal with respect to $\A^h$ containing a proper totally
reflexive and totally symmetric subalgebra; let $\C=(C,F)$ be a
maximal subalgebra of $\A^h$ of that kind.

For $h \leq n \leq \kappa$ define sets $C_n \subseteq A^n$ to
contain all $(a_1,...,a_n)$ for which there exists an $a\in A$
such that for each $(h-1)$-element subset $\{i_1,...,i_{h-1}\}$ of
$\{1,...,n\}$, $(a_{i_1},...,a_{i_{h-1}},a)\in C$. Then $(C_n,F)$
is a subalgebra of $\A^n$ and is totally symmetric as $C$ is.

\begin{lem}\label{q-main-7.1}
    Either $C=C_h$ or $C_h=A^h$.
\end{lem}

\begin{proof}
    Let $(a_1,...,a_h)\in C$. Set $a=a_1$; then by the total
    reflexivity and total symmetry of $C$ we have that for $1\leq
    i\leq h$, $(a_1,...,a_{i-1},a_{i+1},...,a_h,a_1)\in C$ so that
    $(a_1,...,a_h)\in C_h$. Hence, $C\subseteq C_h$ and so by
    the maximality of $C$, $C=C_h$ or $C_h=A^h$.
\end{proof}

\begin{lem}\label{q-main-7.2}
    If $C_h=A^h$, then $C_\kappa=A^\kappa$.
\end{lem}

\begin{proof}
    All $C_n$ are totally symmetric, $h\leq n\leq \kappa$. Thus, by the maximality of
    $h$, if $C_n$ is also totally reflexive
    then $C_n=A^n$. But clearly the definition of $C_n$ implies
    that if $C_n=A^n$, then $C_{n+1}$ is totally reflexive so that
    by induction we get $C_\kappa=A^\kappa$.
\end{proof}

The following lemma states that the case $C_h=A^h$ is impossible.

\begin{lem}
    If $C\neq C_h$, then $C$ is central.
\end{lem}

\begin{proof}
    By Lemmas \ref{q-main-7.1} and \ref{q-main-7.2} our
    hypothesis implies that $C_\kappa=A^\kappa$. Hence,
    $(\alpha_1,...,\alpha_\kappa)\in C_\kappa$. Therefore, by definition of $C_\kappa$, there is an $a\in A$ such that
    for every $h-1$-element subset $\{\alpha_{i_1},...,\alpha_{i_{h-1}}\}$ of
    $\{\alpha_1,...,\alpha_\kappa\}=A$,
    $(\alpha_{i_1},...,\alpha_{i_{h-1}},a)\in C$. Together with the fact that $C$ is
    totally reflexive and totally symmetric this implies that $a$
    is in the center of $C$. But since $C\neq A^h$ this means that
    $C$ is central.
\end{proof}

We may therefore assume that $C=C_h$. Such a $C$ is called
\emph{homogeneous}. Note that in the case $h=2$, if $(a_1,a_3)\in
C$ and $(a_2,a_3)\in C$, we have that also $(a_1,a_2)\in C$ (set
$a=a_3$ in the definition of $C_h$). As $C$ is symmetric this
means that $C$ is transitive and so, together with its
reflexivity, we get that $C$ is a non-trivial equivalence relation
and thus in Rosenberg's list. Therefore, we may assume $h\geq 3$.

In the following, we will make use of the homogeneity of $C$.
Notice therefore that by the definition of $C_h$, to prove that a
tuple $(a_1,...,a_h)$ it is an element of $C$ it suffices to find
an arbitrary $b\in A$ such that for all $1\leq i\leq h$, if we
replace $a_i$ by $b$, then the resulting tuple is in $C$. Such an
element will be referred to as a \emph{replacement element}. The
condition is not only sufficient but also necessary for membership
of $C$.

If $C$ contains all tuples $(a_1,...,a_h)\in A^h$ for which there
exists $(v_1,...,v_h)\in C$ such that
$(a_1,...,a_{i-1},a_{i+1},...,a_h,v_j)\in C$ for all $1\leq
i,j\leq h$ with $i\neq j$, then it is called \emph{strongly
homogeneous}. The tuple $(v_1,...,v_h)$ will be referred to as the
\emph{replacement tuple}. Notice that strongly homogeneous
immediately implies homogeneous if one considers the replacement
tuple containing a replacement element at every coordinate. We
will see that the maximality of $C$ implies that it is strongly
homogeneous; from that we will derive that $C$ is $h$-regularly
generated so that it belongs to Rosenberg's list.

For $h\leq r \leq \kappa$, define $C_r\subseteq A^r$ by
$$
    C_r=\{(a_1,...,a_r) | \forall E \subseteq \{1,...,r\} \,(|E|=h\rightarrow \pi_E(a_1,...,a_r)\in
    C)\}.
$$
Clearly for all $r$, $(C_r,F)$ is a subalgebra of $\A^r$,
$C_r$ is totally symmetric since $C$ is, and $C_h=C$. For $h\leq
 r \leq \kappa$, define $D_r\subseteq A^r$ by
\newline
\begin{eqnarray*}
    \begin{split}
    D_r=\{&(a_1,...,a_r)|\, \exists (b_1,...,b_r)\in C_r\\
          &\forall \,1\leq j\leq r \, \forall \,\{i_1,...,i_{h-2}\}\subseteq\{1,...,r\} \quad(a_{i_1},...,a_{i_{h-2}},a_j,b_j)\in C\}.\\
    \end{split}
\end{eqnarray*}
\newline
Then for all $r$, $(D_r,F)$ is a subalgebra of $\A^r$;
furthermore, $D_r$ is totally symmetric by its symmetric
definition and the total symmetry of $C_r$.

\begin{lem}
    If $C=D_h$ then $C$ is strongly homogeneous.
\end{lem}

\begin{proof}
    Suppose $C=D_h$ and let $(a_1,...,a_h) \in A^h$ and
    $(v_1,...,v_h)\in C$ be given such that
    $(a_1,...,a_{i-1},a_{i+1},...,a_h,v_j)\in C$ for all $1\leq i,j
    \leq h$ with $i\neq j$. Then for
    every $1\leq j\leq h$ and every
    $\{i_1,...,i_{h-2}\}\subseteq\{1,...,h\}$  we have
    that $(a_{i_1},...,a_{i_{h-2}},a_j,v_j)\in C$. This is clear if
    $|\{i_1,...,i_{h-2},j\}|<h-1$ from the total reflexivity of $C$ and
    if not, then there is an $i\neq j$ such that
    $(a_{i_1},...,a_{i_{h-2}},a_j,v_j)\in C$ is  by the total symmetry
    of $C$ equivalent to $(a_1,...,a_{i-1},a_{i+1},...,a_h,v_j)\in C$.
    Since the latter statement is true, setting $(b_1,...,b_h)=(v_1,...,v_h)$ in
    the definition of $D_h$ shows $(a_1,...,a_h)\in D_h=C$. Hence, $C$ is indeed
    strongly homogeneous.
\end{proof}

\begin{lem}\label{q-main-7.3}
    Either $C=D_h$ or $D_h=A^h$.
\end{lem}

\begin{proof}
    Let $(a_1,...,a_h)\in C$ and set in the definition of $D_h$
    $(b_1,...,b_h)$ equal to $(a_1,...,a_h)$. Then $(b_1,...,b_h)\in
    C_h=C$, and so $(a_1,...,a_h)\in D_h$ as $C$ is totally
    reflexive. Therefore, $C\subseteq D_h$ and consequently the maximality
    of $h$ implies $C=D_h$ or $D_h=A^h$.
\end{proof}

\begin{lem}\label{q-main-7.4}
    If $D_h=A^h$, then $D_\kappa=A^\kappa$.
\end{lem}

\begin{proof}
    The proof will be by induction. Suppose that $D_n=A^n$; choose
    an arbitrary $(a_1,...,a_n)\in A^n$ and say it is in $D_n$ via $(b_1,...,b_n)\in
    C_n$. But if
    $(b_1,...,b_n)\in C_n$, then obviously $(b_1,b_1,...,b_n)\in
    C_{n+1}$ and it is easily seen that $(a_1,a_1,...,a_n)\in D_{n+1}$ via
    $(b_1,b_1,...,b_n)$. Since we could have chosen any other
    coordinate instead of the first in that argument, we conclude
    that $D_{n+1}$ is totally reflexive. But $D_{n+1}$ is also
    totally symmetric, and so the maximality of $C$ implies
    $D_{n+1}=A^{n+1}$.
\end{proof}

\begin{thm}\label{q-main-7.5}
    $C$ is strongly homogeneous.
\end{thm}

\begin{proof}
    We will show that $D_\kappa\neq A^\kappa$. Then by Lemma
    \ref{q-main-7.4}, $D_h\neq A^h$, and so by Lemma
    \ref{q-main-7.3}, $C=D_h$ which we know implies that $C$ is
    strongly homogeneous. Suppose
    towards contradiction that $D_\kappa=A^\kappa$: then any vector $(a_1,...,a_\kappa)\in A^\kappa$ is an element
    of
    $D_\kappa$, say via $(b_1,...,b_\kappa)\in C_\kappa$. By definition of $D_\kappa$,
    for every $1\leq j\leq \kappa$ and every
    $\{i_1,...,i_{h-2}\}\subseteq\{1,...,\kappa\}$ we have
    $(a_{i_1},...,a_{i_{h-2}},a_j,b_j)\in C$. We will
    prove by induction that for $0\leq n\leq h$,
    $(a_1,...,a_n,b_{n+1},...,b_h)\in C$. Since $(b_1,...,b_\kappa)\in
    C_\kappa$, $(b_1,...,b_h)\in C_h=C$ and so in the case $n=0$ our
    assertion is true. Suppose it is true for $n<h$ and consider
    $(a_1,...,a_{n+1},b_{n+2},...,b_h)$. If we replace any element
    other that $a_{n+1}$ by $b_{n+1}$, the resulting $h$-tuple
    contains both $a_{n+1}$ and $b_{n+1}$ so that it is in $C$ by the
    preceding discussion. On the other hand, replacing $a_{n+1}$
    by $b_{n+1}$ gives us an element of $C$ by induction
    hypothesis. Hence, the homogeneity of $C$ implies that
    $(a_1,...,a_{n+1},b_{n+2},...,b_h)\in C$ and the induction is
    complete. Now setting $n=h$ yields $(a_1,...,a_h)\in C$. But
    the vector $(a_1,...,a_\kappa)$ was arbitrarily chosen; hence,
    $C=A^h$, contradicting our assumption on $C$.
\end{proof}

\subsubsection{The case of the $h$-regularly generated relations}

We will show that strongly homogeneous $C$ is $h$-regularly
generated. That is, we will find a surjection $\varphi:
A\rightarrow h^m$ such that $C=\varphi^{-1}(\omega_m)$ as defined
in the first chapter. Our first step is to find the equivalence
relation induced by $\varphi$. Define $E\subseteq A^2$ by
$$E=\{(a,b)|\forall (a_1,...,a_{h-2})\in A^{h-2}\,(a_1,...,a_{h-2},a,b)\in
C\};$$ then $E$ is an equivalence relation on $A$. Reflexivity and
symmetry of $E$ immediately follow from the corresponding
properties of $C$. To see $E$ is transitive, let $(a,b), (b,c) \in
E$ be given. Then since $C$ is homogeneous, using $b$ as a
replacement element yields that for all $(a_1,,...,a_{h-2})\in
A^{h-2}$, $(a_1,...,a_{h-2},a,c)\in C$ and thus $(a,c)\in E$.
Suppose now that $E$ has $q$ equivalence classes and assume
without loss of generality that $A'=\{\alpha_1,...,\alpha_q\}$
contains one element from each equivalence class. Let $\gamma:
A\rightarrow A'$ be the function that maps each $a\in A$ to the
element in $A'$ that represents the equivalence class of $a$; that
is, $(a,\gamma(a))\in E$ for all $a\in A$. Define $C^*\subseteq
A^h$ by
$$C^*=\{(a_1,...,a_h)|(\gamma(a_1),...,\gamma(a_h))\in C\}.$$

\begin{lem}\label{q-main-8.1}
    Let $(a,b)\in E$ and $(a_1,...,a_{h-1}) \in A^{h-1}$. Then
    $(a_1,...,a_{h-1},a)\in C$ iff $(a_1,...,a_{h-1},b)\in
    C$.
\end{lem}

\begin{proof}
    Let $(a_1,...,a_{h-1},a)\in C$. By definition of $E$ and the total symmetry of $C$,
    $(a_1,...,a_{i-1},a,a_{i+1},...,a_{h-1},b)\in C$ for all $1\leq
    i\leq h-1$. Thus, if we use $a$ as a replacement element, the
    homogeneity of $C$ implies $(a_1,...,a_{h-1},b)\in C$.
\end{proof}

\begin{thm}\label{q-main-8.2}
    $C^*=C$. That is, membership of $C$ is completely determined by
    the equivalence classes of $E$.
\end{thm}

\begin{proof}
    First, let $(a_1,...,a_h)\in C$. Then by the previous lemma,
    $(\gamma(a_1),a_2,...,a_h)\in C$. Hence by induction, $(\gamma(a_1),\gamma(a_2),...,\gamma(a_h))\in
    C$ so that $(a_1,...,a_h)\in C^*$. Conversely, if $(a_1,...,a_h)\in
    C^*$, then $(\gamma(a_1),\gamma(a_2),...,\gamma(a_h))\in C$, and applying the
    same induction backwards yields $(a_1,...,a_h)\in C$.
\end{proof}

Call $C$ \emph{universal} if there exists a function
$f:h^{h^\kappa}\rightarrow A$ such that for $1\leq j\leq \kappa$,
$f(\pi_j^\kappa)=\alpha_j$ and such that for all $(b_1,...,b_h)\in
\omega_{h^\kappa}$, $(f(b_1),...,f(b_h))\in C$. Our next goal is
to prove that $C$ is universal. For $h\leq i\leq \kappa$ define
$\bar{C}_i\subseteq A^i$ by $(a_1,...,a_i)\in \bar{C}_i$ iff there
is an $f:h^{h^i}\rightarrow A$ such that for $1\leq j\leq i$,
$f(\pi_j^i)=a_j$ and such that for all $(b_1,...,b_h)\in
\omega_{h^i}$, $(f(b_1),...,f(b_h))\in C$. We will prove that
$\bar{C}_\kappa=A^\kappa$ to show that $C$ is universal.

\begin{lem}\label{q-main-8.3}
    $\bar{C}_h=A^h$.
\end{lem}

\begin{proof}
    Let $(a_1,...,a_h)\in A^h$ be given and let $z=(0,...,h-1)\in
    h^h$. Define $\omega:h\rightarrow A$ by $\omega (j)=a_{j+1}$
    for $0\leq j\leq h-1$ and define $f_\omega :h^{h^h}\rightarrow A$
    by $f_\omega (b)=\omega (b(z))$ for all $b\in h^{h^h}$. Then
    $f_\omega (\pi_j^h)=\omega (\pi_j^h(z))=\omega (j-1)=a_j$ for
    $1\leq j\leq h$. Moreover, if $(b_1,...,b_h)\in \omega_{h^h}$,
    then $(f_\omega (b_1),...,f_\omega (b_h))=(\omega
    (b_1(z)),...,\omega (b_h(z)))\in \iota_h^A\subseteq C$ since
    $(b_1(z),...,b_h(z))\in \iota_h^h$ by definition of $\omega_{h^h}$. Thus, $(a_1,...,a_h)\in C_h$.
\end{proof}

\begin{lem}\label{q-main-8.4}
    Let $h\leq i\leq \kappa$. Then $(\bar{C}_i,F)$ is a totally symmetric subalgebra
    of $\A^i$.
\end{lem}

\begin{proof}
    Let $g$ be an $n$-ary operation of $\A$ and for $1\leq j\leq n$ let
    ${a_j}\in \bar{C}_i$ via the function $f_j:h^{h^i}\rightarrow
    A$. Set $f=g(f_1,...,f_n)$ and write
    $g({a_1},...,{a_n})=(d_1,...,d_i)\in A^i$. We will
    show that ${d}=(d_1,...,d_i)\in \bar{C}_i$ via $f$. First note
    that if ${a_j}=(a_{j1},...,a_{ji})$, then
    $f(\pi_j^i)=g(f_1(\pi_j^i),...,f_n(\pi_j^i))=g(a_{1j},...,a_{nj})=d_j$.
    Moreover, if $(b_1,...,b_h)\in \omega_{h^i}$, then
    $(f(b_1),...,f(b_h))=(g(f_1(b_1),...,f_n(b_1)),...,g(f_1(b_h),...,f_n(b_h)))\in
    C$ since already $(f_j(b_1),...,f_j(b_h))\in C$ for $1\leq
    j\leq n$ and since $\C$ is a subalgebra of $\A^h$. Hence,
    ${d}$ is indeed an element of $\bar{C}_i$ via $f$ and thus
    $(\bar{C}_i,F)$ is a subalgebra of $\A^i$. To show that
    $\bar{C}_i$ is totally symmetric, let $\sigma$ be any
    permutation of $\{1,...,i\}$ and let $(a_1,...,a_i)\in
    \bar{C}_i$ via $f$. Then $(a_{\sigma (1)},...,a_{\sigma (i)})\in
    \bar{C}_i$ via $f_\sigma$ if we set $f_\sigma
    (b)=f(\tilde{b})$
    where $\tilde{b}(x_1,...,x_i)=b(x_{\sigma (1)},...,x_{\sigma
    (i)})$.
    For $f_\sigma (\pi_j^i)=f(\pi_{\sigma (j)}^i)=a_{\sigma (j)}$ and if
    $(b_1,...,b_h)\in\omega_{h^i}$, then
    also $(\tilde{b_1},...,\tilde{b_h})\in\omega_{h^i}$ so that
    $(f_\sigma (b_1),...,f_\sigma
    (b_h))=(f(\tilde{b_1}),...,f(\tilde{b_h}))\in C$.
\end{proof}

\begin{lem}\label{q-main-8.5}
    $C$ is universal.
\end{lem}

\begin{proof}
    We will prove by induction that for $h\leq n\leq \kappa$, $\bar{C}_n=A^n$.
    By Lemma \ref{q-main-8.3}, $\bar{C}_h=A^h$. Now assume
    $\bar{C}_n=A^n$; we will show that this implies that $\bar{C}_{n+1}$ is
    totally reflexive. Let $(a_1,...,a_n)\in A^n=\bar{C}_n$ via $f$ and
    define $\tilde{f}:h^{h^{n+1}}\rightarrow A$ by
    $\tilde{f}(b(x_1,...,x_{n+1}))=f(b(x_2,x_2,...,x_{n+1}))$.
    Then it is easy to see that $(a_1,a_1,a_2,...,a_n)\in \bar{C}_{n+1}$
    so that since $(a_1,...,a_n)$ was an arbitrary tuple in $A^n$
    we have that $\bar{C}_{n+1}$ is totally reflexive and must
    therefore
    equal $A^{n+1}$. Now in particular $\bar{C}_\kappa=A^\kappa$ and hence,
    $(\alpha_1,...,\alpha_\kappa)\in\bar{C}_\kappa$ which means
    exactly that $C$ is universal.
\end{proof}

In the light of Theorem \ref{q-main-8.2}, it is natural to
consider $D=C\cap (A')^h$. $D$ is totally reflexive, totally
symmetric and strongly homogeneous, the latter since Lemma
\ref{q-main-8.1} implies that in the definition of strong
homogeneity, we can replace $(v_1,...,v_h)$ with
$(\gamma(v_1),...,\gamma(v_h))$. Fix $f':h^{h^\kappa}\rightarrow
A$ making $C$ universal and set $f=\gamma\circ
f':h^{h^\kappa}\rightarrow A'$. $f$ makes $D$ kind of universal in
the sense that for all $(b_1,...,b_h)\in\omega_{h^\kappa}$,
$(f(b_1),...,f(b_h))\in D$; this is a consequence of Lemma
\ref{q-main-8.2}. However, $D$ is not a subuniverse of $\A^h$ as
$A'$ is not closed under those operations.

We will prove that there is an $\lambda$ such that
$q=|A'|=h^\lambda$. Let $s,t\in h^\kappa$, $j\in h$,
$g:h^\kappa\rightarrow h$. Define $g_j^t:h^\kappa\rightarrow h$ by
$$g_j^t(s)=\begin{cases}g(s)&,s\neq t\\j&,otherwise\end{cases}$$
Define $B_j^t=\{b\in h^{h^\kappa}|b(t)=j\}$.

\begin{lem}\label{q-main-8.6}
    Let $(a_1,...,a_h)\in \iota_h^h$, $t\in h^\kappa$, $g:h^\kappa \rightarrow h$,
    and $b_i\in B_{a_i}^t$ for $1\leq i\leq h-2$. Then
    $(f(b_1),...,f(b_{h-2}),f(g_{a_{h-1}}^t),f(g_{a_h}^t))\in D$.
\end{lem}

\begin{proof}
    Since $D$ is kind of universal via $f$, it is enough to show
    that $(b_1,...,b_{h-2},g_{a_{h-1}}^t,g_{a_h}^t)\in
    \omega_{h^\kappa}$. Let $t'\in h^\kappa$ be given; we evaluate the tuple
    above at $t'$. If $t'=t$, then we get
    $(a_1,...,a_{h-2},a_{h-1},a_{h})\in \iota_h^h$; if $t'\neq t$, then
    we get $(b_1(t'),...,b_{h-2}(t'),g(t'),g(t'))\in \iota_h^h$.
    Therefore $(b_1,...,b_{h-2},g_{a_{h-1}}^t,g_{a_h}^t)\in
    \omega_{h^\kappa}$ by definition of $\omega_{h^\kappa}$.
\end{proof}

\begin{lem}\label{q-main-8.7}
    Let $g:h^\kappa\rightarrow h$, $t\in h^\kappa$, $b_p:h^\kappa\rightarrow h$ for
    $1\leq p \leq h$, and $(f(g_0^t),...,f(g_{h-1}^t))\in D$.
    Consider $(f(b_1),...,f(b_h))$; then for all $1\leq r <
    s \leq h$ the following holds: If we replace the $r$-th
    component, $f(b_r)$, by $f(g_{r-1}^t)$ and the $s$-th
    component, $f(b_s)$, by $f(g_{s-1}^t)$, then the resulting
    tuple is an element of $D$.
\end{lem}

\begin{proof}
    Since the $b_p$ are arbitrarily given and since $D$ is totally
    symmetric, it suffices to consider $r=h-1$ and $s=h$. If
    $(b_1(t),...,b_{h-2}(t),h-2,h-1)\in \iota_h^h$, then by setting
    $(a_1,...,a_h)=(b_1(t),...,b_{h-2}(t),h-2,h-1)$ the result
    follows from the previous lemma. Otherwise, assume without
    loss of generality that for $1\leq p\leq h-2$, $b_p(t)=p-1$.
    We will proof by induction that for $0\leq n\leq h-2$,
    $(f(b_1),...,f(b_n),f(g_n^t),...,f(g_{h-1}^t))\in D$. For
    $n=0$, this is an assumption of the lemma; suppose it
    holds for $n<h-2$. Consider
    $(b_1,...,b_{n+1},g_{n+1}^t,...,g_{h-1}^t)$. Replace any
    component other than $b_{n+1}$ by $g_n^t$. We have that
    $g_n^t(t)=n=b_{n+1}(t)$, whereas for all $t'\neq t$,
    $g_n^t(t')=g_{h-2}^t(t')=g_{h-1}^t(t')=g(t')$. Thus, each of
    the tuples that result from our replacement belongs to
    $\omega_{h^\kappa}$ so that since $D$ is kind of universal,
    applying $f$ to each coordinate of such replacement tuples
    results in a member of $D$. By induction hypothesis, replacing
    $b_{n+1}$ by $g_n^t$ and application of $f$ yields a member of
    $D$ too, so that using the
    homogeneity of $D$ with $g_n^t$ as the replacement element concludes the proof.
\end{proof}

\begin{lem}\label{q-main-8.8}
    Suppose all assumptions of Lemma \ref{q-main-8.7} hold.
    Then $f(g_0^t)=...=f(g_{h-1}^t)$.
\end{lem}

\begin{proof}
    Clearly, it suffices to show $f(g_{h-2}^t)=f(g_{h-1}^t)$.
    Apply Lemma \ref{q-main-8.7} to see that for all $a_1,...,a_{h-2}\in A'$,
    $(a_1,...,a_{h-2},f(g_{h-2}^t),f(g_{h-1}^t))\in D$ by choosing
    $b_1,...,b_{h-2}$ from $\{\pi_1^q,...,\pi_q^q\}$. Then we must have $f(g_{h-2}^t)=f(g_{h-1}^t)$
    because if $(a_1,...,a_{h-2},a,b)\in D$
    for all $a_1,...,a_{h-2}\in A'$, then $(a,b)\in E$; thus,
    $a$ and $b$ represent the same equivalence
    class of $E$ so that they must be equal.
\end{proof}

\begin{lem}\label{q-main-8.9}
    Let $g: h^\kappa\rightarrow h$, $t\in h^\kappa$, $b_p\in B_{p-1}^t$ for $1\leq p\leq
    h$ and $(f(g_0^t),...,f(g_{h-1}^t))\notin D$.
    Consider $(f(b_1),...,f(b_h))$; then for all $1\leq r < s\leq
    h$ the following holds: If we replace the $r$-th component,
    $f(b_r)$, by $f(g_{r-1}^t)$ and the $s$-th component,
    $f(b_s)$, by $f(g_{s-1}^t)$, then the resulting tuple is not
    an element of $D$.
\end{lem}

\begin{proof}
    It suffices to prove the assertion for $r=h-1$ and $s=h$. We
    will apply a similar induction as in the proof of Lemma
    \ref{q-main-8.7}: Consider
    $(b_1,...,b_n,g_n^t,...,g_{h-1}^t)$ where $n<h-2$. Then replacing any component
    other than $g_n^t$ by $b_{n+1}$ yields an element of
    $\omega_{h^\kappa}$; for $b_{n+1}(t)=g_n^t(t)=n$ and for $t'\neq
    t$, $g_n^t(t')=g_{h-2}^t(t')=g_{h-1}^t(t')=g(t')$. Therefore,
    application of $f$ to such a tuple gives us a
    member of $D$ so that if our induction assumption is
    $(f(b_1),...,f(b_{n+1}),f(g_{n+1}^t),...,f(g_{h-1}^t))\in D$,
    then the homogeneity of $D$ implies $(f(b_1),...,f(b_{n}),f(g_{n}^t),...,f(g_{h-1}^t))\in
    D$. But as we assume that $(f(g_0^t),...,f(g_{h-1}^t))\notin D$,
    by our induction we have
    $(f(b_1),...,f(b_{h-2}),f(g_{h-2}^t),f(g_{h-1}^t))\notin D$.
\end{proof}

\begin{lem}\label{q-main-8.10}
    Let $t\in h^\kappa$. Then either for all $g:h^\kappa\rightarrow h$ we
    have $f(g_0^t)=...=f(g_{h-1}^t)$ or for all $g:h^\kappa\rightarrow h$ we
    have $(f(g_0^t),...,f(g_{h-1}^t)\notin D$.
\end{lem}

\begin{proof}
    Suppose that for some $g: h^\kappa\rightarrow h$,
    $(f(g_0^t),...,f(g_{h-1}^t))\in D$.  By Lemma
    \ref{q-main-8.8}, $f(g_0^t)=...=f(g_{h-1}^t)$. Let
    $\tilde{g}:h^\kappa\rightarrow h$. It is easy to verify
    $(\tilde{g}_0^t,...,\tilde{g}_{h-2}^t,g_0^t)\in \omega_{h^\kappa}$,
    and therefore, $(f(\tilde{g}_0^t),...,f(\tilde{g}_{h-2}^t),f(g_0^t))\in
    D$. As $f(g_0^t)=f(g_{h-1}^t)$, this tuple can be written as
    $(f(\tilde{g}_0^t),...,f(\tilde{g}_{h-2}^t),f(g_{h-1}^t))\in
    D$. But then assuming
    $(f(\tilde{g}_0^t),...,f(\tilde{g}_{h-1}^t))\nin D$
    and application of Lemma \ref{q-main-8.9} by replacing
    the first two components of $(f(\tilde{g}_0^t),...,f(\tilde{g}_{h-2}^t),f(g_{h-1}^t))$
    leads to a contradiction: The vector stays the same but is supposed to result
    in a vector not in $D$. Therefore, $(f(\tilde{g}_0^t),...,f(\tilde{g}_{h-1}^t))\in D$ so that
    $f(\tilde{g}_0^t)=...=f(\tilde{g}_{h-1}^t)$.
\end{proof}

Let $T=\{t_1,...,t_\lambda\}$ be the subset of $h^\kappa$
containing all $t$ such that for some (or all)
$g:h^\kappa\rightarrow h$, $(f(g_0^t),...,f(g_{h-1}^t))\notin D$.
Denote by $S$ the complement of $T$ in $h^\kappa$, that is, $t\in
S$ iff for some (or all) $g:h^\kappa\rightarrow h$,
$f(g_0^t)=...=f(g_{h-1}^t)$. For $g:h^\kappa\rightarrow h$, let
$\hat{g}=(g(t_1),...,g(t_\lambda))\in h^\lambda$.

\begin{lem}\label{q-main-8.11}
    Let $g_1, g_2:h^\kappa\rightarrow h$ with $\hat{g}_1=\hat{g}_2$.
    Then $f(g_1)=f(g_2)$.
\end{lem}

\begin{proof}
    Since $\hat{g}_1=\hat{g}_2$, $g_1$ and $g_2$ differ only on
    $S$. But if $s \in S$, we have by the previous lemma that for any $g:h^\kappa\To h$ and
    all $1\leq i \leq h-1$, $f(g_i^s)=f(g)$. Thus, we may alter the values
    of $g_1$ on $S$ to those of $g_2$ without changing its image
    under $f$ and the assertion follows.
\end{proof}

\begin{lem}\label{q-main-8.12}
    Let $g_1, g_2:h^\kappa\rightarrow h$ with $\hat{g}_1\neq\hat{g}_2$.
    Then $f(g_1)\neq f(g_2)$.
\end{lem}

\begin{proof}
    Let $t\in T$ with $g_1(t)\neq g_2(t)$. Assume without loss of
    generality that $g_1(t)=0$ and $g_2(t)=h-1$. Then clearly,
    $(g_1)_0^t=g_1$ and $(g_2)_{h-1}^t=g_2$. By definition of $T$,
    $(f((g_1)_0^t),...,f((g_1)_{h-1}^t))\nin D$ so that
    application of Lemma \ref{q-main-8.9} by replacement of the
    first two components of $(f((g_2)_0^t),...,f((g_2)_{h-1}^t))$
    yields $(f(g_1),f((g_1)_1^t),f((g_2)_2^t),...,f((g_2)_{h-2}^t),f(g_2)) \nin
    D$. Therefore, since $D$ is totally reflexive, $f(g_1)\neq
    f(g_2)$.
\end{proof}

\begin{lem}\label{q-main-8.13}
    $|A'|=h^\lambda$.
\end{lem}

\begin{proof}
    Clearly, $f:h^{h^\kappa}\To A'=\{\alpha_1,...,\alpha_q\}$ is onto since $f(\pi_j^\kappa)=\alpha_j$ for
    $1\leq j\leq q$. Thus, by Lemmas \ref{q-main-8.11} and
    \ref{q-main-8.12}, $|A'|=|\{f(g)|g\in h^{h^\kappa}\}|=|\{\hat{g}|g\in
    h^{h^\kappa}\}|=|h^\lambda|$, the latter equality holding because for
    every tuple $(i_1,...,i_\lambda)\in h^\lambda$ there is a $g\in h^{h^\kappa}$
    such that $g(t_j)=i_j$, $j=1,...,\lambda$.
\end{proof}

We define $\varphi ':A'\To h^\lambda$ as follows: For $a\in A'$,
$\varphi'(a)=\hat{g}$, where $g$ is an arbitrary element of
$h^{h^\kappa}$ satisfying $f(g)=a$. By Lemma \ref{q-main-8.12},
$\varphi'$ is well-defined, by Lemma \ref{q-main-8.11} it is
one-one and so together with Lemma \ref{q-main-8.13} we have that
$\varphi'$ is a bijection.

\begin{lem}\label{q-main-8.14}
    Let $g_1,...,g_h: h^\kappa\To h$ such that there is $t\in T$ with
    $(g_1(t),...,g_h(t)) \nin \iota_h^h$. Then $(f(g_1),...,f(g_h))\nin
    D$.
\end{lem}

\begin{proof}
    Assume without loss of generality that $g_i(t)=i-1$, $1\leq i
    \leq h$. Let $\tilde{0}:h^\kappa\To h$ be the zero function,
    that is, $\tilde{0}(s)=0$ for all $s\in h^\kappa$. Then for all
    $0\leq i,j\leq h-1$ with $i\neq j$,
    $(\tilde{0}_0^t,...,\tilde{0}_{i-1}^t,g_{j+1},\tilde{0}_{i+1}^t,...,\tilde{0}_{h-1}^t)\in\omega_{h^\kappa}$:
    For $s\neq t$, evaluating the tuple at $s$ yields
    $(0,...,0,g_{j+1}(s),0,...,0)\in \iota_h^h$ and evaluating the tuple
    at $t$ results in $(0,...,i-1,j,i+1,...,h-1)\in \iota_h^h$ as $i\neq
    j$. Thus,
    $(f(\tilde{0}_0^t),...,f(\tilde{0}_{i-1}^t),f(g_{j+1}),f(\tilde{0}_{i+1}^t),...,f(\tilde{0}_{h-1}^t))\in
    D$. Now suppose $(f(g_1),...,f(g_h))\in D$; then if we make
    use of the strong homogeneity of $D$ by taking
    $(f(g_1),...,f(g_h))$ as a replacement vector we get that
    $(f(\tilde{0}_0^t),...,f(\tilde{0}_{h-1}^t))\in D$.
    Hence by Lemma \ref{q-main-8.8},
    $f(\tilde{0}_0^t)=...=f(\tilde{0}_{h-1}^t)$. But $t\in T$ and
    obviously $\tilde{0}_0^t(t)=0\neq 1=\tilde{0}_1^t(t)$ so that
    by Lemma \ref{q-main-8.12}, $f(\tilde{0}_0^t)\neq
    f(\tilde{0}_1^t)$, contradiction. Therefore, we must have $(f(g_1),...,f(g_h))\nin
    D$.
\end{proof}

The last step is to show

\begin{lem}\label{q-main-8.15}
    $\varphi '(D)=\omega_\lambda$.
\end{lem}

\begin{proof}
    Let $(a_1,...,a_h)\in D$, and choose $g_i:h^\kappa\To h$ for $1\leq i \leq h$ such
    that $f(g_i)=a_i$. Then by the previous lemma,
    $(g_1(t),...,g_h(t))\in \iota_h^h$ for all $t\in T$ so that
    $(\hat{g}_1,...,\hat{g}_h)=(\varphi '(a_1),...,\varphi
    '(a_h))\in\omega_\lambda$. Conversely every element of $\omega_\lambda$
    can clearly be written as $(\hat{g}_1,...,\hat{g}_h)=(\varphi '(a_1),...,\varphi
    '(a_h))$ for some $g_i:h^\kappa\To h$ and some $a_i\in A'$, $1\leq i\leq
    h$. We may assume that $g_i(s)=0$ for all $s\in S$ and all $1\leq i\leq
    h$. But then for arbitrary $u\in h^\kappa$, $(g_1(u),...,g_h(u))\in
    \iota_h^h$ so that $(g_1,...,g_h)\in\omega_{h^\kappa}$. Hence,
    $(f(g_1),...,f(g_h))=(a_1,...,a_h)\in D$.
\end{proof}

\begin{thm}
    $C$ is $h$-regularly generated.
\end{thm}

\begin{proof}
    By Lemma \ref{q-main-8.2}, $\gamma^{-1}(D)=C^*=C$ so that by
    the previous lemma, $C=(\varphi '\circ \gamma)^{-1}(\omega_\lambda)$.
    Since Lemma \ref{q-main-8.13} implies that $\varphi
    '\circ \gamma: C\To \omega_\lambda$ is onto, the theorem follows.
\end{proof}

\subsubsection{The case $m=2$}

Recall that in the beginning we established that for $m\geq 4$,
$\B$ must be totally reflexive and totally symmetric. We showed
that this case is impossible; now we will consider the case $m=2$.

First note that for $0\in Con(\A)$ we have that either $0\cap B=0$
(that is, $0\subseteq B$) or $0\cap B=\emptyset$; for otherwise,
the projection of that intersection on one coordinate would be a
proper subalgebra of $\A$. In the first case $\B$ is reflexive, in
the second case we call $\B$ \emph{areflexive}. Observe that if
$\B$ is areflexive then, as a binary relation, it can neither have
a least nor a greatest element.

For two binary relations $C_1, C_2$ we denote the relation product
by $C_1\cdot C_2$; we define the inverse relation of $C_1$ to be
$C_1^{-1}=\{(a,b)|(b,a)\in C_1\}$.

\begin{lem}\label{q-main-5.1}
    If $B\cdot B^{-1}=A^2$, then $B$ has a least and
    a greatest element.
\end{lem}

\begin{proof}
    We will prove that $B$ has a greatest element. Since $B\cdot B^{-1}=A^2$ if and only if $B^{-1}\cdot
    B=A^2$ this implies that $B^{-1}$ has a greatest element as well so that $B$ has a
    least element. Define for $2\leq i\leq \kappa$ sets $C_i\subseteq A^i$ by
    $$
        C_i=\{\,(a_1,...,a_i)\,|\,\exists b\in A\,\forall\, 1\leq j\leq i
        \,\,(a_j,b)\in B\}.
    $$
    Then clearly $C_2=B\cdot B^{-1}=A^2$ and $(C_i,F)$ is a
    totally symmetric subalgebra of $\A^i$ for all $2\leq i\leq \kappa$.
    Moreover, $C_i=A^i$ obviously implies that $C_{i+1}$ is
    totally reflexive so that it must be equal to $A^{i+1}$. By
    induction we get $C_\kappa=A^\kappa$ and so $(\alpha_1,...,\alpha_\kappa)\in
    C_\kappa$. But that means there exists $o\in A$ such that
    $(a,o)\in B$ for all $a\in A$; hence, $o$ is a greatest element of $B$.
\end{proof}

Now if $B$ is areflexive, it has no greatest element so that by
the previous lemma, $B\cdot B^{-1}\neq A^2$. Since $B\cdot B^{-1}$
is reflexive, $|B\cdot B^{-1}|\geq |A|$. Furthermore, $|B|>|A|$ so
that there exist $a, b, c\in A$ with $a\neq b$ such that $(a,c)\in
B$ and $(b,c)\in B$. Thus, $(a,b)\in B\cdot B^{-1}$ and so
$|B\cdot B^{-1}|> |A|$. Hence, as ($B\cdot B^{-1},F)$ is a
subalgebra of $\A$, we can consider $B\cdot B^{-1}$ instead of
$B$; we will therefore assume from now on that $B$ is reflexive.

We call $B$ \emph{antisymmetric} iff $B\cap B^{-1}=0$. Suppose $B$
is not antisymmetric; then since $B$ is reflexive, $B\cap B^{-1}$
properly contains $0$. In that case we may as well assume $B$ is
symmetric by replacing $B$ with $B\cap B^{-1}$ which trivially is
symmetric. Hence, we are back in the totally reflexive totally
symmetric case which we have already shown absurd.

So we assume $B$ is antisymmetric as well. Then the following is
true.

\begin{lem}
    $B$ has exactly one least and exactly one greatest element.
\end{lem}

\begin{proof}
    Note that $B$ can have at most one least and one greatest element because of its antisymmetry.
    Since $B\cdot B^{-1}$ contains $0$ properly and since it is
    symmetric, it must equal $A^2$ so that by the last lemma, $B$ has at least one least and at least one greatest
element.
\end{proof}

\begin{lem}\label{q-main-5.2}
    $B\cdot B\neq A^2$.
\end{lem}

\begin{proof}
    Let $a\in A$ and let $o$ be the greatest element of $B$. If $(o,a)\in
    B\cdot B$, then there exists $b\in A$ such that $(o,b)\in B$
    and $(b,a)\in B$. Because $o$ is the greatest element of $B$, $(b,o)\in B$
    and so, since $B$ is antisymmetric, $b=o$. Hence, $(o,a)\in B$
    so that $o=a$. Thus, if $a\neq o$, then $(o,a)\nin B\cdot B$.
\end{proof}

Assume that $\B$ is maximal among the antisymmetric subalgebras of
$\A^2$. We will finish the case $m=2$ and show $B$ is a partial
order with least and greatest element; the only thing that is
missing is the transitivity if $B$.

\begin{lem}
    $B$ is transitive.
\end{lem}

\begin{proof}
    Consider the subalgebra $(B\cdot B,F)$ of $\A^2$. By the previous lemma, $B\cdot
    B\neq A^2$. Suppose $B\cdot B$ is not antisymmetric. Then
    $(B\cdot B)\cap (B\cdot B)^{-1}$ is strictly between $0$ and $1\in Con(\A)$; since it is both symmetric and reflexive, this is
    impossible. Thus, $B\cdot B$ is antisymmetric, and as $B\cdot
    B \supseteq B$, the maximality of $B$ implies $B\cdot B=B$.
    Hence, $B$ is transitive.
\end{proof}

We are left with the case $m=3$ as we have already eliminated the
cases $m\geq 4$ and $m=2$ as possibilities.

\subsubsection{The case $m=3$}
We will show that the case where $\B\leq\A^3$ is impossible as
well to finish the proof. Denote by $\Delta_n(A)$ the diagonal of
$A^n$.

\begin{lem}
    $B_{(1,2,3)}=\Delta_3(A)$, that is, $B$ contains
    $\Delta_3(A)$.
\end{lem}

\begin{proof}
    Since $B_{(1,2)}$ is essentially a subuniverse of $\A^2$ we have that
    $|B_{(1,2)}|$ equals either $|A|^2$ or $|A|$. In the first
    case our assertion follows immediately. In the second case
    consider the subalgebra $(\rho_1(B_{(1,2)}),F)$ of $\A^2$
    and note that $|\rho_3(B)|=|A|^2$
    implies that $\pi_1^2(\rho_1(B_{(1,2)}))=A$. But since
    $|\rho_1(B_{(1,2)})|=|B_{(1,2)}|=|A|$, $\rho_1(B_{(1,2)})$ is
    the graph of an automorphism of $\A$ and must therefore equal
    $\Delta_2(A)$. Hence, $B_{(1,2)}=\Delta_3(A)$ and so also
    $B_{(1,2,3)}=\Delta_3(A)$.
\end{proof}

We will show now that we can assume that $B\cap
\iota_3^A\supsetneqq \Delta_3(A)$. Define a subuniverse $B'$ of
$\A^3$ by
$$
    B'=\{(a_2,a_3,a_3')\,|\,\exists \,a_1\in A \,((a_1,a_2,a_3)\in B
    \wedge(a_1,a_2,a_3')\in B)\}.
$$
Since $|\rho_1(B)|=|A|^2$ by Lemma \ref{q-main-4.1},
$|B_{(2,3)}'|=|A|^2$. Therefore, $B'\cap \iota_3^A\neq
\Delta_3(A)$. It is possible that $B\cap \iota_3^A= \Delta_3(A)$.
But in that case, $|B'_{(1,2)}|=|A|$ and thus $|B'| < |A|^3$;
moreover, since $|\rho_3(B)|=|A|^2$ and $|B|>|A|^2$, $|B'|>|A|^2$.
Hence, we can replace $B$ by $B'$ in our proof and so we will
assume from now on that $B\cap \iota_3^A\neq \Delta_3(A)$. Up to
symmetry, this leaves us with three possibilities (since
$B_{(i,j)}$ is essentially a subuniverse of $\A^2$ and must
therefore have cardinality a power of $|A|$):

\begin{lem}\label{q-main-poss}
    Either
    \begin{enumerate}
        \item{$|B_{(1,2)}|=|A|^2$ and $|B_{(1,3)}|=|B_{(2,3)}|=|A|$ or}
        \item{$|B_{(1,2)}|=|B_{(1,3)}|=|A|^2$ and $|B_{(2,3)}|=|A|$ or}
        \item{$|B_{(1,2)}|=|B_{(1,3)}|=|B_{(2,3)}|=|A|^2$}
    \end{enumerate}
\end{lem}

\begin{lem}
    Possibility 3 is impossible.
\end{lem}

\begin{proof}
    If Possibility 3 was true, then $\iota_3^A\subseteq B$ and so also the
    subuniverse generated by $\iota_3^A$ would be a subset of $B$ and
    thus a proper subuniverse of $\A^3$. But this is impossible as
    that subuniverse is totally reflexive and totally symmetric.
\end{proof}

\begin{lem}
    Possibility 1 is impossible.
\end{lem}

\begin{proof}
    Define for $2\leq i\leq \kappa$ sets $B_i \subseteq A^i$ by
    $$
        B_i=\{(a_1,...,a_i)\,|\,\exists a\in A \,\exists b\in A \,\forall
        \, 1\leq j\leq i \quad (a_j,a,b)\in B\}
    $$
    Clearly, $B_i$ is a totally symmetric subuniverse of $\A^i$. Since
    $\Delta_3(A)\subseteq B$, $\Delta_2 (A)\subseteq B_2$. Also,
    $|B|>|A|^2$ and so there are distinct $a, b, c\in A$ such
    that $(a,b,c)\in B$. But also $(b,b,c)\in B$ so that $(a,b)\in
    B_2$, and therefore $|B_2|>|A|$. Hence, as $|B_2|$ must be a power of $|A|$, we must have
    $|B_2|= |A|^2$. Now since $B_n=A^n$ implies that $B_{n+1}$ is
    totally reflexive, it implies further $B_{n+1}=A^{n+1}$. By induction we get
    $B_\kappa=A^\kappa$ and as a consequence, $(\alpha_1,...,\alpha_\kappa)\in B_\kappa$. But
    that means that there exist $a,b\in A$ such that for all $x\in
    A$ we have $(x,a,b)\in B$. Setting $x=b$ yields $(b,a,b)\in
    B$, and since $|B_{(1,3)}|=|A|$ and $\Delta_3(A)\subseteq B$, we
    conclude $a=b$. But if we choose now $x\neq a$, then we get
    $(x,a,a)\in B$, contradicting $|B_{(2,3)}|=|A|$. Hence, Case 1
    cannot occur.
\end{proof}

We will conclude the proof of Theorem \ref{q-main-3.4} by showing
that Possibility 2 in \ref{q-main-poss} is impossible as well.
This will require more effort than the other cases.

\begin{lem}
    If for $a,b,c\in A$ both $(a,b,c)\in B$ and $(b,a,c)\in B$,
    then $a=b$.
\end{lem}

\begin{proof}
    Define a subuniverse $B'$ of $\A^3$ by
    $$
        B'=\{(b,a,c)\,|\, (a,b,c)\in B\}
    $$
    and consider the subuniverse $B\cap B'$. Because $|B_{(2,3)}|=|A|$, $(a,b,b)\in B$
    implies $a=b$ and clearly $(b,a,b)\in B'$ implies the same; thus,
    $B\cap B'\cap \iota_3^A=\{(a,a,b)\,|\, a, b\in A\}$. Therefore $B\cap
    B'$ satisfies exactly the equalities of Case 1, and so, if $|B\cap
    B'|>|A|^2$, we have a contradiction. Hence, $|B\cap B'|=|A|^2$, or
    equivalently, $B\cap B'=B\cap B'\cap \iota_3^A=\{(a,a,b)\,|\, a, b\in
    A\}$ which means exactly that $(a,b,c)\in B$ and $(b,a,c)\in B$ implies
    $a=b$.
\end{proof}

For an equivalence relation $\sim$ on $\{1,...,n\}$ we define
$\Delta_\sim\subseteq A^n$ by
$$
    \Delta_\sim =\{(a_1,...,a_n)\,|\, i\sim j\To a_i=a_j\}.
$$
In this context, we will denote an equivalence relation by writing
down its equivalence classes. Let $(C,F)$ be the subalgebra of
$\A^4$ generated by
$$
    \Delta_{\{1,2\}\{3,4\}}\cup \Delta_{\{1,3\}\{2,4\}}\cup
    \Delta_{\{1,4\}\{2,3\}};
$$
clearly, $C$ is totally symmetric. Also for $1\leq i\leq 4$,
$\rho_i(C)\supseteq \iota_3^A$, and so $\rho_i(C)= A^3$.

\begin{lem}\label{q-main-6.3}
    For $\{i,j,r,s\}=\{1,2,3,4\}$, we have
    $C_{(i,j)}=C_{(i,j)(r,s)}$.
\end{lem}

\begin{proof}
    Define a subuniverse $B'$ of $\A^4$ by
    \begin{eqnarray*}
        \begin{split}
            B'=\{&(a_1,a_2,a_3,a_4)|\exists a_5\in A \exists a_6\in A \\
                &((a_1,a_2,a_5)\in B\,\wedge \,(a_5,a_3,a_4)\in B\,\wedge
                        \,(a_2,a_1,a_6)\in B\,\wedge\, (a_6,a_3,a_4)\in B)\}\\
        \end{split}
    \end{eqnarray*}
    By our assumptions for Case 2, one
    easily checks that $C\subseteq B'$. Note next that if
    $(a_1,a_2,a_3,a_3)\in B'$, then $a_6=a_5=a_3$ so that
    $(a_1,a_2,a_3)\in B$ and $(a_2,a_1,a_3)\in B$ which implies
    $a_1=a_2$. The same property holds for $C$ since $C\subseteq B'$.
    Hence, $C_{(3,4)}=C_{(3,4)(1,2)}$ and the lemma follows
    from the total symmetry of $C$.
\end{proof}

\begin{lem}\label{q-main-6.4}
    Let $a_1, a_2, a_3\in A$ be given. Then there exists exactly
    one $a_4\in A$ such that $(a_1,a_2,a_3,a_4)\in C$.
\end{lem}

\begin{proof}
    The existence of such an $a_1$ follows immediately from
    $\rho_1(C)=A^3$. We will show that $a_1$ is unique. Define
    $B'\leq A^3$ to be
    $$
        B'=\{(a_1,a_2,a_3)\,|\,\exists a_4\in A\,\exists a_5\in
        A\quad ((a_1,a_3,a_4,a_5)\in
        C\,\wedge\,(a_2,a_3,a_4,a_5)\in C)\}
    $$
    It is obvious that $B'$ is a subuniverse of $\A^3$. Clearly,
    $|B_{(1,2)}'|=|A|^2$. If $(a_1,a_2,a_2)\in B'$, then by the
    definition of $B'$, $(a_2,a_2,a_4,a_5)\in C$ so that by the
    previous lemma, $a_1=a_2$. Thus, $|B_{(2,3)}'|=|A|$, and by
    the same argument, $|B_{(1,3)}'|=|A|$. But as we have already
    proven Case 1 impossible, we must have $|B'|=|A|^2$ which
    implies $B'=B_{(1,2)}'$. So if we have $(a_1,a_2,a_3,a_4)\in
    C$ and $(a_1',a_2,a_3,a_4)\in C$, then $a_1=a_1'$.
\end{proof}

A consequence of the previous lemma is that we can define a
function $f: A^3\To A$ assigning to $(a_2,a_3,a_4)\in A^3$ the
unique $a_1\in A$ such that $(a_1,a_2,a_3,a_4)\in C$.

\begin{lem}\label{q-main-6.4.b}
    For all permutations $\pi$ of $\{a,b,c\}$ we have that
    $f(a,b,c)=f(\pi(a),\pi(b),\pi(c))$. Moreover, $f(a,a,b)=b$,
    and $f(a,b,f(a,b,c))=c$.
\end{lem}

\begin{proof}
    The first assertion holds because $C$ is totally symmetric. The
    second one is a direct consequence of Lemma
    \ref{q-main-6.3}. For the last one, note that by the definition
    of $f$, $(a,b,c,f(a,b,c))\in C$. Thus, again by the definition
    of $f$ and the total symmetry of $C$, $f(a,b,f(a,b,c))=c$.
\end{proof}

\begin{lem}\label{q-main-6.6}
    $f$ satisfies the equation $f(a,b,c)=f(f(a,d,c),d,b)$.
\end{lem}

\begin{proof}
    Define $C'$ to be the subset of $A^4$ containing exactly the
    tuples $(a,b,c,d)$ for which $f$ satisfies the equation of the lemma; $C'$ is a subuniverse of $\A^4$. It is easy to
    check with the properties of $f$ stated in Lemma
    \ref{q-main-6.4.b} that for all $1\leq i<j\leq4$,
    $|C_{(i,j)}'|=|A|^3$. Thus, $|C'|>|A|^3$ and consequently
    $C'=A^4$. This proves the lemma.
\end{proof}

Now choose $0\in A$ arbitrarily and define a binary operation $+$
on $A$ by
$$
    a+b=f(a,b,0).
$$
This will give us a prime affine relation and lead the last
possibility ad absurdum.

\begin{lem}\label{q-main-6-group}
    $(A,+)$ is an abelian $2$-group.
\end{lem}

\begin{proof}
    Checking the associative law, we use Lemmas
    \ref{q-main-6.4.b} and \ref{q-main-6.6} to calculate
    \begin{eqnarray*}
    \begin{split}
        a+(b+c)&=f(a,b+c,0)\\
            &=f(a,f(b,c,0),0)\\
            &=f(f(b,0,c),0,a)\\
            &=f(b,a,c)\\
            &=f(a,b,c)
    \end{split}
    \end{eqnarray*}
    A similar computation yields  $(a+b)+c=f(a,b,c)$ so that
    $a+(b+c)=(a+b)+c$. $0\in A$ is the neutral element since for all
    $a\in A$, $a+0=0+a=f(a,0,0)=a$ by Lemma \ref{q-main-6.4.b}.
    As by Lemma \ref{q-main-6.4.b} we have $a+a=f(a,a,0)=0$, each element is its
    own inverse. Observe that this also means the group $(A,+)$ is
    a $2$-group. Finally, the group is abelian since $a+b=f(a,b,0)=f(b,a,0)=b+a$ for all
    $a,b\in A$.
\end{proof}

Now all tuples in $C$ are of the form $(a,b,c,f(a,b,c))$, which we
know can be written as $(a,b,c,a+b+c)$. If we set $c=0$ then by
the fact that $C$ is a subuniverse of $\A^4$ we get that for any
$n$-ary operation $g$ of $\A$, $g(a)+g(b)+g(0)=g(a+b)$, where
$a,b\in A^n$ are arbitrary. Define $\rho\subseteq A^4$ by
$$
    (a,b,c,d)\in\rho\gdw a+b=c+d.
$$
Then $\rho$ is a subuniverse of $\A^4$ since for $n$-ary $g$ we
have that if $a+b=c+d$, where $a,b,c,d\in A^n$, then
$g(a+b)=g(c+d)$ and so by the preceding discussion
$g(a)+g(b)=g(c)+g(d)$. Hence, $\rho$ is a prime affine relation
with respect to the abelian $2$-group $(A,+)$ and so forbidden by
Rosenberg's list. We have therefore eliminated Possibility 2 as a
possibility in \ref{q-main-poss} and so finally finished the case
$m=3$. Theorem \ref{q-main-3.4} has been proven.

    \section{$\VV(\A)$ is congruence permutable}

We will use the result of the last section, namely that $\A$ has
almost minimal spectrum, to show that the variety generated by
$\A$ is congruence permutable. In the beginning of our proof, we
will follow another result by R. Quackenbush in \cite{Qua71}; for
the second part we will go another way than the one shown there.

\begin{defn}
    An algebra $\A$ is called \emph{congruence permutable} iff for
    all congruences $\psi, \theta$ on $\A$, $\psi\cdot\theta
    =\theta\cdot\psi$. We say a variety is congruence permutable iff
    every algebra in the variety is.
\end{defn}

\begin{thm}[R. Quackenbush \cite{Qua71}]
    Let $\A$ be a finite non-trivial algebra. If $\A$ has
    minimal spectrum, then the
    variety generated by $\A$ is congruence permutable.
\end{thm}

\begin{rem}
    The converse holds under the assumption that $\A$ is simple and
    has no proper subalgebras. For a proof of this consult
    \cite{Qua71}. Note that the assumption of the theorem is that
    $\A$ has \emph{minimal spectrum}, whereas we only know until now that
    $\A$ has \emph{almost minimal spectrum}. For the proof of
    congruence permutability, this is still sufficient.
\end{rem}

\begin{thm}\label{q-qua-1.1}
    Let $\A$ be a finite non-trivial algebra. If $\A$ has
    almost minimal spectrum, then the
    variety generated by $\A$ is congruence permutable.
\end{thm}

For a subdirect product $\B$ of algebras $\A_1,...,\A_n$, the meet
of the kernels of all projections onto the $\A_i$ is clearly the
trivial congruence $0$. Recall that $\B$ is called an
\textit{irreducible} subdirect product iff for every proper subset
of the projections the meet of the kernels of the projections of
that subset is strictly greater than $0$.

\begin{defn}
    A set $\AAA$ of finite algebras is a \emph{direct factor set} iff
    whenever $\B$ is an irreducible subdirect product of algebras
    $\A_1,...,\A_n$ of $\AAA$, then $\B = \A_1 \times ... \times
    \A_n$.
\end{defn}

If we look at a congruence relation $\theta$ on an algebra $\A$ as
a subalgebra of $\A^2$, we will denote this subalgebra by
$\A_\theta$. Note that $\A_\theta$ is a subdirect power of $\A$.

\begin{lem}\label{q-qua-2.2}
    All algebras in a direct factor set are simple.
\end{lem}

\begin{proof}
    Let $\theta$ be a congruence relation on a member $\A$ of a direct
    factor set. If $\A_\theta$ is reducible, then its
    projection onto one coordinate is one-one and hence $\theta = 0$.
    On the other hand, if it is irreducible then by the definition of
    a direct factor set $\theta=1$.
\end{proof}

\begin{thm}\label{q-qua-2.3}
    If an algebra $\A$ has almost minimal spectrum then $\{\A\}$ is a direct
    factor set.
\end{thm}

\begin{proof}
    We must prove that an irreducible subdirect product $\B$ of $n$ copies of $\A$
    is equal to $\A^n$. The proof will be by induction.
    For $n=1$ the assertion is trivial. Assume it is true for
    $n-1$ and let $\B$ be an irreducible
    subdirect product of $n$ copies of $\A$. Set $\B'=\B/ker(\pi_n^n)$, where
    $\pi_n^n$ denotes the projection onto the $n$-th coordinate.
    As one can easily see, $\B'$ is essentially an irreducible subdirect product
    of $n-1$ copies of $\A$ and hence by the induction assumption,
    $\B'\cong\A^{n-1}$. As $\B$ is also irreducible,
    $|\B|>|\B'|=|\A^{n-1}|$. But since $\B\subseteq\A^n$, we have
    that $|\B|\leq |\A^n|$ so that by the fact that $\A$ has
    almost minimal spectrum we have $|\B|=|\A^n|$ and therefore
    $\B =\A^n$.
\end{proof}

\begin{defn}
    A congruence is \textit{uniform} iff all its equivalence classes
    are of the same cardinality. An algebra is said to have uniform congruences iff all its
    congruences are uniform and a variety has uniform congruences iff
    all its algebras do.
\end{defn}

Our next goal is to prove that our algebra has uniform
congruences. We will need the following theorem.

\begin{thm}\label{q-qua-2.4}
    Let $\AAA$ be a direct factor set, $I=\{1,...,n\}$ be a
    finite index set and $(\A_i)_{i\in I}$ be algebras in
    $\AAA$. Let $\A=\prod_{i\in I}\A_i$ and let $\theta$ be a
    non-trivial congruence of $\A$. Then there exists a proper
    subset $J$ of $I$ such that
    $\A_\theta\cong\A\times\prod_{j\in J}\A_j$. In this
    case $\theta$ has $|\prod_{i \in I\setminus J}\A_i|$ equivalence
    classes each of which has $|\prod_{j\in J}\A_j|$ elements.
\end{thm}

\begin{proof}
    Since $\A_\theta\leq\A^2=\prod_{i\in I}\A_i\times \prod_{i\in
    I}\A_i$, we know that $\A_\theta$ is an irreducible subdirect
    product and hence, since $\AAA$ is a direct factor set, a direct product
    of some of the factors of $\A^2$; that is, there are subsets $K$ and $K'$ of $I$
    such that $\A_\theta\cong\prod_{i\in K}\A_i\times\prod_{i\in K'}\A_i$. We claim that each
    $\A_i$ occurs at least once in this direct product
    representation so that $K \cup K'=I$. Given $a_i \in \A_i$, $2\leq i\leq n$,
    we have that
    $(a,a_2,...,a_n,a,a_2,...,a_n) \in \A_\theta$  for every $a \in A$ because $\theta$ is
    reflexive. On the other hand the components of an element in
    $\A_\theta$ corresponding to indices in $K$ and $K'$ uniquely
    determine the other components. Therefore, $i=1$ must be in $K$
    or $K'$ and clearly the same holds for any $i \in I$. Hence, by
    reordering the factors of the direct product representation of
    $\A_\theta$ and setting $J=K\cap K'$, the first assertion of the theorem
    follows.\newline
    Now denote the equivalence classes of $\theta$ by
    $C_1,...,C_l$ and their cardinalities by $c_1,...,c_l$. First we will
    show that $c_k\leq |\prod_{j\in J}\A_j|$, $1\leq k\leq l$. Consider the
    projections
    $$
        \varsigma_k:\quad
        \begin{matrix}
            C_k && \rightarrow && \prod_{j\in J}A_j \\
            (a_1,...,a_n) && \mapsto &&(a_j)_{j\in J}
        \end{matrix}
    $$
    We claim that $\varsigma_k$ is one-one. Let $a,b \in C_k$ with
    $\varsigma_k (a)=\varsigma_k (b)$, that is, $(a_j)_{j\in J}=(b_j)_{j\in J}$;
    if we prove them equal our assertion follows. Define a vector $d$ by
    $$
        d_j=
        \begin{cases}
            b_j & ,j\in K' \\
            a_j & ,j\in I\setminus K'
        \end{cases}
    $$
    As the $b_j$ and
    $a_j$ agree on $J$, our projection
    $$
        \pi:\quad
        \begin{matrix}
            A\times A && \rightarrow && \prod_{j\in K}A_j\times\prod_{j\in K'}A_j\\
            (a_1,...,a_n,a_1',...,a_n') && \mapsto && ((a_j)_{j\in K},(a_j')_{j\in K'})
        \end{matrix}
    $$
    maps $(d,d) \in \theta$ to
    $((a_j)_{j\in K},(b_j)_{j\in K'})$. But $(a,b) \in \theta$ is
    mapped to exactly the same vector and so, as the coordinates
    in $K$ and $K'$ uniquely determine all the others, we have
    that $a_j=b_j$ for $j\in I\setminus K'$. By symmetry we
    conclude that
    $a_j=b_j$ for $j\in I\setminus K$ and hence, $a=b$ follows.
    \newline Finally, our inequality together with the obvious equalities
    $\sum_{k=1}^l c_k=|\A|$ and $\sum_{k=1}^l
    c_k^2=|\A|\,|\prod_{j\in J}\A_j|$ implies that for all $1\leq k\leq l$ we must have
    $c_k=|\prod_{j\in J}\A_j|$.
\end{proof}

Now we can establish that the finite algebras in $\VV(\A)$ have
uniform congruences. Recall that a variety is called
\textit{locally finite} iff every finitely generated algebra in it
is finite. For a set $\AAA$ of algebras of the same type, denote
by $P(\AAA)$ all products, by $S(\AAA)$ all subalgebras, and by
$H(\AAA)$ all homomorphic images of algebras of $\AAA$. Then it is
well-known that $\VV(\AAA)=HSP(\AAA)$.

\begin{thm}\label{q-qua-2.5}
    Let $\AAA$ be a finite direct factor set with the property that
    $S(\AAA)\subseteq P(\AAA)$. Then the finite algebras in $\VV(\AAA)$ have uniform
    congruences.
\end{thm}

\begin{proof}
    First let $\B\in SP(\AAA)$ be finite. If $\B$ is a subdirect product of algebras
    in $\AAA$ then $\B\in P(\AAA)$ because $\AAA$ is a direct
    factor set. If on the other hand the projection $\pi_i$ of $B$ onto
    some coordinate $i$ is not onto, then by the assumption $S(\AAA)\subseteq P(\AAA)$
    we can replace that
    coordinate with a product of algebras in $\AAA$ equal to
    $\pi_i(B)$ and we have the first case again. Hence,
    $SP(\AAA)=P(\AAA)$ and thus by the last theorem, all finite algebras in
    $SP(\AAA)$ have uniform congruences. Now let $\C\in
    \VV(\AAA)=HSP(\AAA)$ be given and assume it is finite; let
    $\psi$ be a congruence on $\C$. We want to show $\psi$ is
    uniform. Clearly, $\C\cong\B/\theta$ for some $\B\in SP(\AAA)$ and some
    (uniform) $\theta\in Con(\B)$. We can assume $\B$ is finite:
    Observe first that $|\AAA|<\aleph_0$ implies that
    $\VV(\AAA)$ is locally finite. Now if $\B$ is infinite,
    replace it by the subalgebra $\D$ generated by any finite subset of $B$ containing at least one
    representative from each $\theta$-class;
    then obviously $\C\cong\D /\tilde{\theta}$ if we set
    $\tilde{\theta}=\theta\cap D^2$. Now $\psi$ induces a congruence
    relation $\zeta$ on $\B$, defined by $a\zeta b\leftrightarrow
    [a]_\theta\psi [b]_\theta$, and every congruence class of $\zeta$
    corresponds to exactly one congruence class of $\psi$. Since $\zeta$ and $\theta$
    are uniform, $\psi$ must be uniform as well: If the size of all $\zeta$-classes
    is $n$ and the size of all $\theta$-classes is $j$, then the size of all $\psi$-classes must
     be $\frac{n}{j}$. Hence $\psi$ is uniform.
\end{proof}

Finally we have also established what we wanted earlier: $\A$ has
really minimal spectrum. Still, it is worth mentioning.

\begin{thm}
    Let $\A$ satisfy the hypotheses of Theorem \ref{Ros1}. Then $\A$ has minimal spectrum.
\end{thm}

\begin{proof}
    This follows from the proof of the previous lemma and Theorem
    \ref{q-qua-2.4}.
\end{proof}

Our next goal is to show that $\VV(\A)$ has coherent congruences.

\begin{defn}
    An algebra $\A$ is \textit{congruence coherent} iff for every
    subalgebra $\B$ of $\A$ and every congruence $\theta$ on $\A$ it is true that if
    $\B$ contains a congruence class of $\theta$, then $\B$ is the
    union of congruence classes of $\theta$; in other words, if
    $[a]_\theta \subseteq B$ for some $a \in B$ implies $[b]_\theta
    \subseteq B$ for every $b \in B$. A variety is called \textit{congruence coherent} or simply
    \textit{coherent} iff all of its members are.
\end{defn}

The following two lemmas are due to M. Clark and P. Krauss
\cite{CK76}.

\begin{lem}\label{q-ck-3.9}
    If the finite algebras in a variety $\VV$ are congruence uniform
    then the finite algebras in $\VV$ are congruence coherent.
\end{lem}

\begin{proof}
    Let $\A \in \VV$ be finite, and let $\B$ be a subalgebra and $\theta$ be a congruence of
    $\A$. If $X \subseteq B$ is an equivalence class of $\theta$,
    then it is also an equivalence class of $\theta \cap B^2$. Now
    let $Y$ be another congruence class of $\theta \cap B^2$. Then
    there exists a congruence class $Z$ of $\theta$ such that
    $Y=Z\cap B$. But by our hypothesis, $|Y|=|X|=|Z|$ and hence,
    since $Z$ is finite, $Y=Z$.
\end{proof}

\begin{lem}\label{q-ck-3.10}
    If $\VV$ is a locally finite variety and the finite algebras in
    $\VV$ are congruence coherent then $\VV$ is congruence coherent.
\end{lem}

\begin{proof}
    Let $\A \in \VV$ and let $\B$ be a subalgebra and $\theta$ be a congruence of
    $\A$. Consider an equivalence class $X$ of $\theta$, $X
    \subseteq
    B$, and consider $[b]_\theta$ for some $b \in
    B$. Now if $a\theta b$ is given, choose $x\in
    X$ and consider the restrictions to the subuniverse $[a,b,x]$ of $\A$ generated by $\{a,b,x\}$: $X\cap [a,b,x]$ is a congruence
    class of $\theta\cap [a,b,x]$ on $\A\cap [a,b,x]$ and $X\cap [a,b,x]=X\cap[a,b,x]\cap
    B$. Clearly, $a\equiv b \,(\theta \cap [a,b,x])$ and so by the
    hypothesis $a \in B \cap [a,b,x]$. Thus we have that
    $[b]_\theta \subseteq B$ and the lemma follows.
\end{proof}

\begin{cor}\label{q-ck-3.11}
    If $\VV$ is a locally finite variety and the finite algebras in
    $\VV$ are congruence uniform then $\VV$ is congruence coherent.
\end{cor}

\begin{proof}
    This is an immediate consequence of Lemmas \ref{q-ck-3.9} and
    \ref{q-ck-3.10}.
\end{proof}

We will use a version of a result on g-coherence from \cite{CE01}
to finish our proof.

\begin{lem}\label{q-ce-1}
    If a variety $\VV$ is coherent, then for some $n$ there exist
    ternary terms $t_1,...,t_n$ and an $n+1$-ary term $\omega$
    such that the following identities hold in $\VV$:
    \begin{eqnarray*}
    \begin{split}
        &t_i(x,x,z) = z, \quad i=1,...,n \\
        &y= \omega (x,t_1(x,y,z),...,t_n(x,y,z)).
    \end{split}
    \end{eqnarray*}
\end{lem}

\begin{proof}
    Consider the free algebra with three generators determined by $\VV$, $\F_3(\VV)$,
    and call the generators $x, y, z$. Let $\theta=\theta(x,y)$ be
    the congruence on $\F_3(\VV)$ generated by $\{(x,y)\}$ and let
    $\B$ be the subalgebra of $\F_3(\VV)$ generated by the set
    $\{x\}\cup [z]_\theta$. Clearly, $[z]_\theta \subseteq B$ and
    hence by the coherence of $\VV$, $[b]_\theta \subseteq B$ for
    all $b \in B$. Since $x \in B$ and $y\theta x$, also $y \in B$
    and so there exists a term $\omega$ in the language of $\VV$ such
    that $y=\omega(x,c_1,...,c_n)$, where $c_1,...,c_n \in [z]_\theta$.
    As elements of $\F_3(\VV)$ the $c_i$ have representations as terms
    $t_i(x,y,z)$ so that $y=\omega (x,t_1(x,y,z),...,t_n(x,y,z))$. Furthermore,
     since $t_i(x,y,z) \in [z]_\theta$ and since
    $\theta$ is the congruence identifying $x$ and $y$, one can
    easily derive that $t_i(x,x,z)=z$ is an identity of $\VV$ for $i=1,...,n$.
\end{proof}

The following theorem is a well-known criterion for congruence
permutability by A. Mal'cev.

\begin{thm}\label{q-ce-6}
    A variety $\VV$ is congruence permutable iff there exists a
    ternary term $p(x,y,z)$ of $\VV$ such that the
    identities
        $$p(x,x,z)=p(z,x,x)=z$$
    can be derived in $\VV$. The term $p$ is called a Mal'cev
    term.
\end{thm}

\begin{proof}
    First assume that $\VV$ is congruence permutable. Consider the
    congruences $\theta=\theta(x,y)$ and $\psi=\psi(y,z)$ generated by $\{(x,y)\}$
    and $\{(y,z)\}$, respectively, on the free algebra $\F_3(\VV)$
    with generators $x,y,z$. Clearly, $(x,z)\in \theta\cdot\psi$ and so by
    our assumption $(x,z)\in \psi\cdot\theta$. Hence, there is a
    term $p(x,y,z)$ in $\F_3(\VV)$ such that $(x,p(x,y,z))\in
    \psi$ and $(p(x,y,z),z)\in \theta$ which by the definition of those congruences
    yields $x=p(x,z,z)$ and $p(x,x,z)=z$.

    Conversely, let $p(x,y,z)$ be a term of $\VV$ satisfying
    $p(x,x,z)=p(z,x,x)=z$, let $\A$ be an arbitrary algebra of
    $\VV$ and let $\theta,\,\psi\in Con(\A)$ be any two
    congruences on $\A$. If $(a,b) \in \theta\cdot\psi$ then there
    exists $c\in A$ with $(a,c)\in\theta$ and $(c,b)\in\psi$.
    Since trivially $(a,a), (b,b) \in \theta$ (and in $\psi$),
    $(p(a,c,b),p(a,a,b))=(p(a,c,b),b)\in\theta$ and
    $(p(a,b,b),p(a,c,b))=(a,p(a,c,b))\in\psi$ proving
    $(a,b)\in\psi\cdot\theta$.
\end{proof}

\begin{thm}\label{q-ce-2}
    If a variety $\VV$ is coherent, then its congruences permute.
\end{thm}

\begin{proof}
    Set $p(x,y,z)=\omega(z,t_1(y,x,z),...,t_n(y,x,z))$, where
    $\omega$ and $t_1,...,t_n$ are the term operations of Lemma
    \ref{q-ce-1}. Then we have
        $$p(x,z,z)=\omega(z,t_1(z,x,z),...,t_n(z,x,z))=x$$
    and
        \begin{eqnarray*}
        \begin{split}
            p(x,x,z)&=\omega(z,t_1(x,x,z),...,t_n(x,x,z))\\
                    &=\omega(z,z,...,z)\\
                    &=\omega(z,t_1(z,z,z),...,t_n(z,z,z))\\
                    &=z
        \end{split}
        \end{eqnarray*}
    Hence $p(x,y,z)$ is a Mal'cev term of $\VV$ and so by the last
    theorem $\VV$ is congruence permutable.
\end{proof}

Now we can prove this section's main theorem.

\begin{proof}[Proof of Theorem \ref{q-qua-1.1}]
    Let $\A$ have almost minimal spectrum. Then Lemma
    \ref{q-qua-2.3} says that $\{\A\}$ is a direct factor set.
    Theorem \ref{q-qua-2.5} implies that the finite algebras in $\VV(\A)$ have uniform
    congruences and so by Corollary \ref{q-ck-3.11}, $\VV(\A)$ is
    coherent. Finally reference to Theorem \ref{q-ce-2} concludes
    the proof.
\end{proof}

    \section{$\VV(\A)$ is congruence distributive}

In this section we will show that if $\A \times \A$ has a skew
congruence, then $\A$ is prime affine, and if not, then $\A$
generates a congruence distributive equational class.
\begin{defn}
    An algebra is called \emph{congruence distributive} iff it has a distributive congruence lattice.
    We say a variety is congruence distributive iff all of its
    members are.
\end{defn}
For algebras $(\A_i)_{i\in I}$ of the same type there is a natural
embedding of the product of their congruence lattices to the
congruence lattice of their product, namely

$$
    \epsilon:\quad \begin{matrix}\prod_{i\in I}Con(\A_i)&&\To &&
    Con(\prod_{i\in I}\A_i)\\\prod_{i\in I}\theta_i && \mapsto &&
    \{({a},{b})|\forall i\in I
    \,((a_i,b_i)\in\theta_i)\}\end{matrix}
$$

\begin{defn}
    A congruence on a product of algebras of the same type is
    called \emph{factor congruence} iff it is the product of
    congruences on those algebras as defined before; otherwise, it is called
    \emph{skew}.
\end{defn}

We will now concentrate on the case where $\A\times\A$ has a skew
congruence $\theta$; we will follow a result of H. P. Gumm in
\cite{Gum79} to show that in this case $\A$ is affine with respect
to an abelian $p$-group for some prime $p$.

\begin{defn}
    A lattice $\Lfr$ is called \emph{modular} iff it satisfies the
    equation
    $$
        x\cap ((x\cap y)\cup z)=(x\cap y)\cup (x\cap z).
    $$
\end{defn}

\begin{rem}
    It is easy to check that a lattice $\Lfr$ is modular iff in $\Lfr$, $y\leq x$ implies
    $x\cap (y\cup z)=y \cup (x\cap z)$. The 5-element lattice
    $\NN_5$ (over the set $\{0,a,b,c,1\}$ it is defined by $0\leq a\leq
    b\leq 1$ and $0\leq c\leq 1$ and no other elements are
    comparable) is nonmodular:
    $$
        b\cap (a\cup c)=b\neq a=a\cup (b\cap c).
    $$
    Hence, every lattice containing $\NN_5$ is nonmodular.
    Conversely, every nonmodular lattice contains a sublattice
    isomorphic to $\NN_5$: For if $x,y,z$ do not satisfy the modular
    law, then it is easily verified that the identification $(0,a,b,c,1)=(x\cap z, y\cup (x\cap
    z), x\cap (y\cup z),z, y\cup z)$ is such an isomorphism.
    Therefore, if we have a modular lattice $\Lfr$ and two arbitrary elements $b,t\in L$ with $b\leq t$,
    then the length of every path from $b$ to $t$ in the Hasse
    diagram of $\Lfr$ is the same. In that light, the following definition makes
    sense.
\end{rem}

\begin{defn}
    For a cardinal $\alpha$, by $\MM_\alpha$ we understand the modular
    lattice with least and greatest element and $\alpha$ atoms and no other elements.
\end{defn}

The reason why we defined all this is the following:

\begin{lem}\label{q-cp-mod}
    If $\A$ is an algebra with permuting congruences, then
    $Con(\A)$ is modular.
\end{lem}

For the proof of the lemma as well as for later proofs, we need to
recall the following well-known fact.

\begin{lem}
    An algebra $\A$ has permuting congruences iff for all $\psi,
    \theta\in Con(\A)$, $\psi\cup\theta =\psi\cdot\theta$.
\end{lem}

\begin{proof}
    Clearly, $\psi\cdot\theta$ is a congruence on $\A$, the symmetry provided
    by the permutability of $\psi$ and $\theta$. If
    $(a,b)\in\psi$, then since trivially $(b,b)\in\theta$ we have
    that $(a,b)\in\psi\cdot\theta$. Hence, $\psi\cdot\theta
    \geq\psi$ and by the same argument $\psi\cdot\theta
    \geq\theta$. If $(a,b)\in\psi\cdot\theta$, then there exists
    $c\in A$ such that $(a,c)\in\psi$ and $(c,b)\in\theta$.
    Therefore any congruence $\vartheta$ with $\psi\leq\vartheta$
    and $\theta\leq\vartheta$ must by its transitivity contain
    $(a,b)$ so that $\psi\cdot\theta\leq\vartheta$. This concludes
    the proof of one direction; the other one is obvious.
\end{proof}

\begin{proof}[Proof of Lemma \ref{q-cp-mod}]
    Let $\theta_1, \theta_2, \theta_3\in Con(\A)$ with
    $\theta_1\leq \theta_2$. It must be shown that $\theta_2\cap
    (\theta_1\cup\theta_3)=\theta_1\cup (\theta_2\cap\theta_3)$,
    or equivalently, that $\theta_2\cap
    (\theta_1\cup\theta_3)\leq\theta_1\cup
    (\theta_2\cap\theta_3)$. Let $(a,c)\in\theta_2\cap
    (\theta_1\cup\theta_3)$; since $(a,c)\in
    (\theta_1\cup\theta_3)$ and since $\A$ has permuting
    congruences, there exists $b\in A$ such that
    $(a,b)\in\theta_1$ and $(b,c)\in\theta_3$. Moreover,
    $(a,b)\in\theta_2$ as $\theta_1\leq \theta_2$, and
    since also $(a,c)\in\theta_2$, we have that
    $(b,c)\in\theta_2\cdot\theta_2=\theta_2$. Thus,
    $(b,c)\in\theta_2\cap\theta_3$; consequently,
    $(a,c)\in\theta_1\cdot (\theta_2\cap\theta_3)=\theta_1\cup
    (\theta_2\cap\theta_3)$.
\end{proof}

Let us return to our algebra, which we know now has a modular
congruence lattice. Since $\A$ is simple and since
$Con(\A^2)/\ker(\pi_1)\cong Con(\A^2)/\ker(\pi_2)\cong Con(\A)$,
the intervals $[ker \pi_i, 1]$ are equal to $\{ker \pi_i,1\}$. As
$Con(\A\times\A)$ is modular we conclude that the sublattice
generated by $\{ker \pi_1,ker \pi_2,\theta\}$ is isomorphic to
$\MM_3$. Furthermore, the greatest (resp. least) element in
$\MM_3$ coincides with the greatest (resp. least) element in
$Con(\A\times\A)$. We say that $\MM_3$ is a \emph{0-1-sublattice}
of $Con(\A\times\A)$. We summarize: $\A\times\A$ has three
congruences $\theta_1,\theta_2,\theta_2$ satisfying $\theta_i\cdot
\theta_j=1$ and $\theta_i\cap\theta_j=0$ for $1\leq i,j\leq 3$. In
the following, we will investigate an abstraction of this
situation.

Let $S$ be a set, $|S|\geq 4$, and let
$\theta_1,\theta_2,\theta_3$ be equivalence relations on $S$
satisfying $\theta_i\cdot \theta_j=1$ and $\theta_i\cap\theta_j=0$
for $1\leq i,j\leq 3$, $i\neq j$. Then we call the quadruple
$\SSS=(S,\theta_1,\theta_1,\theta_1)$ an \emph{S-3-System}. A
geometrical interpretation of an S-3-system, the so-called
\emph{\"{A}quivalenzklassengeometrie}, will prove useful: Call the
elements of $S$ \emph{points} and the equivalence classes of the
relations \emph{lines}. Two lines are \emph{parallel} iff they are
classes of the same equivalence relation. A point lies on a line
iff it is an element of the line. With these definitions we have:

\begin{lem}\label{q-hpg-2.1}
    The \"{A}quivalenzklassengeometrie of an S-3-system has the following properties:
    \begin{itemize}
    \item[(S1)]{There are three classes of parallel lines.}
    \item[(S2)]{Each point lies on exactly one line of each
    parallel-class.}
    \item[(S3)]{Two non-parallel lines intersect in exactly one point, that is, they have
    exactly one point in common.}
    \end{itemize}
\end{lem}

\begin{proof}
    (S1) and (S2) are trivial. For (S3), let $l_1,l_2$ be two
    non-parallel lines, and assume without loss of generality they
    are equivalence classes of $\theta_1$ and $\theta_2$,
    respectively. Let $x\in l_1$ and $y\in l_2$ be arbitrary
    points. Then since $\theta_1\cdot\theta_2=1$, there is $z\in
    S$ such that $x\theta_1 z$ and $z\theta_2 y$. Hence, $z\in
    l_1\cap l_2$. Suppose there is another $u\in l_1\cap l_2$.
    Then $u\theta_1 z$ and $u\theta_2 z$ so that $u=z$ since
    $\theta_1\cap\theta_2=0$.
\end{proof}

\begin{defn}
    An algebra $\Q=(Q,\cdot)$ with one binary operation $\cdot$
    is called a \emph{quasigroup} iff
    for all $a, c\in Q$ the equations $c\cdot x=a$ and $y\cdot
    c=a$ have unique solutions $x, y\in Q$.
\end{defn}

We will show now that quasigroups give rise to S-3-systems, and
conversely, from S-3-systems we can construct quasigroups. Let
$\Q=(Q,\cdot)$ be a quasigroup. Set $S=Q\times Q$ and define
$\theta_1$, $\theta_2$ and $\theta_2$ on $S$ by
$$
    \begin{matrix}(x,y)\theta_1(x',y')&&\gdw && x=x'\\
                    (x,y)\theta_2(x',y')&&\gdw && y=y'\\
                    (x,y)\theta_3(x',y')&&\gdw && x\cdot y=x'\cdot y'
    \end{matrix}
$$

\begin{lem}\label{q-hpg-quasiS3}
    $(S,\theta_1,\theta_2,\theta_3)$ is an S-3-system.
\end{lem}

\begin{proof}
    Obviously $\theta_i$ is an equivalence relation, $i=1,2,3$.
    Also $\theta_1\cap\theta_2 =0$ and $\theta_1\cdot\theta_2 =1$
    is clear. If $(x,y)(\theta_1\cap\theta_3)(x',y')$, then $x\cdot
    y=x'\cdot y'=x\cdot y'$; since $\Q$ is a quasigroup
    this implies $y=y'$. Thus, $\theta_1\cap\theta_3
    =0$. Let $(x,y)$ and $(x',y')$ be arbitrary elements of $S$.
    There exists $y''\in S$ such that $x\cdot y''=x'\cdot y'$.
    Hence, $(x,y)\theta_1 (x,y'')\theta_3 (x',y')$ so that we have
    $\theta_1\cdot\theta_3 =1$. As the situation with $\theta_2$
    is analogous this concludes the proof.
\end{proof}

For the inverse process start with an S-3-system
$(S,\theta_1,\theta_2,\theta_3)$. Since $\theta_i\cdot\theta_j =
1$ and $\theta_i\cap\theta_j = 0$ for $i\neq j$, we have that
$S\cong S/\theta_i\times S/\theta_j$ as sets (see also Lemma
\ref{q-fp-factoriso}). Hence,
$$S/\theta_1\times S/\theta_2\cong S/\theta_1\times S/\theta_3\cong S/\theta_2\times
S/\theta_3$$ and thus $S/\theta_1\cong S/\theta_2 \cong
S/\theta_3$. Therefore, if we set $Q=|S/\theta_1|$, we know there
are bijections $f_i:S/\theta_i\To Q$, $i=1,2,3$. Define a function
$$g:\quad
    \begin{matrix}
        S&&\To && Q\times Q\\
        s&&\mapsto && (f_1([s]_{\theta_1}),f_2([s]_{\theta_2}))
    \end{matrix}
$$
Then $g$ is a bijection: For it is onto since in the corresponding
\"{A}quivalenzklassengeometrie two non-parallel lines have an
intersection point and it is one-one since this intersection is
unique. Fix an arbitrary element $e\in S$ and an arbitrary element
$1\in Q$. We may suppose we have chosen $f_1,f_2$ such that
$f_1([e]_{\theta_1})=f_2([e]_{\theta_2})=1$. Furthermore, we
suppose that if $a\theta_1 e$, $b\theta_2 e$ and $(a,b)\nin
\theta_3$ then $f_2([a]_{\theta_2})\neq f_1([b]_{\theta_1})$. This
is legitimate since there are $|Q|^2-|Q|$ possibilities to choose
an ordered pair of two distinct equivalence classes in $\theta_3$;
but the other assumptions on $a$ and $b$ already uniquely
determine the representatives of those equivalence classes. Hence,
there are $|Q|^2-|Q|$ possibilities to choose $a$ and $b$
satisfying all conditions which is exactly the number of ordered
pairs of unequal values of $f_1$ and $f_2$. Note that the
assumption implies $g^{-1}(1,x)\theta_3 g^{-1}(x,1)$ for all $x\in
Q$. Define a binary operation $\cdot$ on $Q$ in the following way:
For $x,y$ in $Q$ set $s=g^{-1}(x,y)$. Let $t$ be the intersection
of the $\theta_3$-line through $s$ with the $\theta_2$-line
through $e$; then $z=x\cdot y=g(t)$. More formally,
$$
    x\cdot y= z\gdw g^{-1}(x,y)\theta_3 g^{-1}(z,1).
$$

\begin{defn}
    A \emph{loop} is a quasigroup $\Lfr=(L,\cdot)$ which has an
    element $1\in L$ such that $x\cdot 1=1\cdot x=x$ for all $x\in L$.
\end{defn}

\begin{lem}\label{q-hpg-S3quasi}
    $\Q=(Q,\cdot,1)$ is a loop.
\end{lem}

\begin{proof}
    $x\cdot 1=x$ since trivially $g^{-1}(x,1)\theta_3 g^{-1}(x,1)$;
    $1\cdot x=x$ since we chose $g$ such that $g^{-1}(1,x)\theta_3
    g^{-1}(x,1)$. To find the right-side inverse of an element
    $x\in Q$, let $s$ be the intersection of the $\theta_1$-line
    $f_1^{-1}(x)$ with the $\theta_3$-line through $e$. Then for
    $y=f_2([s]_{\theta_2})$ we have $x\cdot y=1$:
    $g^{-1}(x,y)=s\,\theta_3\, e=g^{-1}(1,1)$. The left-side inverse
    can be found in a similar way.
\end{proof}

In the following, we will identify $S$ with $Q\times Q$. Then one
can easily verify
$$
    \begin{matrix}(x,y)\theta_1(x',y')&&\gdw && x=x'\\
                    (x,y)\theta_2(x',y')&&\gdw && y=y'\\
                    (x,y)\theta_3(x',y')&&\gdw && x\cdot y=x'\cdot y'
    \end{matrix}
$$
Hence, if we start out with an S-3-system, construct a loop as
shown before, and construct from that an S-3-system again, we end
up with the system we started with. We summarize this connection
in the following theorem.

\begin{thm}\label{q-hpg-2.3}
    Let $\SSS =(S,\theta_1,\theta_2,\theta_3)$ be an S-3-system and let $e\in S$ arbitrary. Then
    there exist a loop $\Lfr =(L,\cdot,1)$ and a bijection $g:L\times
    L\To S$ such that $e=g(1,1)$ and for all $x,y,x',y'\in L$ we
    have
    $$
        \begin{matrix}(x,y)\theta_1(x',y')&&\gdw && x=x'\\
                        (x,y)\theta_2(x',y')&&\gdw && y=y'\\
                        (x,y)\theta_3(x',y')&&\gdw && x\cdot y=x'\cdot y'
        \end{matrix}
    $$
    if we identify the elements of $S$ with those of $L\times L$
    via $g$.
\end{thm}

\begin{rem}
    Note that $\theta_3$ need not be a congruence of $\Lfr\times \Lfr$
    whereas $\theta_1$ and $\theta_2$ obviously are.
\end{rem}

Now let us return to our algebra $\A$. Since $\A$ generates a
congruence permutable variety, there exists a Mal'cev term on
$\A$, that is, there exists a ternary term $p$ satisfying the
equations $p(x,x,y)=y$ and $p(x,y,y)=x$. In our case, $p$ is
unique:

\begin{thm}\label{q-hpg-3.2}
    Let $\SSS=(S,\theta_1,\theta_2,\theta_3)$ be an S-3-system and
    let $p$ be a Mal'cev function on $S$ preserving $\theta_1,\theta_2,$ and
    $\theta_3$. Then $p$ is uniquely determined.
\end{thm}

\begin{proof}
    Let $x,y,z\in S$ be given. If $x=y$ or $y=z$ then $p(x,y,z)$
    is determined by the equations of a Mal'cev function. Suppose
    $x\neq y$ and $y\neq z$.
    Assume first that $x$ and $y$ lie on one line $l_1$ and $y$
    and $z$ lie on a line $l_2$ and $l_1\neq l_2$. Since $y$ is
    the intersection of those lines, for some $i\neq k$ we have
    $l_1=[y]_{\theta_i}$ and $l_2=[y]_{\theta_k}$ so that by
    compatibility $p(x,y,z)\theta_i p(x,x,z)=z$ and $p(x,y,z)\theta_k
    p(x,y,y)=x$. Hence, $p(x,y,z)$ is the intersection of the
    $\theta_i$-line through $z$ with the $\theta_k$-line through
    $x$ so that it must be unique.
    In a next step, assume that $x,y,z$ lie on one line $l$ and
    say without loss of generality $l$ is a $\theta_1$-line.
    Denote by $x'$ the intersection of the $\theta_2$-line through
    $y$ with the $\theta_3$-line through $x$. As $x'$ and $y$ lie
    on one line and $y$ and $z$ on another line, we know from the
    first step of the proof that $p(x',y,z)$ is uniquely
    determined. Since $x,y,z$ lie on one $\theta_1$-line we have
    $x=p(x,x,x)\theta_1 p(x,y,z)$; hence, $p(x,y,z)$ lies on $l$
    as well. But $x\theta_3 x'$ implies $p(x,y,z)\theta_3
    p(x',y,z)$. Thus, $p(x,y,z)$ is the intersection of the
    $\theta_3$-line through $p(x',y,z)$ with $l$ and so it is
    uniquely determined.
    To finish the proof, let $x,y,z$ be arbitrary. Consider an
    arbitrary $\theta_1$-line $l_1$ and an arbitrary
    $\theta_2$-line $l_2$. Denote the intersections of the
    $\theta_2$-lines through $x,y,z$ with $l_1$ by $x',y',z'$ and
    the intersections of the $\theta_1$-lines through $x,y,z$ with
    $l_2$ by $x'',y'',z''$. By the second step of our proof,
    $p(x',y',z')$ and $p(x'',y'',z'')$ are uniquely determined.
    Clearly, $p(x,y,z)\theta_2 p(x',y',z')$ and $p(x,y,z)\theta_1
    p(x'',y'',z'')$. Hence, $p(x,y,z)$ is the unique intersection
    of the $\theta_2$-line through $p(x',y',z')$ with the
    $\theta_1$-line through $p(x'',y'',z'')$.
\end{proof}

\begin{cor}\label{q-hpg-3.3}
    Let $\SSS=(S,\theta_1,\theta_2,\theta_3)$ be an S-3-system.
    If $p$ is a Mal'cev function on $S$ compatible with
    $\theta_1,\theta_2,\theta_3$, then it satisfies the equation
    $p(x,y,z)=p(z,y,x)$.
\end{cor}

\begin{proof}
    Set $\tilde{p}(x,y,z)=p(z,y,x)$. Then $\tilde{p}$ is
    apparently a Mal'cev function on $S$ preserving
    $\theta_1,\theta_2,\theta_3$. Therefore, it must equal $p$ and
    so $\tilde{p}(x,y,z)=p(z,y,x)=p(x,y,z)$ for all $x,y,z\in S$.
\end{proof}

\begin{lem}\label{q-hpg-3.4}
    Let $\SSS=(S,\theta_1,\theta_2,\theta_3)$ be an S-3-system.
    If there is a compatible Mal'cev function on $S$, then the loop associated
    with $\SSS$ satisfies:
    $$
        (x_1\cdot
        y_1=x_2\cdot y_2 \wedge x_1\cdot y_3=x_2\cdot y_4 \wedge
        x_3\cdot y_1=x_4\cdot y_2)\To x_3\cdot y_3=x_4\cdot y_4.
    $$
\end{lem}

\begin{proof}
    Recall that in terms of the S-3-system and its congruence $\theta_3$ our hypothesis says
    $$
        (x_1,y_1)\theta_3 (x_2,y_2) \wedge (x_1,y_3)\theta_3
        (x_2,y_4)\wedge (x_3,y_1)\theta_3 (x_4,y_2).
    $$
    Since $p$ is compatible with $\theta_1$ and $\theta_2$, it satisfies the Mal'cev
    conditions componentwise. Hence,
    $(x_3,y_3)=p((x_1,y_3),(x_1,y_1),(x_3,y_1))\theta_3
    p((x_2,y_4),(x_2,y_2),(x_4,y_2))=(x_4,y_4)$.
\end{proof}

\begin{lem}\label{q-hpg-3.5}
    Let $\SSS$ be an S-3-system with a compatible Mal'cev function. Then the loop $\Lfr$ associated with $\SSS$ is
    associative, i.e. a group.
\end{lem}

\begin{proof}
    The previous lemma applies; so to check the associative law
    for arbitrary $x,y,z\in L$ set $x_1=y_2=1$, $x_2=y$, $x_3=x$,
    $x_4=x\cdot y$, $y_1=y$, $y_3=y\cdot z$, $y_4=z$. Then all
    hypotheses of the lemma are satisfied and it yields
    $x\cdot(y\cdot z)=(x\cdot y)\cdot z$.
\end{proof}

So we know that if we have a Mal'cev operation compatible with an
S-3-system, the associated loop is in fact a group. We will show
now that this group is even abelian.

\begin{lem}\label{q-hpg-3.6}
    Let the S-3-system $\SSS$ admit the Mal'cev function $p$. Then we can
    calculate $p$ by
    $$
        p((x_1,y_1),(x_2,y_2),(x_3,y_3))=(x_1\cdot x_2^{-1}\cdot
        x_3,y_1\cdot y_2^{-1}\cdot y_3).
    $$
\end{lem}

\begin{proof}
    We will calculate $p(x,y,z)$ following the construction of
    that point in the proof of Theorem \ref{q-hpg-3.2}. Set
    $x=(x_1,y_1)$, $y=(x_2,y_2)$ and $z=(x_3,y_3)$. Let $l_1$ be
    the $\theta_1$-line and $l_2$ be the
    $\theta_2$-line through $x$. Then, using the same notation
    as in that proof, we have $x'=(x_1,y_1)$, $y'=(x_1,y_2)$,
    $z'=(x_1,y_3)$ and $x''=(x_1,y_1)$, $y''=(x_2,y_1)$,
    $z''=(x_3,y_1)$. Thus, if we write $p$ also for the functions
    $p$ induces on the components,
    \begin{eqnarray*}
    \begin{split}
        p(x',y',z')&=p((x_1,y_1),(x_1,y_2),(x_1,y_3))\\
        &=(p(x_1,x_1,x_1),p(y_1,y_2,y_3))\\
        &=:(x_1,\bar{p})
    \end{split}
    \end{eqnarray*}
    and similarly
    $p(x'',y'',z'')=(p(x_1,x_2,x_3),y_1)=:(\bar{\bar{p}},y_1)$.
    Since $p(x,y,z)\theta_2 p(x',y',z')$ and $p(x,y,z)\theta_1
    p(x'',y'',z'')$ we get $p(x,y,z)=(\bar{\bar{p}},\bar{p})$. For
    the computation of $p(x',y',z')$ we can use the second step in
    the proof of \ref{q-hpg-3.2} since $x',y',z'$ lie on one line $l_1$: Let
    $s=(u,y_2)$ be the intersection of the $\theta_3$-line through
    $x'$ with the $\theta_2$-line through $y'$; the definition of
    $\theta_3$ immediately yields the equation $x_1\cdot y_1=u\cdot y_2$.
     Let $t=(u,y_3)$ be the intersection of the $\theta_1$-line through $s$ with
    the $\theta_2$-line through $z'$. Then, as
    $p(x',y',z')=(x_1,\bar{p})$ is the intersection of the
    $\theta_3$-line through $t$ with $l_1$, we get $x_1\cdot
    \bar{p}=u\cdot y_3$ which we can solve to
    $$
        \bar{p}=x_1^{-1}\cdot u\cdot y_3\\
        =x_1^{-1}\cdot x_1\cdot y_1\cdot y_2^{-1}\cdot y_3\\
        =y_1\cdot y_2^{-1}\cdot y_3.
    $$
    Similarly, $\bar{\bar{p}}=x_1\cdot x_2^{-1}\cdot x_3$
    and so the proof is complete.
\end{proof}

\begin{cor}\label{q-hpg-3.7}
    If $\SSS=(S,\theta_1,\theta_2,\theta_3)$ is an S-3-system
    which allows a Mal'cev function on $S$, then the associated group $\G$ is abelian.
\end{cor}

\begin{proof}
    Combining Lemma \ref{q-hpg-3.6} with Corollary \ref{q-hpg-3.3}
    yields $(x_1\cdot x_2^{-1}\cdot x_3,y_1\cdot y_2^{-1}\cdot
    y_3)=(x_3\cdot x_2^{-1}\cdot x_1,y_1\cdot y_2^{-1}\cdot
    y_3)$ for all $x_i,y_i\in G$, $i=1,2,3$. Therefore, if we set
    $x_2=1$, we have $x_1\cdot x_3=x_3\cdot x_1$.
\end{proof}

Since we know now that we are dealing an abelian group, we will
change our notation to an additive one; that is, we will write $+$
for the binary group operation and $0$ for the neutral element.
Furthermore, we will identify the base set $S$ of an S-3-system
that admits a Mal'cev operation with $G\times G$, where
$\G=(G,+,0,-)$ is the associated abelian group. Observe that for
the Mal'cev operation $p$ on $G\times G$ we have $p(x,y,z)=x-y+z$,
where $+$ is calculated componentwise.

Now let $f:{(G\times G)}^n\To G$ be an $n$-ary function on
$G\times G$ compatible with $\theta_1,\theta_2,\theta_3$. Since
$f$ is compatible with $\theta_1$ and $\theta_2$, it is the
product of two mappings $f_1, f_2:G^n\To G$, i.e.
$$
    f((x_1,y_1),...,(x_n,y_n))=(f_1(x_1,...,x_n),f_2(y_1,...,y_n)).
$$
\begin{lem}\label{q-hpg-3.8}
    For $x,y,x',y'\in G^n$ the following
    holds:
    \begin{enumerate}
        \item[(i)]{$x+y=x'+y'\To f_1(x)+f_2(y)=f_1(x')+f_2(y')$}
        \item[(ii)]{$f_1(x)+f_2(0)=f_1(0)+f_2(x)$}
        \item[(iii)]{$f_1(x)+f_2(y)=f_1(x+y)+f_2(0)$}
    \end{enumerate}
\end{lem}

\begin{proof}
    The hypothesis of (i) says that for $1\leq k\leq n$,
    $x_k+y_k=x_k'+y_k'$. Thus, by definition of $\theta_3$,
    $(x_k,y_k)\theta_3 (x_k',y_k')$ for $1\leq k\leq n$ so that by
    the compatibility of $f$ with $\theta_3$,
    $f((x_1,y_1),...,(x_n,y_n))\theta_3
    f((x_1',y_1'),...,(x_n',y_n'))$. Hence,
    $(f_1(x),f_2(y))\theta_3
    (f_1(x'),f_2(y'))$ and so
    $f_1(x)+f_2(y)=f_1(x')+f_2(y')$. (ii) is trivial with (i). For
    (iii) set $x'=x+y$ and $y'=0$ and apply (i).
\end{proof}

\begin{defn}
    Let $A$ be a set and $f:A^n\To A$ be an $n$-ary function on
    $A$. If a binary operation
    $+$ can be defined on $A$ such that $(A,+)$ is an abelian
    group and for all ${x},{y}\in A^n$ we have
    $$
        f(x)+f(y)=f(x+y)+f(0),
    $$
    then we say that $f$ is \emph{affine} with respect to $(A,+)$.
    An algebra $\A$ is \emph{affine} iff every fundamental
    operation is affine with respect to the same abelian group
    over $A$.
\end{defn}

\begin{rem}
    It is obvious that an algebra $\A$ is affine if and only if there is an
    affine relation $\rho\subseteq A^4$ as defined in the first
    chapter which is preserved by all operations of $\A$.
\end{rem}

\begin{lem}\label{q-hpg-4.3}
    Let $\SSS=(S,\theta_1,\theta_2,\theta_3)$ be an S-3-system with a compatible Mal'cev
    function. Then every mapping on $S$ which is compatible with $\theta_1$, $\theta_2$, $\theta_3$ is affine with
    respect to $\G\times \G$, where $\G$ is the abelian group associated with $\SSS$.
\end{lem}

\begin{proof}
    Suppose $f:S^n\To S$ is compatible with $\theta_1$, $\theta_2$, $\theta_3$. Then for
    $x,y\in S$, if we write $x=(x',x''),y=(y',y'')$, we can
    compute by the previous lemma:
    \begin{eqnarray*}
    \begin{split}
        f(x)+f(y)&=(f_1(x'),f_2(x''))+(f_1(y'),f_2(y''))\\
        &=(f_1(x')+f_1(y'),f_2(x'')+f_2(y''))\\
        &=(f_1(x')+f_1(0)+f_2(y')-f_2(0),f_2(x'')+f_2(0)+f_1(y'')-f_1(0))\\
        &=(f_1(x'+y')+f_1(0),f_2(x''+y'')+f_2(0))\\
        &=(f_1(x'+y'),f_2(x''+y''))+(f_1(0),f_2(0))\\
        &=f(x+y)+f(0).
    \end{split}
    \end{eqnarray*}
\end{proof}

In terms of algebras we have established:

\begin{thm}\label{q-hpg-4.6}
    Let $\A$ be an algebra in a congruence permutable variety and
    let $p$ be a Mal'cev term of $\A$. If $\MM_3$ is a
    0-1-sublattice of $Con(\A)$, then there is an abelian group
    $\G=(G,+,0)$ such that $\A$ is isomorphic as sets to
    $\G\times\G$ and such that the following holds:
    $p(x,y,z)=x-y+z$ and every term operation $f$ is affine with
    respect to $\G\times\G$ and of the form $f_1\times f_2$ where
    $f_1,f_2: G^n\To G$ if $f$ is $n$-ary.
\end{thm}

The connection to skew congruences is the following:

\begin{thm}\label{q-hpg-4.7}
    Let $\A$ be a simple algebra in a congruence permutable variety. If
    $\A\times\A$ has a skew congruence, then
    $\A$ is affine.
\end{thm}

\begin{proof}
    If $\A\times\A$ has a skew congruence $\theta$, then $\theta$ is by our
    discussion at the beginning of this section
    a complement of $ker(\pi_1^2)$ and of $ker(\pi_2^2)$.
    Hence, the last theorem applies to $\A\mult\A$. But it follows from the
    construction of $\G\mult\G$ by means of the congruences $ker(\pi_1^2)$ and
    $ker(\pi_2^2)$ that the canonical coordinate representation of an element of
    $\A\mult\A$ is exactly the same as the representation with
    respect to the factorization $\G\mult\G$. We conclude that every term operation on $\A$
    is affine with respect to $\G$ as the corresponding term
    operation on $\A\mult\A$ is affine with respect to
    $\G\mult\G$. Thus, $\A$ is affine with respect to $\G$.
\end{proof}

Our next goal is to show that for a simple non-trivial affine
algebra the underlying abelian group is in fact a $p$-group for
some prime $p$. Define for arbitrary $n\geq 1$ a binary relation
$\Delta_n$ on $\A$ by
$$
    x\Delta_n y\gdw n(x-y)=0.
$$

\begin{lem}\label{q-hpg-4.9}
    Let $\A$ be an affine algebra. Then $\Delta_n$ is a
    congruence on $\A$ for all $n\geq 1$.
\end{lem}

\begin{proof}
    For an arbitrary $k$-ary fundamental operation of $\A$, if
    $x_i\Delta_n y_i$ for $1\leq i\leq k$, we have
    \begin{eqnarray*}
    \begin{split}
        n(f(x_1,...,x_k)-f(y_1,...,y_k))
        &=n(f(x_1-y_1,...,x_k-y_k)-f(0,...,0))\\
        &=n f(x_1-y_1,...,x_k-y_k)-n f(0,...,0)\\
        &=f(n(x_1-y_1),...,n(x_k-y_k))-f(0,...,0)\\
        &=f(0,...,0)-f(0,...,0)\\
        &=0.
    \end{split}
    \end{eqnarray*}
\end{proof}

\begin{lem}\label{q-hpg-5.5}
    Let $\A$ be a simple non-trivial affine algebra.
    Then the underlying group $\G$ is either torsion-free or a
    $p$-group for some prime $p$.
\end{lem}

\begin{proof}
    Suppose $\G$ is not
    torsion-free. Then there exists a smallest positive number $p$ such
    that $p \, a=0$ for some $a\in G$, $a\neq 0$. Obviously $p$ is a
    prime. Consider the congruence $\Delta_p$ of the previous lemma.
    Since $a\Delta_p 0$ and $a\neq 0$, $\Delta_p$ must equal $1\in Con(\A)$
    as $\A$ is simple. But then it readily follows that $p \, a=0$ for all $a\in G$.
\end{proof}

\begin{thm}\label{q-hpg-6.1}
    Let $\A$ be a simple non-trivial algebra in a permutable variety. If
    $\A\mult\A$ has a skew congruence, then $\A$ is affine with
    respect to a torsion-free abelian group or with respect to an
    abelian $p$-group.
\end{thm}

\begin{proof}
    Follows from Theorem \ref{q-hpg-4.7} and the previous lemma.
\end{proof}

The following theorem (see H. Werner \cite{Wer74} for a slightly
more general result) tells us that if we have a skew congruence on
a higher power of $\A$ we can still use our results and $\A$ must
be prime affine as well.

\begin{thm}\label{q-wer-4}
    Let $\A_1,...,\A_{n}$ be algebras in a permutable variety. Then
    $\A_1\times ... \times \A_n$ has a skew congruence if and only if $\A_i\times \A_j$ has a skew
    congruence for some $1\leq i,j\leq n$ with $i\neq j$.
\end{thm}

Before we can prove the theorem, we need a definition. For
$\theta\in Con(\A\times\B)$ and $a\in A$, we define an equivalence
relation $\theta_a$ on $\B$ by
$$
    \theta_a=\{(b_1,b_2)|(a,b_1)\theta (a,b_2)\}.
$$
In a congruence permutable variety, $\theta_a$ is a congruence
relation and independent of the choice of $a$:
\begin{lem}\label{q-wer-2}
    Let $\A, \B$ be algebras in a variety with permutable
    congruences and let $\theta\in Con(\A\times\B)$. Then for all
    $a,a'\in A$ and all $b\in B$
        \begin{enumerate}
            \item[(i)]{ $\theta_a=\theta_{a'}$ }
            \item[(ii)]{$\theta_a\in Con(\B)$}
            \item[(iii)]{$\theta_b\times\theta_a\subseteq\theta $}
        \end{enumerate}
\end{lem}

\begin{proof}
    To prove (i), denote the Mal'cev term of the variety by $p$. We have to show
    for arbitrary $a,a'\in A$ and $b,b'\in B$ that
    $(a,b)\theta (a,b')$ implies $(a',b)\theta
    (a',b')$. But as trivially $(a',b')\theta
    (a',b')$ and $(a,b)\theta (a,b)$, application of $p$
    to the three elements of $\theta$ yields
    $$
        (a',b')=(p(a',a,a),p(b',b,b))\theta
        (p(a',a,a),p(b',b',b))=(a',b)
    $$
    so that (i) is indeed true.
    Now (ii) is an immediate consequence of (i) since if $f$ is an
    $n$-ary operation of $\B$ and $(b_i,b'_i)\in\theta_a$, $1\leq i\leq n$, then
    $(f(b_1,...,b_n),f(b'_1,...,b'_n))\in\theta_{f(a,...,a)}=\theta_a$.
    For (iii), let $((c,d),(c',d'))\in\theta_b\times\theta_a$ be
    given. By (i), we can replace $b$ by $d$ and $a$ by $c'$. Then
    we get by the definitions of $\theta_d$ and $\theta_{c'}$ that $(c,d)\theta (c',d)$ and $(c',d)\theta
    (c',d')$ so that $((c,d),(c',d'))\in\theta$.
\end{proof}

\begin{cor}\label{q-wer-3}
    Let $\A, \B$ be algebras in a variety with permutable
    congruences and let $\theta\in Con(\A\times\B)$. Then the
    following conditions are equivalent:
    \begin{enumerate}
        \item[(i)]{$\theta$ is not skew.}
        \item[(ii)]{$\theta_b\times\theta_a =\theta $.}
        \item[(iii)]{$(a,b)\theta (a',b')$ implies $(a,b)\theta (a,b')$.}
        \item[(iii*)]{$(a,b)\theta (a',b')$ implies $(a,b)\theta (a',b)$.}
    \end{enumerate}
\end{cor}

\begin{proof}
    (i) $\Rightarrow$ (ii): Let $\theta =\theta_1\times\theta_2$. If $(a,a')\in\theta_1$, then
    since trivially $(b,b)\in\theta_2$ we have
    $((a,b),(a',b))\in\theta_1\times\theta_2=\theta$ so that
    $(a,a')\in\theta_b$. As the situation for $\theta_2$ is the
    same we conclude $\theta_b\times\theta_a \supseteq\theta$ and
    so $\theta_b\times\theta_a =\theta$.
    (ii) $\Rightarrow$ (i) and (ii) $\Rightarrow$ (iii) are
    trivial. To see that (iii) implies (ii), let $(a,b)\theta
    (a',b')$. By hypothesis, $(a,b)\theta (a,b')$ so that
    $(b,b')\in\theta_a$. But the transitivity of $\theta$ implies $(a',b')\theta
    (a,b')$. Hence, $(a,a')\in\theta_{b'}=\theta_b$.
\end{proof}

Now we are ready to prove the theorem.

\begin{proof}[Proof of Theorem \ref{q-wer-4}]
    We have to show that if $\A_1\times ... \times \A_n$ has a skew congruence then
    $\A_i\times \A_j$ has a skew congruence for some $1\leq i,j\leq n$ with $i\neq
    j$ as the other direction is obvious. To achieve this we will
    prove for algebras $\A$, $\B$, $\C$ of our variety
    that if $\A\times\B, \A\times\C,$ and $\B\times\C$ have only
    factor congruences, then $\A\times\B\times\C$ can have only
    factor congruences as well. The rest will follow by induction.
    Let $\theta\in Con(\A\times\B\times\C)$. Define $\phi_A=\theta_{(b,c)}, \,\phi_B=\theta_{(a,c)}$ and
    $\phi_C=\theta_{(a,b)}$. By Lemma \ref{q-wer-2} (iii),
    $\phi_A\times\phi_B\times\phi_C\subseteq\theta$.

    \textbf{Claim.} $(a,b,c)\theta (a',b',c')\wedge a\phi_A a'\To
    b\phi_B b'\wedge c\phi_C c'$.

    \emph{Proof.} Since
    $\phi_A\times\phi_B\times\phi_C\subseteq\theta$,
    $(a,b,c)\theta (a',b',c')\theta (a,b',c')$. Thus,
    $(b,c)\theta_a (b',c')$. As $\theta_a$ is not skew,
    $(b,c)\theta_a (b,c')$ by Corollary \ref{q-wer-3} (iii).
    Therefore, $(a,b,c)\theta (a,b,c')$ and so $(c,c')\in\phi_C$.
    Hence, $(a,b,c)\theta (a',b',c') \theta (a,b',c)$ so that
    $(b,b')\in\phi_B$.

    Denote by $\pi$ the projection of $A\times B\times C$ onto
    $A\times B$. Then $\tilde{\theta}=\pi [\theta]$ is a congruence on $\A\times
    \B$: The only property which is not obvious is the
    transitivity of $\tilde{\theta}$. Assume $(a,b)\tilde{\theta}(a',b')$ and
    $(a',b')\tilde{\theta}(a'',b'')$; then there exist $c,
    c', d, d'\in C$ such that $(a,b,c)\theta (a',b',c')$ and
    $(a',b',d)\theta (a'',b'',d')$. If $p$ is a Mal'cev term of
    the variety, then since trivially $(a'',b'',d')\theta (a'',b'',d')$ we
    get
    $$
        (p(a,a'',a''),p(b,b'',b''),p(c,d',d'))\theta
        (p(a',a',a''),p(b',b',b''),p(c',d,d'))
    $$
    This yields $(a,b,c)\theta (a'',b'',p(c',d,d'))$ and so $(a,b)\tilde{\theta}(a'',b'')$. Hence,
    $\tilde{\theta}$ is a congruence
     which is not skew by assumption. Now if $(a,b,c)\theta
    (a',b',c')$, then $(a,b)\tilde{\theta}(a',b')$ so that also
    $(a,b)\tilde{\theta}(a,b')$. Hence, there exist $d, d'\in C$
    such that $(a,b,d)\theta (a,b',d')$. Since trivially
    $(a,a)\in\phi_A$, by our claim we get that $(b,b')\in\phi_B$
    and again by the claim that $(a,a')\in\phi_A$ and
    $(c,c')\in\phi_C$. Therefore,
    $\theta\subseteq\phi_A\times\phi_B\times\phi_C$ which is
    exactly what we wanted to show.
\end{proof}

\begin{thm}\label{cl-3}
    Let $\A$ be a finite simple algebra in a permutable variety.
    Then $\A$ is either prime affine or its powers have only
    (trivial) factor congruences.
\end{thm}

\begin{proof}
    If $\A\times\A$ has a skew congruence, then by Theorem
    \ref{q-hpg-6.1} and the fact that $\A$ is finite we get that
    $\A$ is prime affine. Otherwise Theorem \ref{q-wer-4} applies and all
    powers of $\A$ have only factor congruences.
\end{proof}

To use this result in the proof of primality of our algebra $\A$,
recall that since $\{\A\}$ is a direct factor set,
$\VV(\A)=HP(\A)$. Now by the last theorem, if a power of $\A$ has
a skew congruence, then $\A$ is prime affine which is forbidden by
Rosenberg's list. Thus, all powers of $\A$ have only trivial
factor congruences: Since $\A$ is simple, they are products of $0$
and $1\in Con(\A)$. But this implies also that up to isomorphism
$\VV(\A)=P(\A)$ and so the variety generated by $\A$ is obviously
congruence distributive.

    \section{$\A$ is primal}
We will use a special case of a result on \emph{semi-primal}
algebras by A. Foster and A. Pixley in \cite{FP64} to show that
our hypotheses on $\A$ imply it is primal.

\begin{thm}\label{q-fp3.1}
    Let $\A$ be an algebra, $1<|\A|<\aleph_0$. Assume also that $\A$
    is simple, has no proper subalgebras and no proper automorphisms
    and that it generates a congruence permutable and congruence
    distributive variety. Then $\A$ is primal.
\end{thm}

\begin{rem}
    Non-trivial congruences are obviously exactly class three in $RBL$;
    proper subalgebras are central relations and thus in class five of
    $RBL$. Moreover, assume $\phi$ is a proper automorphism of $\A$. Now either
    all cycles of $\phi$ have the same length $n$; then for any
    prime factor $p$ of $n$, $\phi^{\frac{n}{p}}$ has only cycles of the same prime length $p$ and
    hence its graph belongs to class two of $RBL$. If there are
    cycles of different length, then denote the length of the shortest cycle by $n$; clearly, $\phi^n$ is not the
    identity but has at least one fixed point. But the set of all fixed points of
    $\phi^n$ is a proper subalgebra of $\A$ and therefore in class five of $RBL$. Hence our algebra
    fulfills the hypotheses and is primal.
\end{rem}

We will need a couple of rather basic lemmas on subdirect
products; a good standard reference with more details on the
subject is \cite{Bir48}.

\begin{lem}\label{q-fp2.1}
    An algebra $\A$ is isomorphic to a subdirect product of the
    algebras $\{\A_i | i \in I\}$ iff for each $i \in I$ there is a
    homomorphism $h_i$ from $\A$ onto $\A_i$ such that
    $\bigwedge_{i\in I}ker(h_i)=0$
\end{lem}

\begin{proof}
    For one implication, consider as homomorphisms the projections
    $\pi_i$ of elements of the subdirect product onto the $i$-th
    coordinate. The $\pi_i$ are obviously homomorphisms and one can
    easily verify the assertion on the kernels. Conversely, consider the
    mapping
    $$
        \phi:\quad
        \begin{matrix}
            \A&&\rightarrow&& \prod_{i\in I}\A_i\\
            a &&\mapsto &&(h_i(a))_{i\in I}\\
        \end{matrix}
    $$
    Then $\phi$ is clearly a homomorphism
    onto a subdirect product of $\{\A_i | i \in I\}$ which is one-one
    and hence an isomorphism as $\bigwedge_{i\in I}ker(h_i)=0$.
\end{proof}

\begin{lem}\label{q-fp2.2}
    Let $\A$ and $\B$ be algebras of the same type and let
    $h:\A\rightarrow \B$ be a homomorphism. Then $Con(h(\A))$ is
    isomorphic to the sublattice of $Con(\A)$ over the set $\{\theta \in
    Con(\A) | \theta \geq ker(h) \}$.
\end{lem}

\begin{proof}
    By the Homomorphism Theorem, $\A/ker(h)\cong h(\A)$ so that
    $Con(\A/ker(h))\cong Con(h(\A))$. But as one can easily see, the
    sublattice of $Con(\A)$ over the set $\{\theta \in Con(\A) | \theta
    \geq ker(h) \}$ is isomorphic to $Con(\A/ker(h))$ and the lemma
    follows.
\end{proof}

\begin{lem}\label{q-fp-factoriso}
    Let $\theta_1, \, \theta_2$ be permutable congruences on an
    algebra $\A$ satisfying $\theta_1 \cap \theta_2 = 0 $ and
    $\theta_1 \cup \theta_2=1$. Then $\A \cong \A/\theta_1 \times
    \A/\theta_2$.
\end{lem}

\begin{proof}
    Consider the mapping
    $$
        \phi:\quad
        \begin{matrix}
            \A && \rightarrow && \A/\theta_1 \times \A/\theta_2\\
            a  && \mapsto     && ([a]_{\theta_1},[a]_{\theta_2})\\
        \end{matrix}
    $$
    Apparently, $\phi$ is a homomorphism and since $\theta_1 \cap \theta_2 = 0
    $ it is one-one. Now let any $([a]_{\theta_1},[b]_{\theta_2})
    \in \A/\theta_1 \times \A/\theta_2$ be given. As a consequence
    of the permutability of the two congruences, $\theta_1\cdot
    \theta_2 = \theta_1 \cup \theta_2 = 1$ so that $(a,b) \in
    \theta_1 \cdot\theta_2$. Hence, there exists $z$ such that
    $(a,z) \in \theta_1$ and $(z,b) \in \theta_2$. Thus,
    $$([a]_{\theta_1},[b]_{\theta_2})=([z]_{\theta_1},[z]_{\theta_2})=\phi(z)$$
    and $\phi$ is onto.
\end{proof}

\begin{lem}\label{q-fp2.3}
    An algebra $\A$ is isomorphic to the direct product of algebras
    $\{\A_1,...,\A_n\}$ iff
    \begin{enumerate}
        \item{for each $1 \leq i\leq n$ there is a homomorphism $h_i$ from
        $\A$ onto $\A_i$ such that}
        \item{the set of all intersections of the kernels of the $h_i$ consists of pairwise permutable congruence relations and}
        \item{$ker(h_1)\cap ...\cap ker(h_n)=0$ and for $2\leq i\leq n$,
        $(ker(h_1)\cap ...\cap ker(h_{i-1}))\cup ker(h_i)= 1$}
    \end{enumerate}
\end{lem}

\begin{proof}
    If $\A \cong \A_1 \times ... \times \A_n$, consider as
    homomorphisms again the projections; the asserted properties of their
    kernels are easy to verify.
    Conversely, by the previous lemma we have that $\A \cong \B_n
    \times \A_n$, where $\B_n = \A/(\theta_1 \cap ... \cap
    \theta_{n-1})$. A straightforward induction finally shows
    that $\B_n\cong \A_1 \times ... \times \A_{n-1}$ and completes the proof.
\end{proof}

\begin{thm}\label{q-fp2.4}
    Let $\A$ be an algebra isomorphic to a subdirect product of
    finitely many simple algebras $\A_1,...,\A_n$. If the congruences
    of $\A$ permute, then $\A$ is isomorphic to the direct product of
    a subset of the $\A_1,...\A_n$.
\end{thm}

\begin{proof}
    Let $h_1,...,h_n$ be the homomorphisms given by Lemma
    \ref{q-fp2.1} and denote their kernels by
    $\theta_1,...,\theta_n$. Lemma \ref{q-fp2.2} together with the
    simplicity of the $\A_i$ implies that the $\theta_i$ are maximal
    in $Con(\A)$ (if we assume all the $\A_i$ are non-trivial; if not,
    we simply leave the trivial ones away; if all are trivial then the
    theorem is as well). Since $\theta_1 \cap ... \cap
    \theta_n=0$, we can extract a minimal subset of the
    $\theta_i$ having the same property. Assume without loss of
    generality the first $k$ congruences form such a subset, that is, $\theta_1 \cap ... \cap
    \theta_k=0$.
    Trivially,
    \[ (\theta_1 \cap ... \cap \theta_{i-1})\cup \theta_i \geq
    \theta_i, \quad 2\leq i\leq k \] and thus by the maximality of
    $\theta_i$, $(\theta_1 \cap ... \cap \theta_{i-1})\cup \theta_i$
    must be equal to $\theta_i$ or to $1$. But if it was equal to
    $\theta_i$, we could conclude that $(\theta_1 \cap ... \cap
    \theta_{i-1}) \leq \theta_i$ and then
    \[\theta_1 \cap ... \cap \theta_{i-1}\cap \theta_{i+1} \cap ...
    \cap \theta_k=0,\] contradicting the minimality of the set
    $\{\theta_1,...,\theta_k\}$. Consequently, $(\theta_1 \cap ...
    \cap \theta_{i-1})\cup \theta_i= 1$ for $2\leq i\leq k$ and so Lemma
    \ref{q-fp2.3} completes the proof.
\end{proof}

The following lemma from lattice theory will help us using the
congruence distributivity of $\VV(\A)$ in our proof. A meet in a
lattice is \textit{irredundant} iff it cannot be written as a meet
of a subset of its elements. An element in a lattice is
\textit{meet irreducible} iff it is not the meet of two elements
not equal to itself.

\begin{lem}\label{q-fp-uniqueFact}
    In a distributive lattice $\Lfr$, the representation of an element as
    an irredundant meet of meet-irreducible elements is unique (and
    dually).
\end{lem}

\begin{proof}
    Let $a$ be an element of $\Lfr$ such that
    $$
        a= x_1 \cap ... \cap x_r = y_1 \cap ... \cap y_s.
    $$
    Then for a
    fixed $x_i$ we can observe the following: Clearly, $x_i \geq y_1
    \cap ... \cap y_s$. Thus $x_i = x_i \cup (y_1 \cap ... \cap y_s) =
    (x_i \cup y_1) \cap ... \cap (x_i \cup y_s)$ by the distributivity
    of $\Lfr$. But if $x_i$ is meet irreducible, the above representation
    yields $x_i \cup y_j = x_i$ and hence $y_j \leq x_i$ for some $j$.
    Similarly for $y_j$ we have $x_k \leq y_j$ for some $k$ so that
    $x_k \leq x_i$ and therefore $x_k = x_i = y_j$ because the
    representation was assumed to be non-redundant. So the $x_i$ and
    $y_j$ are equal in pairs, $r=s$ and the representation is indeed
    unique.
\end{proof}

We shall now obtain some results concerning the structure of free
algebras. Let $S$ be an algebra type. For a set of $S$-identities
$\Psi$ we let $\F_k(\Psi)$ denote the free algebra with $k$
generators determined by $\Psi$. If $\A$ is an algebra then
$\Sigma(\A)$ will denote the equations satisfied by $\A$. Finally,
$\F_k(\A)$ is short for $\F_k(\Sigma(\A))$. Recall the following
important fact:

\begin{lem}\label{q-fp4.1}
    $\F_k(\A) \in \VV(\A)$.
\end{lem}

Let $\A$ be a non-trivial finite algebra of order $n$ and let
$G=\{\xi_1,...,\xi_k\}$ be a set of $k$ indeterminates. Then there
exist $n^k$ functions $e_1,...,e_{n^k}$ from $G$ to $\A$. Each of
the $e_i$ induces a subuniverse $S_i$ of $\A$; but as we assume
that $\A$ has no proper subalgebras, all of the $S_i$ are equal to
$A$.

Construct $\F_k(\A)$ over $G$. All the $e_i$ induce in a canonical
way a homomorphism $h_i$ from $\F_k(\A)$ onto $\A$: For a class
$\Phi$ of equivalent expressions in $\F_k(\A)$,
$$
    h_i(\Phi)=\phi(e_i(\xi_1),...,e_i(\xi_n)),
$$
where $\phi$ is an arbitrary $S$-expression in $\Phi$. The
function is well-defined since for $\phi_1, \phi_2$ in $\Phi$,
$\phi_1(e_i(\xi_1),...,e_i(\xi_n))=\phi_2(e_i(\xi_1),...,e_i(\xi_k))$.
It is easily seen that $h_i$ is indeed a homomorphism onto $\A$.
Hence, $\F_k(\A)/ker(h_i)\cong \A$.

If $\Phi_1 \equiv  \Phi_2 \,(\bigwedge_{1 \leq i \leq n^k}
\ker(h_i))$ then for all $1 \leq i \leq n^k$ we have that
$h_i(\Phi_1)=h_i(\Phi_2)$ and hence for all $\phi_1 \in \Phi_1$
and all $\phi_2 \in \Phi_2$,
$\phi_1(e_i(\xi_1),...,e_i(\xi_n))=\phi_2(e_i(\xi_1),...,e_i(\xi_k))$.
Since this holds for all possible mappings $e_i$ into $\A$,
$\phi_1=\phi_2$ must be an identity of $\Sigma(\A)$ and thus
$\Phi_1=\Phi_2$. Therefore,
\begin{equation}\label{q-fp-representation}
    \bigwedge_{1 \leq i\leq n^k} \ker(h_i)=0,
\end{equation}
and using Lemma \ref{q-fp2.1} we conclude:

\begin{lem}\label{q-fp4.2}
    Let $\A$ be a non-trivial finite algebra of order $n$ having no
    proper subalgebras. Then $\F_k(\A)$ is isomorphic to a subdirect product of $n^k$ copies of
    $\A$ via the mapping
    $$
    \begin{matrix}
        \F_k(\A) && \To     && \A^{n^k} \\
        \Phi     && \mapsto &&
        [\phi(e_1(\xi_1),...,e_1(\xi_k)),...,\phi(e_{n^k}(\xi_1),...,e_{n^k}(\xi_k))]
    \end{matrix}
    $$
    where $\phi$ is an arbitrary term in $\Phi$.
\end{lem}

With the additional assumption that $\A$ is simple and generates a
congruence permutable variety, Theorem \ref{q-fp2.4} and Lemma
\ref{q-fp4.1} imply

\begin{lem}\label{q-fp4.3}
    Let $\A$ be a simple non-trivial finite algebra of order $n$ having no
    proper subalgebras, and which generates a congruence permutable
    equational class. Then there exists $1\leq r \leq n^k$ such that
    $\F_k(\A)$ is isomorphic to $\A^r$.
\end{lem}

If all of the assumptions on $\A$ in Theorem \ref{q-fp3.1} hold,
then all of the factors occur in the representation of the free
algebra $\F_k(\A)$.

\begin{lem}
    If $\A$ is an algebra satisfying the assumptions of Theorem
    \ref{q-fp3.1}, then $\F_k(\A)\cong\A^{n^k}$.
\end{lem}

\begin{proof}
    First note that the kernels $ker(h_i)$ must be distinct. For
    assume that $ker(h_j)=ker(h_i)$ for some $j\neq i$. Then
    $$\A\cong\F_k(\A)/ker(h_j)=\F_k(\A)/ker(h_i)\cong\A$$ and so, since
    $h_i$ and $h_j$ are different homomorphisms, we have found a non-trivial
    automorphism on $\A$ contrary to our assumption. Now since the kernels are
    maximal by Lemma
    \ref{q-fp2.2}, they are meet irreducible. Therefore, equation
    (\ref{q-fp-representation}) provides a representation of $0$ as a
    meet of meet irreducible elements. Since $\A$ generates a
    congruence distributive variety and since $\F_k(\A) \in \VV(\A)$, the congruence lattice of $\F_k(\A)$ is
    distributive.
    Now assume the representation of $0$ as the meet of
    the kernels of our homomorphisms can be shortened; say
    $$
        \bigwedge_{\genfrac {}{}{0pt}{1}{1\leq l \leq n^k}{l\neq i}}ker(h_l)=\bigwedge_{1\leq l \leq n^k}ker(h_l)=0
    $$
    for some $1\leq i\leq n^k$. Then, just like in the proof of Theorem
    \ref{q-fp-uniqueFact},
    we get that there exists $j\neq
    i$ such that $ker(h_i)\geq ker(h_j)$, contradicting either the
    maximality or the distinctiveness of the kernels.
\end{proof}

Now we can prove $\A$ primal.

\begin{proof}[Proof of Theorem \ref{q-fp3.1}]
    Let $f$ be any $k$-ary function on $A$. Then by the preceding
    lemma, there exists a class $\Phi$ in $\F_k(\A)$ such that for
    every $\phi \in \Phi$ and for all $1\leq i \leq n^k$ the identity
    \[\phi(e_i(\xi_1),...,e_i(\xi_k))=f(e_i(\xi_1),...,e_i(\xi_k))\]
    holds. But since $(e_i(\xi_1),...,e_i(\xi_k))$ runs through all
    $k$-tuples of elements of $\A$, the term operation $\phi$ is
    identical with $f$ and so $\A$ is primal.
\end{proof}

    \chapter{The clones from $RBL$ are maximal}

In the previous chapter we demonstrated that all maximal clones
over a finite set are of the form $Pol(\rho)$, where $\rho$ is a
relation in $RBL$, and provided a characterization of primal
algebras. This chapter is devoted to the proof of the converse
statement, namely that all clones of that kind are indeed maximal.
It will therefore result in the aim of this work, the
characterization of maximal clones or preprimal algebras
respectively. We will consider the six types of relations of our
main theorem one after the other. For the first three of them, the
same technique (see M. Goldstern and S. Shelah \cite{GS02}) will
be used for the proof. Each of the other three requires a special
treatment; in those cases, we will essentially follow the original
proof of I. G. Rosenberg and include a result of J. S\l upecki on
functional completeness.

\begin{defn}
    For a set of functions $F\subseteq \FF$, we define the
    \emph{closure} $\mathopen{<}F\mathclose{>}$ of $F$ to be the smallest clone containing
    $F$.
\end{defn}

It is thus our goal to prove $\mathopen{<}
Pol(\rho)\cup\{g\}\mathclose{>} =\FF$ for all relations $\rho$ in
$RBL$ and all $g\nin Pol(\rho)$.

\begin{nota}\label{not-ch3}
    Throughout this chapter, the letter $\kappa$ will be reserved to
    denote the cardinality of our finite base set $A$. Moreover,
    $\{\alpha_1,...,\alpha_\kappa\}=A$ will be a fixed enumeration of
    $A$.
\end{nota}

\section{Partial orders with least and greatest element}

Let $\rho\subseteq A^2$ be a partial order with least and greatest
element. For $a, b\in A^n$ we write $a\leq b$ iff
$(a_i,b_i)\in\rho$ for all $1\leq i\leq n$ and $a<b$ iff $a\leq b$
and $a\neq b$. Let $g\nin Pol(\rho)$ be an $n$-ary non-monotone
function, that is there exist $a,b\in A^n$ such that $a\leq b$ but
$g(a)\nleq g(b)$. Since $\rho$ has a greatest element (and
$|A|>1$), it is non-trivial and $g$ exists.

\begin{thm}\label{thm-po}
    If $\rho\subseteq A^2$ is a partial order with least and greatest element,
    then $Pol(\rho)$ is a maximal clone.
\end{thm}

\begin{lem}\label{lem-po}
    For any $k$ and all $c,d\in A^k$, $c< d$, there exists $f_{cd}\in
    \mathopen{<}Pol(\rho)\cup\{g\}\mathclose{>}$ such that $f_{cd}(c)\nleq f_{cd}(d)$.
\end{lem}

\begin{proof}
    Our first step is to see that for given
    $c,d\in A$, $c< d$, we can construct an unary $f_{cd}\in
    \mathopen{<}Pol(\rho)\cup\{g\}\mathclose{>}$ satisfying $f_{cd}(c)\nleq f_{cd}(d)$: There
    are unary monotone functions $f^i_{cd}(x)$ mapping $c$ to $a_i$
    and $d$ to $b_i$, $1\leq i\leq n$. This is because we can map all elements
    $s$ with $s\leq c$ to $a_i$, and all other elements to $b_i$.
    Now if $f^i_{cd}(y)\nleq f^i_{cd}(z)$, then $f^i_{cd}(z)=a_i$
    which implies $z\leq c$. But as $f^i_{cd}(y)$ must equal $b_i$,
    $y\nleq c$ and so $y\nleq z$. Hence, the functions $f^i_{cd}(x)$
    are indeed monotone. Now set
    $f_{cd}(x)=g(f^1_{cd}(x),...,f^n_{cd}(x))$. Then
    $f_{cd}(c)=g(a)\nleq g(b)=f_{cd}(d)$ which is exactly what we
    wanted. \newline Next note that we can do the same thing for
    arbitrary tuples $c,d\in A^k$ with $c< d$: Choose $1\leq i\leq n$ such that
    $c_i\neq d_i$. Since $c_i< d_i$,
    we can construct $f_{c_id_i}$ as shown before and then define
    $f_{cd}=f_{c_id_i}\circ\pi^k_i$.
\end{proof}

\begin{proof}[Proof of Theorem \ref{thm-po}]
    Let $h$ be an arbitrary $k$-ary function. Using the
    functions we just constructed in the preceding lemma, we will show
    $h\in \mathopen{<}Pol(\rho)\cup\{g\}\mathclose{>}$. Consider the set $S=\{f_{cd}|c,d\in
    A^k, c< d\}$, denote it for reasons of simpler notation by $\{f_i|i\in
    I\}$, and define a mapping
    $$
        ext:\quad
        \begin{matrix}
            A^k      && \To     && A^{k+|I|} \\
            x        && \mapsto && (x,(f_i(x))_{i\in I}).
        \end{matrix}
    $$
    Then for all distinct $x,y\in A^k$ we have that $ext(x)\nleq ext(y)$. This
    is trivial if $x\nleq y$, and if otherwise, then the function
    $f_{xy}$ satisfying $f_{xy}(x)\nleq f_{xy}(y)$ is an element of
    $S$ so that by the definition of $ext$, $ext(x)\nleq ext(y)$. Now
    define an operation $H$ on the range $\{ext(x)|x\in A^k\}$ of
    $ext$ by $H(ext(x))=h(x)$. $H$ respects $\rho$ as on its domain no
    elements are comparable. We can find a monotone continuation
    $\tilde{H}$ of $H$ by setting for all $x$ not in the range of $ext$
    $$
        \tilde{H}(x)=
        \begin{cases}
            o & ,\exists y\in A^k (ext(y)\leq x)\\
            z & , otherwise
        \end{cases}
    $$
    where $o$ is the greatest element and $z$ the least element of $\rho$. But $\tilde{H}\in
    Pol(\rho)$, and so, as obviously $h(x)=\tilde{H}(x,(f_i(x))_{i\in I})$, we
    get that $h\in \mathopen{<}Pol(\rho)\cup\{g\}\mathclose{>}$.
\end{proof}

\section{Non-trivial equivalence relations}
Let $\rho\subseteq A^2$ be a non-trivial equivalence relation on
$A$. For $a,b\in A^n$ we write $a\sim b$ iff $(a_i,b_i)\in\rho$
for all $1\leq i\leq n$. Obviously $\sim$ is an equivalence
relation on $A^n$. Let $g\nin Pol(\rho)$ be an $n$-ary function
not preserving $\rho$, that is, there are $a,b\in A^n$ such that
$a\sim b$ but $g(a)\nsim g(b)$. As $\rho$ is non-trivial, $g$
exists.

\begin{thm}\label{thm-er}
    If $\rho\subseteq A^2$ is a non-trivial equivalence relation,
    then $Pol(\rho)$ is a maximal clone.
\end{thm}

Just like with partial orders, the following lemma is a fact.

\begin{lem}\label{lem-er}
    For any $k$ and all distinct $c,d\in A^k$, $c\sim d$, there exists $f_{cd}\in
    \mathopen{<}Pol(\rho)\cup\{g\}\mathclose{>}$ such that $f_{cd}(c)\nsim f_{cd}(d)$.
\end{lem}

\begin{proof}
    As in Lemma \ref{lem-po}, for arbitrary distinct
    $c,d\in A$, $c\sim d$, we construct an unary $f_{cd}\in
    \mathopen{<}Pol(\rho)\cup\{g\}\mathclose{>}$ such that $f_{cd}(c)\nsim f_{cd}(d)$. Define functions $f^i_{cd}(x)$,
    $i=1,...,n$ by mapping $c$ to $a_i$ and all
    other elements to $b_i$. Obviously, as $a_i\sim b_i$, $f^i_{cd}\in Pol(\rho)$ for all $1\leq i\leq
    n$, and setting
    $f_{cd}(x)=g(f^1_{cd}(x),...,f^n_{cd}(x))$ yields the desired function. \newline For
    arbitrary distinct tuples $c,d\in A^k$ with $c\sim d$, we define
    $f_{cd}=f_{c_id_i}\circ\pi^k_i$, where $i$ is arbitrary with $c_i\neq d_i$, and the lemma follows.
\end{proof}

\begin{proof}[Proof of Theorem \ref{thm-er}]
    Let $h$ be an arbitrary $k$-ary function. Following the proof of Theorem \ref{thm-po}, we define the
    functions $ext$ and $H$. Again, since no elements in the image of $ext$ are equivalent with
    respect to $\sim$, we can extend $H$ to $\tilde{H}\in
    Pol(\rho)$ by mapping all members of an equivalence class $e$ to
    a fixed element $x_e$ of $A$. The element $x_e$ is determined if $ext(x)\in e$ for some $x\in A^k$; otherwise,
    it can be chosen arbitrarily. Hence,
    $h(x)=\tilde{H}(x,(f_i(x))_{i\in I})\in \mathopen{<}Pol(\rho )\cup\{g\}
    \mathclose{>}$.
\end{proof}

\section{Prime permutations}

Let $\rho\subseteq A^2$ be the graph of a prime permutation $\pi$
 on $A$. For an element $a$ of $A^n$ and $l\geq 1$ we write $a+l$ for the $n$-tuple
 $(\pi^l(a_1),...,\pi^l(a_n))$. Then on $A$, $(a,b)\in\rho$ means exactly that $a+1=b$.
  We call two elements $a,b\in A^n$
 \emph{parallel} iff there is an $l\geq 1$
 such that $a+l=b$. Clearly, by that notion an equivalence
 relation is defined on $A^n$ for every $n$. Let $g\nin Pol(\rho)$ be an
$n$-ary function not preserving $\rho$, that is, there are $a,b\in
A^n$ such that $a+1=b$ but $g(a)+1\neq g(b)$.

\begin{thm}\label{thm-pp}
    If $\rho\subseteq A^2$ is a prime permutation,
    then $Pol(\rho)$ is a maximal clone.
\end{thm}

Similarly to the preceding two cases we have:

\begin{lem}\label{lem-pp}
    Let $k\geq 1$ and $c\in A^k$ with $c+l=d$ for some $1\leq l\leq p-1$. Then there
     exists $f_{cd}\in
    \mathopen{<}Pol(\rho)\cup\{g\}\mathclose{>}$ such that $f_{cd}(c)+l\neq f_{cd}(d)$.
\end{lem}

\begin{proof}
    Our first assertion is that there are $\tilde{a},\tilde{b}\in A^n$ with $\tilde{a}+l=\tilde{b}$ but
    $g(\tilde{a})+l\neq g(\tilde{b})$. For assume $g(\tilde{a}+l)=g(\tilde{a})+l$
    for all $\tilde{a}\in A^n$; then if we add $l$ to $a$ for $l^{-1}$ times,
    where $l^{-1}$ is the multiplicative inverse of $l$ modulo
    $p$, we get that
    $g(a+1)=g(a+l^{-1}\,l)=g(a)+l^{-1}\,l=g(a)+1$, contradiction.
    Now if $c,d\in A$, it is clear that there are functions $f^i_{cd}\in Pol(\rho)$ such that
    $f^i_{cd}(c)=\tilde{a}_i$ and $f^i_{cd}(d)=f^i_{cd}(c+l)=f^i_{cd}(c)+l=\tilde{b}_i$  for all $1\leq i\leq n$.
    Thus,
    $f_{cd}(x)=g(f^1_{cd}(x),...,f^n_{cd}(x))\in \mathopen{<}Pol(\rho)\cup\{g\}\mathclose{>}$ satisfies the assertion of the lemma.
    In the case of tuples $c,d\in A^k$, we do as before
    and set
    $f_{cd}=f_{c_1d_1}\circ\pi^k_1$.
\end{proof}

\begin{proof}[Proof of Theorem \ref{thm-er}]
    Let $h$ be an arbitrary $k$-ary function. Again we define the
    functions $ext$ and $H$. Now obviously no elements in the image of $ext$ are parallel. Since
    the value of an element under a function in $Pol(\rho)$ determines only the values of its parallel class, we find an
    extension $\tilde{H}$ of $H$ such that $\tilde{H}\in
    Pol(\rho)$. Therefore, as
    $h(x)=\tilde{H}(x,(f_i(x))_{i\in I})\in \mathopen{<}Pol(\rho )\cup\{g\}
    \mathclose{>}$, it follows that $Pol(\rho)$ is a maximal clone.
\end{proof}

\section{Central relations}

We will show that every central relation $\rho\subseteq A^h$
yields a maximal clone via $Pol$. We distinguish the possibilities
$h=1$, in which case $\rho$ is just a proper subset of $A$, and
$h\geq 2$. In the first case, the method we used so far can be
applied once again; however, in all other cases the issue is more
complicated. As before we denote by $g\nin Pol(\rho)$ the $n$-ary
function not preserving $\rho$; $g$ exists as the center of a
central relation is non-trivial by definition. Thus, there exist
$a_1,...,a_n \in \rho$ such that
$(g(a_{11},...,a_{n1}),...,g(a_{1h},...,a_{nh}))\nin\rho$. The
following theorem does the case $h=1$.

\begin{thm}\label{thm-su}
    If $\rho\subseteq A$ is a proper subset of $A$, then
    $Pol(\rho)$ is a maximal clone.
\end{thm}

\begin{lem}\label{lem-su}
    For every $c\in\rho^k$ there is a $f_c\in \mathopen{<}Pol(\rho)\cup\{g\}\mathclose{>}$ with
     $f_{c}(c)\nin\rho$.
\end{lem}

\begin{proof}
    There are $a_1,...,a_n\in\rho$ such that
    $g(a_1,...,a_n)\nin\rho$. If $c\in\rho$, then there are obviously
    mappings $f^i_c\in Pol(\rho)$ with $f^i_c(c)=a_i$. Setting
    $f_c=g(f^1_c,...,f^n_c)$ proves the lemma for this case. If
    $c$ is a $k$-tuple, define $f_c=f_{c_1}\circ \pi_1^k$ as
    usually.
\end{proof}

\begin{proof}[Proof of Theorem \ref{thm-su}]
    Take any $k$-ary function $h$ and define for every $x\in A^k$
    the tuple $ext(x)$ by $ext(x)=(x,(f_c(x))_{c\in\rho^k})$. On
    the range of $ext$, set $H(ext(x))=h(x)$. Clearly, as tuples
    of the form $ext(x)$ can never have all their components in
    $\rho$, we can extend $H$ to $\tilde{H}\in Pol(\rho)$ like in the
    previous sections and as
    $h(x)=\tilde{H}(x,(f_c(x))_{c\in\rho^k})$ the theorem has been
    proven.
\end{proof}

\subsubsection{A completeness criterion}

For the central relations as well as the $h$-regularly generated
relations we will need a completeness criterion due to J. S\l
upecki
 saying that for $|A|\geq 3$, if $F\subseteq \FF$ contains all unary
functions and a function which takes all values of $A$ and which
depends on at least two variables, then
$\mathopen{<}F\mathclose{>}=\FF$. This criterion will be proven
now; we will essentially follow a proof by J. W. Butler in
\cite{But60}. The restriction $|A|\geq 3$ does not matter to us:
For central relations we use the criterion only for the case
$2\leq h < \kappa$; in the case of $h$-regularly generated
relations, $3\leq h\leq\kappa$ by definition.

\begin{defn}
    An $n$-ary function $f(x_1,...,x_n)$ \emph{depends} on the
    $j$-th variable, $1\leq j\leq n$, iff there exist $a\in A^n$
    and $u\in A$ such that $f(a_1,...,a_n)\neq
    f(a_1,...,a_{j-1},u,a_{j+1},...,a_n)$. We call $f$
    \emph{irreducible} iff it depends on at least two variables
    and \emph{reducible} iff it does not.
\end{defn}

We denote the \emph{range} of a function $f$ by $\Re(f)$.

\begin{lem}\label{but-1}
    Let $\kappa\geq 3$ and let $f$ be an irreducible function of $n$
    arguments, $n\geq 3$, which is onto. Then there is an irreducible function $g$
    of two variables in $\mathopen{<}\FF_1\cup\{f\}\mathclose{>}$ which is onto.
\end{lem}

\begin{proof}
    There are $1\leq q\leq n$, $a\in A^n$ and $u\in
    A$ such that $f(\tilde{a})\neq
    f(a)$ if we set $\tilde{a}=(a_1,...,a_{q-1},u,a_{q+1},...,a_n)$. Say $f(a)=\alpha_1$ and
    $f(\tilde{a})=\alpha_2$ and choose $y_i\in A^n$, $3\leq i\leq
    \kappa$ such that $f(y_i)=\alpha_i$. Note next that there exist
    $w,z\in A^n$ with $w_q=z_q$ but $f(w)\neq f(z)$, for otherwise
    $f$ would depend only on its $q$-th argument and would
    therefore be reducible. We distinguish two cases: First, such $w$
    and $z$ exist with the additional property that $f(w)\neq
    f(a)$ and $f(w)\neq f(\tilde{a})$, and second, no such $w$ and $z$
    fulfill this additional assumption.\newline
    In the first case, say without loss of generality
    $f(w)=f(y_3)$. We define $n-1$ unary functions $h_i$, $1\leq
    i\leq n$, $i\neq q$ by
    $$
        h_i(x)=
        \begin{cases}
            a_i&,x=\alpha_1\\
            z_i&,x=\alpha_2\\
            w_i&,x=\alpha_3\\
            y_{ji}&,x=\alpha_j\wedge j\nin\{1,2,3\}\\
        \end{cases}
    $$
    and $h_q$ by
    $$
        h_q(x)=
        \begin{cases}
            a_q&,x=\alpha_1\\
            u&,x=\alpha_2\\
            z_q&,x=\alpha_3\\
            y_{jq}&,x=\alpha_j\wedge j\nin\{1,2,3\}\\
        \end{cases}
    $$
    and set
    $g(x,y)=f(h_1(x),...,h_{q-1}(x),h_q(y),h_{q+1}(x),...,h_n(x))$.
    Then $g$ is onto since $g(\alpha_1,\alpha_1)=f(a)=\alpha_1$, $
    g(\alpha_1,\alpha_2)=f(\tilde{a})=\alpha_2$, $g(\alpha_3,\alpha_3)=f(w)=\alpha_3$, and
    $g(\alpha_i,\alpha_i)=f(y_i)=\alpha_i$ for $i> 3$. Moreover, $g$ is not
    reducible as $g(\alpha_1,\alpha_1)\neq g(\alpha_1,\alpha_2)$ and $g(\alpha_2,\alpha_3)\neq
    g(\alpha_3,\alpha_3)$.\newline
    In the second case, we choose for $1\leq i\leq n$, $i\neq q$ functions $h_i$ satisfying
    $$
        h_i(x)=
        \begin{cases}
            a_i&,x=\alpha_1\\
            w_i&,x=\alpha_2\\
            z_i&,x=\alpha_3.\\
        \end{cases}
    $$
    Define
    $g(x,y)=f(h_1(x),...,h_{q-1}(x),y,h_{q+1}(x),...,h_n(x))$. Now
    the condition of this case implies that if $s\in A^n$ and
    $s_q=y_{iq}$ for some $3\leq i\leq \kappa$, then $f(s)=f(y_i)$
    since $f(y_i)\neq f(a)$ and $f(y_i)\neq f(\tilde{a})$. Thus,
    $g(\alpha_1,a_q)=g(a)=\alpha_1$, $g(\alpha_1,u)=f(\tilde{a})=\alpha_2$, and
    $g(\alpha_m,y_{iq})=f(y_i)=\alpha_i$ for any $1\leq m\leq n$
    and $3\leq i\leq \kappa$. Hence, $g$ is onto. Moreover,
    $g(\alpha_1,a_q)=g(a)\neq g(\tilde{a})=g(\alpha_1,u)$ and
    $g(\alpha_2,w_q)=g(w)\neq g(z)=g(\alpha_3,w_q)$ and so $g$ is
    irreducible.
\end{proof}

\begin{lem}\label{but-2}
    If $f\in\FF_2$ is an irreducible function of two variables
    which takes at least three distinct values, then there exist
    $a,b,c,d\in A$ such that $f$ takes three distinct values on
    $\{(a,c),(a,d),(b,c),(b,d)\}$.
\end{lem}

\begin{proof}
    Assume first that there is an $a\in A$ such that $f$ takes at least three values on
     $\{(a,x)|x\in A\}$. Since $f$ is
    irreducible, there must be $b,c\in A$ such that $f(a,c)\neq
    f(b,c)$. As $f$ takes at least three values with $a$ as the
    first argument, there is $d\in A$ with $f(a,d)\neq f(a,c)$ and
    $f(a,d)\neq f(b,c)$.\newline
    Consider now the case where there is no such $a$. The
    irreducibility of $f$ implies there is $a\in A$ such that $f$
    takes two values with $a$ as the first argument.
    It follows from the assumption for this case
    that there is $w$ in the range of $f$ such that $f(a,x)\neq w$ for all $x\in
    A$; say $w=f(b,c)$, $b\neq a$. Hence, $f(a,c)\neq f(b,c)$. Now
    take any $d\in A$ with $f(a,c)\neq f(a,d)$ to finish the proof.
\end{proof}

\begin{lem}\label{but-3}
    If $f\in\FF_2$ is an irreducible function of two variables
    with $|\Re(f)|=p$, $p\geq 3$, then there exist two
    unary functions $h_1,h_2\in\FF_1$ which both take at most
    $p-1$ elements such that for every $x\in\Re(f)$
    we have $f(h_1(x),h_2(x))=x$.
\end{lem}

\begin{proof}
    Let $a,b,c,d$ be provided by Lemma \ref{but-2}. Assume without
    loss of generality that $f(a,c)=u, f(a,d)=v, f(b,c)=w$ are all
    different. We define $h_1,h_2\in \FF_1$ as follows: For
    $u,v,w$ we set
    $h_1(u)=a$, $h_1(v)=a$, $h_1(w)=b$ and $h_2(u)=c$, $h_2(v)=d$, $h_2(w)=c$;
    for $x\nin\Re(f)$, we define $h_1(x)=a$ and $h_2(x)=c$; and
    for $x\in\Re(f)\setminus\{u,v,w\}$ we choose any values for
    $h_1(x)$ and $h_2(x)$ such that the requirement
    $f(h_1(x),h_2(x))=x$ is satisfied. Clearly, $h_1,h_2$ have all
    desired properties.
\end{proof}

To proof the completeness criterion, we want to construct the
function returning the maximum of to elements with respect to some
total ordering of the elements of $A$. Therefore we will for the
rest of this section replace $A$ by the set of natural numbers
$\kappa=\{0,...,\kappa-1\}$ and use the standard notions of
$\leq$, $\vee$, $+$ and $-$ on that set.
\begin{lem}\label{but-4}
    Let $f\in\FF_2$, $p<\kappa$, and assume there exist $i,j,l\in
    A$ such that for all $y<p$ we have $f(i,y)=y$ and $f(j,y)=l$.
    Then there is a function of two variables $g\in
    \mathopen{<}\FF_1\cup\{f\}\mathclose{>}$ such that $g(x,y)=x\vee y$ for $x,y<p$.
\end{lem}

\begin{proof}
    We may assume without loss of generality that $i,j,l<p$. This
    is legitimate as we can shift them with unary functions. The
    proof will be by induction on $p$. First, let $p=2$; then
    there are four possibilities: either $i=0,j=1,l=1$ or
    $i=0,j=1,l=0$ or $i=1,j=0,l=0$ or $i=1,j=0,l=1$. In multiplication tables of the restriction of $f$ to $\{0,1\}^2$, these scenarios
    look like this:
    $$
    (i)\quad
    \begin{tabular}{c|cc}
        &0&1\\
        \hline
        0&0&1\\
        1&1&1\\
    \end{tabular}
    \qquad (ii) \quad
    \begin{tabular}{c|cc}
        &0&1\\
        \hline
        0&0&1\\
        1&0&0\\
    \end{tabular}
    \qquad (iii) \quad
    \begin{tabular}{c|cc}
        &0&1\\
        \hline
        0&0&0\\
        1&0&1\\
    \end{tabular}
    \qquad (iv) \quad
    \begin{tabular}{c|cc}
        &0&1\\
        \hline
        0&1&1\\
        1&0&1\\
    \end{tabular}
    $$\newline
     In the first
    case, $f$ is the maximum function on $\{0,1\}$ and we can take
    $f$ itself for $g$. In the other cases we use any unary
    function $h$ exchanging $0$ and $1$ to define $g(x,y)$ to be
    $h(f(x,h(y)))$ or $h(f(h(x),h(y)))$ or $f(h(x),y)$,
    respectively.\newline
    Now assume the lemma is true for $p-1$, and let $g'$ be a
    function in $\mathopen{<}\FF_1\cup\{f\}\mathclose{>}$ satisfying $g'(x,y)=x\vee y$
    for $x,y<p-1$. Choose functions $h_1,h_2\in\FF_1$ such that
    $$
        h_1(x)=
        \begin{cases}
            i&,x<p-1\\
            j&,x=p-1
        \end{cases}
        \qquad\qquad
        h_2(x)=
        \begin{cases}
            p-1&,x=l\\
            l&,x=p-1\\
            x&,otherwise
        \end{cases}
    $$
    and construct $f'\in\FF_2$ as $f'(x,y)=h_2(f(h_1(x),h_2(y)))$.
    It is easy to check that
    $$
        f'(x,y)=
        \begin{cases}
            y&,x<p-1\wedge y<p\\
            p-1&,x=p-1\wedge y<p.
        \end{cases}
    $$
    Now we define $g$ by $g(x,y)=f'(f'(x,y),g'(x,y))$. One readily
    verifies that for $x,y<p-1$, $g(x,y)=g'(x,y)=x\vee y$; for
    $x=p-1$, $y<p$, $g(x,y)=f'(p-1,g'(p-1,y))=p-1$; and for $x<p-1$,
    $y=p-1$, $g(x,y)=f'(p-1,g'(x,p-1))=p-1$. Hence, $g(x,y)=x\vee y$
    for $x,y<p$.
\end{proof}

\begin{lem}\label{but-5}
    If $f\in\FF_2$ is irreducible and $\Re(f)=\{0,...,p-1\}$, $3\leq p\leq \kappa$, then there is a binary
    function $g\in\mathopen{<}\FF_1\cup\{f\}\mathclose{>}$ such that $g(x,y)=x\vee y$ for $x,y<p$.
\end{lem}

\begin{proof}
    The proof will be by induction on $p$. If $p=3$, by Lemma \ref{but-2} there are $a,b,c,d$
    such that $f$ takes at least
    three distinct values on $\{(a,c),(a,d),(b,c),(b,d)\}$. By
    shifting those elements and their values under $f$ with unary
    functions, we may assume that $a=c=0$, $b=d=1$, $f(0,0)=0$,
    $f(0,1)$=1, and $f(1,0)=2$. This leaves us essentially with two
    possible multiplication tables:\newline
    $$
    (i)\quad
    \begin{tabular}{c|cc}
        &0&1\\
        \hline
        0&0&1\\
        1&2&2\\
    \end{tabular}
    \qquad (ii) \quad
    \begin{tabular}{c|cc}
        &0&1\\
        \hline
        0&0&1\\
        1&2&0\\
    \end{tabular}
    $$\newline
    In the first case, we choose functions
    $h_1,h_2\in\FF_1$ with
    $$
        \begin{matrix}
            h_1(0)=0&\qquad h_2(0)=0\\
            h_1(1)=0&\qquad h_2(1)=1\\
            h_1(2)=1&\qquad h_2(2)=1.
        \end{matrix}
    $$
    Then we can construct $g$ as
    $$
        g(x,y)=f(h_2\circ f(h_1(x),h_1(y)),h_2\circ f(h_2(x),h_2(y))).
    $$
    To construct $g$ in the other case we choose additional functions $h_3,h_4 \in\FF_1$ with
    $$
        \begin{matrix}
            h_3(0)=2\\
            h_3(1)=0\\
            h_3(2)=1
        \end{matrix}\qquad\qquad
        \begin{matrix}
            h_4(0)=1\\
            h_4(1)=0
        \end{matrix}
    $$
    and define $g' \in\FF_2$ by
    \begin{eqnarray*}
        g'(x,y)=h_3\circ f(y,h_2\circ f(x,h_4(y))).
    \end{eqnarray*}
    It is boring but possible to verify that $g'$ agrees on $\{0,1\}$ with the $f$ of the first case which we already
    treated.\newline
    Assuming our assertion is true for $p-1$, we prove it for $p$,
    $3<p\leq n$. First we construct a function $f''$ from $f$ satisfying
    the hypotheses of the lemma for $p-1$; we need to restrict the
    range of $f$ to $\{0,...,p-2\}$ without making $f$ reducible.
    To do this we apply Lemma \ref{but-2}, taking $a,b,c,d\in A$
    such that at least three distinct values $u,v,w\in A$ are represented
    among
    $\{f(a,c),f(a,d),f(b,c),f(b,d)\}$. Since $p>3$ there is
    $z\in\Re(f)\setminus\{u,v,w\}$. Define $h\in\FF_1$ by
    $$
        h(x)=
        \begin{cases}
            u&,x=z\\
            x&,otherwise.
        \end{cases}
    $$
    Then $h(f(x,y))\in\FF_2$ is not reducible and has $p-1$
    elements in its range. By permuting the elements of $A$ with
    an unary function we produce a function $f''\in\FF_2$
    satisfying the hypotheses of the lemma for $p-1$, and hence by induction hypothesis we
    get a function $g''\in\FF_2$ such that $g(x,y)=x\vee y$ for
    $x,y<p-1$.\newline
    Next by Lemma \ref{but-3} there exist functions
    $h_1,h_2\in\FF_1$ with $\Re(h_1),\Re(h_2)$ consisting of at
    most $p-1$ elements such that $f(h_1(x),h_2(x))=x$ for $x<p$.
    There exist permutations $h_3,h_4\in\FF_1$ such
    that $h_3(x)<p-1$ for all $x\in\Re(h_1)$ and $h_4(x)<p-1$ for
    all $x\in\Re(h_2)$. Define $h_5,h_6\in\FF_1$ and $f'\in\FF_2$
    by
    \begin{eqnarray*}
    \begin{split}
        h_5&=h_3\circ h_1\\
        h_6&=h_4\circ h_2\\
        f'(x,y)&=f(h_3^{-1}(x),h_4^{-1}(y)).
    \end{split}
    \end{eqnarray*}
    Then obviously $f'(h_5(x),h_6(x))=x$ for all $x<p$ and
    $\Re(h_5),\Re(h_6)$ are subsets of $\{0,...,p-2\}$. We define
    $g'\in\FF_2$ by $g'(x,y)=f'(g''(x,h_5(y)),g''(x,h_6(y)))$.
    Then $g'$ satisfies the equation
    $g'(0,y)=f'(h_5(y),h_6(y))=y$ for $y<p$; moreover, for $y<p$,
    $g'(p-2,y)=f'(p-2,p-2)$ and is therefore constant. Hence by
    Lemma \ref{but-4} we can generate a function $g\in\FF_2$ such
    that $g$ agrees with the maximum function for arguments
    smaller than $p$.
\end{proof}

\begin{thm}\label{but-6}
    Assume $|A|\geq 3$ and let $f\in\FF$ be an irreducible
    function with $\Re(f)=A$. Then $\mathopen{<}\FF_1\cup\{f\}\mathclose{>}=\FF$.
\end{thm}

\begin{proof}
    By Lemma \ref{but-1} we may assume $f$ is a function of two
    variables. Lemma \ref{but-5} then implies that
    $\mathopen{<}\FF_1\cup\{f\}\mathclose{>}$ contains the maximum function with respect
    to some total order of the elements of $A$. But is well-known
    from the results of E. L. Post
    and easily verified that the unary functions together with the
    maximum function already generate all functions of arbitrary
    arity over $A$; thus, $\mathopen{<}\FF_1\cup\{f\}\mathclose{>}=\FF$.
\end{proof}

\subsubsection{Totally reflexive and totally symmetric relations}

The following lemmas hold for totally reflexive and totally
symmetric relations. They will help us with both the central
relations with $h\geq 2$ and the $h$-regularly generated
relations. The first lemma implies that we can assume without loss
of generality that the function $g$ not preserving $\rho$ is
unary.

\begin{lem}\label{ros-9.3.1}
    Let $\rho\neq \iota_h^A$ be a totally reflexive and totally symmetric $h$-ary relation. If $g\nin Pol(\rho)$
    then there is an unary $f\in \mathopen{<}Pol(\rho)\cup\{g\}\mathclose{>}$ that does not preserve $\rho$.
\end{lem}

\begin{proof}
    Let $a_1,...,a_n \in \rho$ such that
    $(g(a_{11},...,a_{n1}),...,g(a_{1h},...,a_{nh}))\nin\rho$. Choose
    $(c_1,...,c_h)\in\rho$, $(c_1,...,c_h)\nin \iota_h^A$. Define for  $1\leq i\leq n$ unary functions $f_i$ by
    $f_i(c_j)=a_{ij}$, $1\leq j\leq h$, and $f_i(x)=a_{i1}$ for
    all other elements $x\in A$. The operations $f_i$ preserve
    $\rho$ as they map just any tuple to a tuple in $\rho$:
    If an $h$-tuple consisting of function values of $f_i$
    has two identical entries, then the tuple is an element of
    $\rho$ as $\rho$ is totally reflexive; if
    otherwise, then the definition of $f_i$ implies that the tuple contains the values
    $a_{i1},...,a_{in}$ in some order and is thus in $\rho$ by its
    total symmetry. Now $f(x)=g(f_1(x),...,f_n(x))\in
    \mathopen{<}Pol(\rho)\cup\{g\}\mathclose{>}$ maps $(c_1,...,c_h)\in\rho$ to
    $(g(a_{11},...,a_{n1}),...,g(a_{1h},...,a_{nh}))\nin\rho$.
\end{proof}

\begin{lem}\label{ros-9.3.2}
    Let $\rho\neq \iota_h^A$ be a totally reflexive and totally symmetric $h$-ary relation, $1\leq h\leq \kappa$. If $g$ is
    an unary function not preserving $\rho$, then there is a
    subset $D=\{d_1,...,d_h\}$ of $A$ such that
    $(d_1,...,d_h)\nin\rho$ and $\mathopen{<}Pol(\rho)\cup\{g\}\mathclose{>}$ contains
    all unary functions which take only values in $D$.
\end{lem}

\begin{proof}
    There is $(a_1,...,a_h)\in\rho$ such that
    $(g(a_1),...,g(a_h))\nin\rho$. If we set $d_i=g(a_i)$, $1\leq i\leq
    h$, then $(d_1,...,d_h)\nin\rho$. Let $h$ be an unary function that takes only values in
    $D$. Define a function $l$ by $l(x)=a_i$ whenever $h(x)=d_i$.
    Then by the same argument as in the preceding lemma for $f_i$,
    $l\in Pol(\rho)$. Hence, $h=g\circ l\in \mathopen{<}Pol(\rho)\cup\{g\}\mathclose{>}$.
\end{proof}

\begin{thm}\label{ros-9.3.3}
    Let $\rho\neq \iota_h^A$ be a totally reflexive and totally symmetric non-trivial $h$-ary relation,
    where $2\leq h\leq \kappa$. If for every $D=\{d_1,...,d_h\}$
    with $(d_1,...,d_h)\nin\rho$ an $n$-ary function $q\in Pol(\rho)$
    exists which takes all values of $A$ on $D^n$, then
    $Pol(\rho)$ is a maximal clone.
\end{thm}

\begin{proof}
    Take an unary $g\nin Pol(\rho)$. Let $D$ be provided by Lemma
    \ref{ros-9.3.2}. By our hypothesis, there are
    $q(x_1,...,x_n)\in Pol(\rho)$ and $a_i\in D^n$, $1\leq i\leq \kappa$, such
    that $q(a_i)=\alpha_i$. Let $h\in\FF_1$ be given. Define for
    $1\leq j\leq n$ functions $g_j\in\FF_1$ by $g_j(x)=a_{ij}$ whenever $h(x)=\alpha_i$. As the
    $g_j$ obviously take only values in $D$, we have $g_j\in
    \mathopen{<}Pol(\rho)\cup\{g\}\mathclose{>}$ and so the same holds for
    $q(g_1,...,g_n)$. But it is easily verified that $q(g_1(x),...,g_n(x))=h(x)$ for
    all $x\in A$, and so $h\in \mathopen{<}Pol(\rho)\cup\{g\}\mathclose{>}$. We have thus shown that $\mathopen{<}Pol(\rho)\cup\{g\}\mathclose{>}$
    contains all unary functions. Now assume
    $q$ depends only on one variable. Then, as $q$ takes all
    values of $A$, we necessarily have that $D=A$ and so $h=\kappa$.
    Therefore, $(\alpha_1,...,\alpha_\kappa)\nin\rho$. But this implies
    $\rho=\iota_h^A$, contradictory to our assumption. Hence, $q$ depends
    on at least two variables, and we can apply Theorem \ref{but-6} to obtain that $\mathopen{<}Pol(\rho)\cup\{g\}\mathclose{>}=\FF$.
\end{proof}

\subsubsection{Central relations with $h\geq 2$}

We make use of the preceding results for the remaining case $h\geq
2$.

\begin{lem}\label{ros-9.3.4}
    Let $\rho$ be an $h$-ary central relation, $2\leq h\leq \kappa$,
    $(d_0,...,d_{h-1})\nin\rho$, and let $u\in A$ be a central element of $\rho$.
    Enumerate all functions in $h^\kappa$ by
    $\{p_1,...,p_{h^\kappa}\}$. For $i\in \kappa$ set
    $b_i=(d_{p_1(i)},...,d_{p_{h^\kappa}(i)})$. We define a
    $h^\kappa$-ary function $q$ by $q(b_i)=\alpha_{i+1}$, $i\in \kappa$,
    and for all other elements $x\in A^{h^\kappa}$ we set $q(x)=u$.
    Then $q$ preserves $\rho$.
\end{lem}

\begin{proof}
    We first show that for distinct $i_0,...,i_{h-1}\in\kappa$,
    $(b_{i_0},...,b_{i_{h-1}})\nin \rho^{h^\kappa}$. Take any function $r\in h^\kappa$
    with $r(i_j)=j$, $j\in h$. There is an $1\leq l\leq h^\kappa$ such that
    $r=p_l$. Thus
    $(b_{i_0l},...,b_{i_{h-1}l})=(d_{p_l(i_0)},...,d_{p_l(i_{h-1})})=(d_0,...,d_{h-1})\nin\rho$
    and so $(b_{i_0},...,b_{i_{h-1}})\nin \rho^{h^\kappa}$. Let
    $(a_1,...,a_h)\in\rho^{h^\kappa}$. If $a_i=a_j$ for some $i\neq j$,
    then $(q(a_1),...,q(a_h))\in\rho$ since $\rho$ is totally
    reflexive. Otherwise, as $(a_1,...,a_h)\in\rho^{h^\kappa}$, but
    $(b_{i_0},...,b_{i_{h-1}})\nin \rho^{h^\kappa}$ for distinct
    $i_0,...,i_{h-1}$, there is an $1\leq l\leq h$ such that $a_l$ is
    not equal to any of the $b_i$. But then $q(a_l)=u$ by definition of
    $q$; hence, $(q(a_1),...,q(a_h))\in\rho$.
\end{proof}

\begin{thm}
    If $\rho$ is a central relation, then $Pol(\rho)$ is a maximal
    clone.
\end{thm}

\begin{proof}
    Follows from the previous lemma together with Theorems \ref{thm-su} and \ref{ros-9.3.3}.
\end{proof}

\section{$h$-regularly generated relations}

Our next step is to show that $h$-regularly relations generate
maximal clones. As those relations are obviously totally reflexive
and totally symmetric, the results at the beginning of the last
section apply. Our goal is therefore to show the hypotheses of
Theorem \ref{ros-9.3.3} are satisfied.

\begin{lem}\label{ros-10.5.1}
    Let $D=\{d_1,...,d_h\}$ be a subset of $A$ with
    $(d_1,...,d_h)\nin\rho$. Then there is an unary $f\in
    Pol(\rho)$ satisfying $f(d_i)=\alpha_i$, $1\leq i\leq h$.
\end{lem}

\begin{proof}
    Denote by $\varphi$ the surjection from $A$ onto $h^\lambda$
    such that
    $\rho=\varphi^{-1}(\omega_\lambda)$. Set $b_i=\varphi(d_i)$ for $1\leq
    i\leq h$; then $(b_1,...,b_h)\nin \omega_\lambda$. That means there
    is $1\leq j\leq \lambda$ such that all $b_{ij}$ are distinct,
    $1\leq i\leq h$. Thus, the function $s(x)=b_{xj}$ is a
    bijection from $\{1,...,h\}$ onto $h$. We define $f$ as follows: If
    $s^{-1}(\varphi(x)_j)=l$, then $f(x)=\alpha_l$. Then, as
    $s^{-1}(\varphi(d_i)_j)=s^{-1}(b_{ij})=i$, we get
    $f(d_i)=\alpha_i$, $1\leq i\leq h$. If $(a_1,...,a_h)\in\rho$,
    then $(\varphi(a_1),...,\varphi(a_h))\in\omega_\lambda$ and so
    $\varphi(a_i)_j=\varphi(a_l)_j$ for some $i\neq l$. Thus,
    $f(a_i)=f(a_l)$ so that $(f(a_1),...,f(a_h))\in\rho$ since
    $\rho$ is totally reflexive. Hence, $f$ preserves $\rho$.
\end{proof}

To keep the notation simple, we will sometimes identify the set of
tuples $h^\lambda$ with its interpretation as a set of natural
numbers, sometimes not, whichever is simpler. In our
interpretation, we let the numbers $i\in h$ correspond to the
tuples $(i,0,...,0)\in h^\lambda$. Observe that this implies
$(l_i)_1=l_i$ for every $l\in h^\lambda$ and every $1\leq i\leq
\lambda$.

\begin{lem}\label{ros-10.5.2}
    If we choose any enumeration $\{\beta_0,...,\beta_{\kappa-1}\}$ of
    $A$ such that $\varphi(\beta_i)=i$, $i\in h^\lambda$, then there
    is a function $r\in Pol(\rho)$ which takes all values of $A$
    on $\{\beta_0,...,\beta_{h-1}\}$.
\end{lem}

\begin{proof}
    Set $F_l=\{x\in A|\,\varphi(x)=l\}$ for $l\in h^\lambda$ and denote the
    elements of $F_l$ by $c_{l0},...,c_{ln_l}$. Set further $n^*=max\{n_l|\,l\in h^\lambda\}$ and define
    for $j\leq n^*$ the $n^*$-tuple $code(j)$ to contain
    $\beta_1$ at its $j$-th component and $\beta_0$ in all other components.
    We write $n=\lambda+n^*$ and for $i\in h^\lambda$ and $j\leq n_i$ we define
     $d_{ij}$ to be the $n$-tuple
    $(\beta_{i_{1}},...,\beta_{i_{m}},code(j))$. The set of
    all $d_{ij}$ we call $D$. We define an $n$-ary $r$ on $D$ by
    $r(d_{ij})=c_{ij}$, and for $a\in A^n\setminus D$ we set
     $r(a)=\beta_l$ where $l=(\varphi(a_1)_1,...,\varphi(a_\lambda)_1)\in h^\lambda\leq\kappa$. As every element
     $a\in A$ is for some $l\in h^\lambda$ an element of $F_l$, we have that $a=c_{lj}$ for some
     $j\leq n_l$ and so $r$ is onto. We claim
    that for all $a\in A^n$ we have that
    $\varphi(r(a))=(\varphi(a_1)_1,...,\varphi(a_\lambda)_1)$. If
    $a\in A^n\setminus D$, then this is a direct consequence of
    our assumption that $\varphi(\beta_i)=i$ for $i\in h^\lambda$. Now if
    $a=d_{ij}\in D^n$ for some $i\in h^\lambda$ and $j\leq n_i$, then $r(a)=r(d_{ij})=c_{ij}$ so that $\varphi(r(d_{ij}))=i$ by the
    definition of $F_i$. On the other hand, $a_l=\beta_{i_l}$, $1\leq
    l\leq \lambda$. Thus, again by our assumption on $\varphi$,
    $\varphi(a_l)_1=\varphi(\beta_{i_l})_1=(i_{l})_1=i_{l}$. Hence,
    $(\varphi(a_1)_1,...,\varphi(a_\lambda)_1)=(i_{1},...,i_\lambda)=i=\varphi(r(a))$. We have proven our claim.\newline
    We show that $r\in Pol(\rho)$. Let
    $(r(a_1),...,r(a_h))\nin\rho$, $a_i\in A^n$, $1\leq i\leq h$.
    Then by the definition of a $h$-regularly generated relation,
    $(\varphi(r(a_1)),...,\varphi(r(a_h)))\nin\omega_\lambda$ which means there is
    $1\leq j\leq \lambda$ such that all $\varphi(r(a_i))_j$ are distinct. By
    our last claim we have $\varphi(r(a_i))_j=\varphi(a_{ij})_1$ and so all
    $\varphi(a_{ij})_1$ are distinct, $1\leq i\leq h$. Hence, by the definition of
    $\omega_\lambda$, $(\varphi(a_{1j}),...,\varphi(a_{hj}))\nin \omega_\lambda$ and
    so $(a_{1j},...,a_{hj})\nin\rho$. Thus,
    $(a_1,...,a_h)\nin\rho^n$ and we conclude that $r$ preserves
    $\rho$.
\end{proof}

\begin{lem}\label{ros-10.5.3}
    If $D=\{d_1,...,d_h\}$ is a subset of $A$ with the property
    that
    $(d_1,...,d_h)\nin\rho$, then there is an $n$-ary $q\in Pol(\rho)$
    which takes all values of $A$ on $D^n$.
\end{lem}

\begin{proof}
    Let $\{\beta_0,...,\beta_{\kappa-1}\}$ be an enumeration of $A$
    with $\varphi(\beta_i)=i$, $i\in h^\lambda$, and let $r\in Pol(\rho)$
    be the function from Lemma \ref{ros-10.5.2}. By Lemma
    \ref{ros-10.5.1} there is a function $g\in Pol(\rho)$ with
    $g(d_i)=\beta_{i-1}$, $1\leq i\leq h$. Setting
    $q(x_1,...,x_n)=r(g(x_1),...,g(x_n))$ proves the lemma.
\end{proof}

\begin{lem}\label{ros-10.5.4}
    Let $3\leq h\leq \kappa$. If $\iota_h^A$ is $h$-regularly generated,
    then
    $\lambda=1$ and $h=\kappa$.
\end{lem}

\begin{proof}
    If $\lambda\geq 2$, then for the vectors $b_i=(i,0,...,0)$, $1\leq
    i\leq h-1$, and $b_h=(1,1,0,...,0)$ we have
    $(b_1,...,b_h)\in\omega_\lambda$. But as those tuples are all
    distinct, it is impossible that $\iota_h^A=\varphi^{-1}(\omega_\lambda)$,
    contradiction. Assume $h<\kappa$. Then, as $\lambda=1$, $\varphi$ is not
    one-one and hence there exist distinct $a_1,...,a_h\in A$ such that
    $(\varphi(a_1),...,\varphi(a_h))\in\omega_1$. But this implies
    $(a_1,...,a_h)\in \iota_h^A$, contradiction.
\end{proof}

\begin{lem}\label{a_h}
    Let $\kappa\geq 3$. If $\rho=\iota_\kappa^A$, then $Pol(\rho)$ is a maximal clone.
\end{lem}

\begin{proof}
    $Pol(\rho)$ contains all unary functions. For if $f\in \FF_1$
    and $(a_1,...,a_\kappa)\in A^\kappa$ has two identical components, then the
    same holds for $(f(a_1),...,f(a_\kappa))$. Therefore, if $g\nin
    Pol(\rho)$, $g$ must depend on at least two variables. But in order to
    produce a tuple not in $\iota_\kappa^A$, $g$ must take all values in $A$.
    Hence Theorem \ref{but-6} yields $\mathopen{<}Pol(\rho)\cup
    \{g\}\mathclose{>}=\FF$. Observe that $\kappa\geq 3$ is necessary as
    otherwise $Pol(\rho)=\FF$.
\end{proof}

\begin{thm}\label{ros-10.5.5}
    Let $3\leq h\leq \kappa$. If $\rho$ is a $h$-regularly generated
    relation, then $Pol(\rho)$ is a maximal clone.
\end{thm}

\begin{proof}
    If $h=\kappa$, then $\rho=\iota_\kappa^A$ and $Pol(\rho)$ is maximal by Lemma \ref{a_h}. Otherwise $\rho\neq \iota_h^A$
    by Lemma \ref{ros-10.5.4} and application of Lemma
    \ref{ros-10.5.3} together with Theorem \ref{ros-9.3.3} proves the theorem.
\end{proof}

\section{Prime affine relations}

Let $\rho$ be a prime affine relation with respect to $(A,+)$.
Recall that by definition $(A,+)$ is an abelian group and every
$a\in A$ has order $p$, where $p$ is a prime. It is a basic fact
from the theory of abelian groups that in this case $|A|=p^m$ for
some $m>1$. Moreover, $(A,+)$ is isomorphic to the additive group
of the field $\frak{G}\frak{F}(p^m)$ with $p^m$ elements. It is
for this reason that we can define a multiplication $\cdot$ on $A$
so that $(A,+,\cdot)$ is isomorphic to $\frak{G}\frak{F}(p^m)$.
$(A,+,\cdot)$ has a primitive element which we call $e$. The
neutral elements of $+$ and $\cdot$ we denote by $0$ and $1$
respectively. In this context, we understand a \emph{polynomial}
to be a function in $\mathopen{<}\{+,\cdot,(a)_{a\in
A}\}\mathclose{>}$. Naturally enough, our approach for proving
$\rho$ maximal will be to construct polynomials. We recall the
following fact:

\begin{lem}\label{ros-7.2.5}
    Every $f\in\FF$ is a polynomial. Furthermore, $f(x_1,...,x_n)$
    can be uniquely expressed as
    \begin{eqnarray}\label{pol}
        f(x_1,...,x_n)=\sum_{(l_1,...,l_n)\in \kappa^n}
        a_{l_1...l_n}\, x_1^{l_1}...x_n^{l_n}.
    \end{eqnarray}
\end{lem}

\begin{proof}
    It is well-known that every function over a finite field is a
    polynomial, and it is trivial that $f$ can then be expressed
    in the form (\ref{pol}). For the uniqueness, note that there are
    $\kappa^{\kappa^n}$ polynomials of that form which is exactly the number
    of $n$-ary functions over $A$.
\end{proof}

\begin{lem}\label{ros-7.2.7}
    The constant functions, the functions $h(x)=a\cdot x$, $a\in A$,
    and the operations $+$ and $-$ are affine.
\end{lem}

\begin{proof}
    This is trivial.
\end{proof}

\begin{lem}\label{ros-7.2.8}
    The functions $h(x)=x^{p^i}$, $0\leq i\leq m-1$ are affine.
\end{lem}

\begin{proof}
    We calculate
    $h(x+y)=(x+y)^{p^i}=\sum_{j=0}^{p^i}\binom{p^i}{j}\,x^j\,y^{p^i-j}$. But
    $\binom{p^i}{j}\equiv 0 \,(p)$ for $1\leq j\leq p^i-1$. Thus,
    $h(x+y)=x^{p^i}+y^{p^i}=h(x)+h(y)=h(x)+h(y)-h(0)$ and $h$ is
    affine.
\end{proof}

\begin{cor}\label{ros-7.2.9}
    The functions of the form
    \begin{eqnarray}\label{pol2}
        f(x_1,...,x_n)=a_0+\sum_{i=1}^n \sum_{j=0}^{m-1}
        a_{ij}\, x_i^{p^j}
    \end{eqnarray}
    are affine.
\end{cor}

\begin{lem}\label{ros-7.2.10}
    If $g\in\FF$ is not a function defined by (\ref{pol2}), then
    $\mathopen{<}Pol(\rho)\cup\{g\}\mathclose{>}$ contains a function
    $h(x,y)=\sum_{i,j=0}^{p^m-1} a_{ij}\,x^i\,y^j$ with at least
    one coefficient $a_{st}\neq 0$, where $1\leq s,t\leq p^m-1$.
\end{lem}

\begin{proof}
    We write $g$ as a polynomial: $g(x_1,...,x_n)=\sum_{(l_1,...,l_n)\in \kappa^n}
    a_{l_1...l_n}\, x_1^{l_1}...x_n^{l_n}$. If for one of the
    coefficients $a_{l_1...l_n}\neq0$ there are $1\leq i,j\leq n$,
    $i\neq j$ such that $l_i$ and $l_j$ are not zero, then setting
    all variables except $x_i$ and $x_j$ to $1$ yields the desired
    function. If on the other hand all non-zero coefficients have
    the form $a_{0...0 l_i0\dots 0}$, then $g(x_1,...,x_n)=g(x_1,0,...,0)+g(0,x_2,0,...,0)+...+g(0,...,0,x_n)+c$.
    Thus there is
    $1\leq q\leq n$ with the property that $g(0,...,0,x_q,0,...,0)$ has not the
    form (\ref{pol2}), for otherwise $g$ would be of that form as
    well
    which it is not. Set $f(x)=g(0,...,0,x_q,0,...,0)$ and write
    $f(x)=\sum_{i=0}^{p^m-1} b_i\,x^i$. Let $d$ be the greatest
    index in that sum such that $d$ is not a power of $p$ and
    $b_d\neq 0$; $d=d_1\cdot p^t$, where $t\geq 0$ and $d_1\geq 2$
    is not divisible by $p$. Set
    $h(x,y)=f(x+y)=\sum_{i,j=0}^{p^m-1} a_{ij}\,x^i\,y^j$. Then
    $a_{d-p^t,p^t}=\binom{d}{p^t}\,b_d$. We show
    that $\binom{d}{p^t}$ is not divisible by
    $p$:
    $$
         \binom{d}{p^t}=\binom{d_1\cdot p^t}{p^t}=\frac{d_1\,p^t\,
        (d_1\,p^t-1)...(d_1\,p^t-p)...(d_1\,p^t-2p)...(d_1\,p^t-p^t+1)}{p^t\,
        (p^t-1)...(p^t-p)...(p^t-2p)...(p^t-p^t+1)}
    $$
    One readily
    checks
    that all factors in the enumerator divisible by powers of $p$
    have corresponding factors in the denominator divisible by the
    same power of $p$. Hence, $a_{d-p^t,p^t}=\binom{d}{p^t}\, b_d$ is not
    $0$
    modulo $p$ and the lemma has been proven.
\end{proof}

\begin{lem}\label{ros-7.2.11}
    If $g\in\FF$ is not a function defined by (\ref{pol2}), then
    $\mathopen{<}Pol(\rho)\cup\{g\}\mathclose{>}$ contains the function $c(x,y)=x^s\,y^t$
    for some $1\leq s,t\leq p^m-1$.
\end{lem}

\begin{proof}
    Let $h(x,y)$ be provided by Lemma \ref{ros-7.2.10}. If all of
    the $a_{ij}$, $(i,j)\neq (s,t)$, $0\leq i,j,\leq p^m-1$, are
    $0$, we are finished by setting $c(x,y)=a_{st}^{-1}\,h(x,y)$.
    So let $a_{uv}\neq 0$ for $(u,v)\neq (s,t)$, and assume
    without loss of generality that $u\neq s$. Set
    $r(x,y)=e^u\,h(x,y)-h(e \,x,y)$, where $e$ is the primitive element
    of $(A,+,\cdot)$. Then
    $$
        r(x,y)=\sum_{i,j=0}^{p^m-1}(e^u-e^i)\,a_{ij}\,x^iy^j=\sum_{i,j=0}^{p^m-1}a'_{ij}\,x^iy^j.
    $$
    Obviously, $a'_{uv}=(e^u-e^u)\,a_{uv}=0$. Furthermore, if
    $a_{ij}=0$, then also $a'_{ij}=0$. On the other hand, as
    $a_{st}\neq 0$ and $u\neq s$, we have that
    $a'_{st}=(e^u-e^s)\,a_{st}\neq 0$. But this implies that
    iteration of this process yields a function
    $d(x,y)=d_{st}\,x^sy^t$ with $d_{st}\neq 0$. Hence,
    we can set $c(x,y)=d_{st}^{-1}\,d(x,y)$, and as all operations we used in the process
    were affine we are finished.
\end{proof}

\begin{lem}\label{ros-7.2.12}
    Let $g\in\FF$ have not the form (\ref{pol2}). Then
    $\mathopen{<}Pol(\rho)\cup\{g\}\mathclose{>}=\FF$.
\end{lem}

\begin{proof}
    We will show that the function $c(x,y)=x\cdot y$ is an element of $\mathopen{<}Pol(\rho)\cup\{g\}\mathclose{>}$, for then
    the definition of the polynomials and Lemma \ref{ros-7.2.7} imply the assertion.
    By Lemma \ref{ros-7.2.11} $\mathopen{<}Pol(\rho)\cup\{g\}\mathclose{>}$ contains
    $c(x,y)=x^sy^t$ with $1\leq s,t\leq p^m-1$. Write $s=hp^u$ and
    $t=lp^v$, where $u,v\geq 0$ and $h,l\geq 1$ are not divisible
    by $p$. By Lemma \ref{ros-7.2.8} the functions
    $a(x)=x^{p^{m-u}}$ and $b(y)=y^{p^{m-v}}$ are affine; thus,
    the function
    $w(x,y)=c(x^{p^{m-u}},y^{p^{m-v}})=(x^{p^{m-u}})^s(y^{p^{m-v}})^t=(x^{p^{m-u}})^{hp^u}(y^{p^{m-v}})^{lp^v}=(x^{p^m})^h(y^{p^m})^l=x^hy^l$
    is affine as well. Consider
    $q(x,y)=w(x+1,y+1)\in\mathopen{<}Pol(\rho)\cup\{g\}\mathclose{>}$. We write
    $q(x,y)=(x+1)^h(y+1)^l=\sum_{i,j=0}^{p^m-1}a_{ij}\,x^iy^j$. As
    $a_{11}= \binom{h}{1}\,\binom{l}{1}=h\,l$
    and $h,l$ are not divisible by $p$, we conclude that
    $a_{11}\neq 0$. Now $c(x,y)=x\cdot y\in\mathopen{<}Pol(\rho)\cup\{g\}\mathclose{>}$ is an immediate consequence
    of the proof of Lemma \ref{ros-7.2.11}.
\end{proof}

\begin{lem}\label{ros-7.2.13}
    Let $f\in \FF$. Then $f$ is affine if and only if it has the
    form (\ref{pol2}).
\end{lem}

\begin{proof}
    If $f$ has the form (\ref{pol2}), then it is affine by Lemma
    \ref{ros-7.2.9}. If conversely an $f$ existed which is affine
    but has not the form (\ref{pol2}), then $\mathopen{<}Pol(\rho)\mathclose{>}=\mathopen{<}Pol(\rho)\cup\{f\}\mathclose{>}=\FF$ by Lemma
    \ref{ros-7.2.12}, which is absurd.
\end{proof}

\begin{thm}
    If $\rho$ is a prime affine relation, then $\mathopen{<}Pol(\rho)\mathclose{>}$ is a
    maximal clone.
\end{thm}

\begin{proof}
    This is the consequence of Lemmas \ref{ros-7.2.12} and
    \ref{ros-7.2.13}.
\end{proof}

    \bibliographystyle{amsplain}
    \bibliography{cl_bibl}

\end{document}